\documentclass[10pt]{amsart}

\usepackage{amssymb,amsmath,bbm,url,xr,a4wide,amsopn}
\usepackage[mathscr]{euscript}
\usepackage{mathrsfs}
\usepackage[hidelinks]{hyperref}
\usepackage[margin=1.2in]{geometry}

\usepackage[all]{xy}
\usepackage{tikz-cd}
\tikzset{
  commutative diagrams/.cd,
  arrow style=tikz,
  diagrams={>=latex}}
\makeatletter
\tikzset{
  column sep/.code=\def\pgfmatrixcolumnsep{\pgf@matrix@xscale*(#1)},
  row sep/.code   =\def\pgfmatrixrowsep{\pgf@matrix@yscale*(#1)},
  matrix xscale/.code=%
    \pgfmathsetmacro\pgf@matrix@xscale{\pgf@matrix@xscale*(#1)},
  matrix yscale/.code=%
    \pgfmathsetmacro\pgf@matrix@yscale{\pgf@matrix@yscale*(#1)},
  matrix scale/.style={/tikz/matrix xscale={#1},/tikz/matrix yscale={#1}}}
\def\pgf@matrix@xscale{1}
\def\pgf@matrix@yscale{1}
\makeatother

\input cyracc.def
\DeclareFontFamily{U}{russian}{}
\DeclareFontShape{U}{russian}{m}{n}
        { <5><6> wncyr5
        <7><8><9> wncyr7
        <10><10.95><12><14.4><17.28><20.74><24.88> wncyr10 }{}
\DeclareSymbolFont{Russian}{U}{russian}{m}{n}
\DeclareSymbolFontAlphabet{\mathcyr}{Russian}
\makeatletter
\let\@math@cyr\mathcyr
\renewcommand{\mathcyr}[1]{\@math@cyr{\cyracc #1}}
\makeatother

\newtheorem{thm}{Theorem}[subsection]
\newtheorem{prop}[thm]{Proposition}
\newtheorem{cor}[thm]{Corollary}
\newtheorem{lm}[thm]{Lemma}
\newtheorem*{thmi}{Theorem}

\theoremstyle{definition}
\newtheorem{df}[thm]{Definition}
\newtheorem{ex}[thm]{Example}

\newtheorem{num}[thm]{}

\theoremstyle{remark}
\newtheorem{rem}[thm]{Remark}
\newtheorem{notat}[thm]{Notation}

\numberwithin{equation}{thm}

\newcommand{\changelocaltocdepth}[1]{%
  \addtocontents{toc}{\protect\setcounter{tocdepth}{#1}}%
  \setcounter{tocdepth}{#1}}

\usepackage{xspace}

\usepackage{enumitem}
\setlist[enumerate,1]{label={(\arabic*)},itemsep=\parskip} %,leftmargin=0pt
\setlist[itemize,1]{itemsep=\parskip} %,leftmargin=0pt
\newlist{thmlist}{enumerate}{2}
\setlist[thmlist,1]{label={\em(\roman*)},ref={(\roman*)},%
  itemsep=\parskip,leftmargin=*,align=left}
\setlist[thmlist,2]{label={\em(\alph*)},ref={(\alph*)},%
  itemsep=\parskip,leftmargin=*,align=left,topsep=0.1cm}
\newlist{remlist}{enumerate}{2}
\setlist[remlist,1]{label={(\roman*)},ref={(\roman*)},itemsep=\parskip,%
  leftmargin=*,align=left}
\setlist[remlist,2]{label={(\alph*)},ref={(\alph*)},itemsep=\parskip,%
  leftmargin=*,align=left,topsep=0.1cm}

\title{Fundamental classes in motivic homotopy theory}

\author{Fr\'ed\'eric D\'eglise}
\address{Institut Math\'ematique de Bourgogne\\
UMR 5584, Universit\'e de Bourgogne\\
9 avenue Alain Savary, BP 47870\\
21078 Dijon Cedex\\
France}
\email{\href{mailto:frederic.deglise@ens-lyon.fr}{frederic.deglise@ens-lyon.fr}}
\urladdr{\url{http://perso.ens-lyon.fr/frederic.deglise/}}

\author{Fangzhou Jin}
\address{Fakult\"at f\"ur Mathematik\\
Universit\"at Duisburg-Essen\\
Thea-Leymann-Strasse 9\\
45127 Essen\\
Germany}
\email{\href{mailto:fangzhou.jin@uni-due.de}{fangzhou.jin@uni-due.de}}
\urladdr{\url{https://sites.google.com/site/fangzhoujin1/}}

\author{Adeel A. Khan}
\address{Fakult\"at f\"ur Mathematik\\
Universit\"at Regensburg\\
Universit\"atsstr. 31\\
93040 Regensburg\\
Germany}
\email{\href{mailto:adeel.khan@mathematik.uni-regensburg.de}{adeel.khan@mathematik.uni-regensburg.de}}
\urladdr{\url{https://www.preschema.com}}

\date{\today}

\newcommand{\ZZ} {\mathbb Z}
\newcommand{\QQ} {\mathbb Q}
\renewcommand{\AA} {\mathbb A}
\newcommand{\GG} {\mathbb G_m}

\newcommand{\PP} {\mathbb P}

\newcommand{\T}{\mathscr T} %catégorie triangulée motivique

\newcommand{\C}{\mathscr C}

\newcommand{\cO}{\mathcal O}
\newcommand{\cE}{\mathcal E}
\newcommand{\E}{\mathbb E} %spectres
\newcommand{\F}{\mathbb F} %spectres
\newcommand{\un}{\mathbb S} %unité de SH
\newcommand{\unT}{\mathbbm 1} %unité de \T
\newcommand{\Sym}{\mathcal S}

\newcommand{\base}{\mathscr S} %schémas de base
\newcommand{\SH}{S\mathscr{H}}

\newcommand{\Hot}{\mathscr{H}_\bullet}

\newcommand{\Tri}{\mathscr{T}\mathit{ri}}

\DeclareMathOperator{\Th}{Th}

\newcommand{\fdl}{\eta}

\newcommand{\Ebiv}{\E}
\newcommand{\Ec}{\E_{\mathrm{c}}} %compact support
\newcommand{\Ebivc}{\E^{\mathrm{c}}} %compact support

\newcommand{\vb}[1]{\langle #1\rangle}%virtual vector bundles

\newcommand{\MGL}{\mathbf{MGL}}
\newcommand{\KGL}{\mathbf{KGL}}
\newcommand{\HB}{\mathbf H_{\mathcyr{B}}}

\newcommand{\HM}[1]{\mathbf H{#1}}
\DeclareMathOperator{\HMW}{\mathbf H_{MW}}
\newcommand{\BO}{\mathbf{BO}}

\DeclareMathOperator{\CH}{CH}

\DeclareMathOperator\Pic{Pic}
\DeclareMathOperator\Vect{Vect}
\DeclareMathOperator\K{K}
\DeclareMathOperator\GW{GW}
\DeclareMathOperator\wCH{\widetilde{CH}}
\newcommand{\spec}[1] {\operatorname{\mathrm{Spec}}(#1)}
\newcommand{\specx}[2] {\operatorname{\mathrm{Spec}_{#1}}\left(#2\right)}
\DeclareMathOperator{\Hom}{Hom}
\DeclareMathOperator{\Maps}{Maps}
\DeclareMathOperator{\uHom}{\underline{Hom}}

\DeclareMathOperator{\Id}{Id}
\DeclareMathOperator{\Ho}{Ho} %homotopy category

\newcommand{\pur} { \mathfrak p } %isomorphisme de pureté

\newcommand{\tr}{\mathrm{tr}}
\newcommand{\cotr}{\mathrm{cotr}}

\newcommand{\et}{\text{\'et}}

\newcommand{\op}{\mathrm{op}}

\newcommand{\term}[1]{#1\xspace}

\newcommand{\inftyCat}{\term{$\infty$-category}}
\newcommand{\inftyCats}{\term{$\infty$-categories}}

\newcommand{\inftyGrpd}{\term{$\infty$-groupoid}}
\newcommand{\inftyGrpds}{\term{$\infty$-groupoids}}

\makeatletter
\newcommand{\colim@}[2]{%
  \vtop{\m@th\ialign{##\cr
    \hfil$#1\operator@font colim$\hfil\cr
    \noalign{\nointerlineskip\kern1.5\ex@}#2\cr
    \noalign{\nointerlineskip\kern-\ex@}\cr}}%
}
\newcommand{\colim}{%
  \mathop{\mathpalette\colim@{\rightarrowfill@\textstyle}}\nmlimits@
}
\makeatother

\begin{document}

\begin{abstract}
We develop the theory of fundamental classes in the setting of motivic homotopy theory.
Using this we construct, for any motivic spectrum, an associated twisted bivariant theory,
 extending the formalism of Fulton and MacPherson.
We import the tools of Fulton's intersection theory into this setting: (refined) Gysin maps, specialization maps, and formulas for excess of intersection, self-intersections, and blow-ups.
We also develop a theory of Euler classes of vector bundles in this setting.
For the Milnor--Witt spectrum recently constructed by D\'eglise--Fasel, we get a bivariant theory extending the Chow--Witt groups of Barge--Morel, in the same way the higher Chow groups extend the classical Chow groups.
As another application we prove a motivic Gauss--Bonnet formula, computing Euler characteristics in the motivic homotopy category.
\end{abstract}

\maketitle

\tableofcontents

\parskip 0.1cm

\changelocaltocdepth{1}
%!TEX root = ../virtual.tex

\section{Introduction}

\subsection*{Historical background.}
From Poincar\'e to Grothendieck, duality has been a central component in the study of (co)homology theories.
 It led Grothendieck, building on Serre's duality, to what can nowadays be considered 
 as the summit of such a theory, the \emph{formalism of the six operations}.
  This formalism appeared in two flavours in Grothendieck's works in algebraic geometry:
   that of coherent sheaves and that of \'etale sheaves.

In the coherent setting, abstract duality was realized through
 the adjunction of the exceptional functors $(f_!,f^!)$ (see \cite[Chap.~5]{Hart}).
 The concept of dualizing complex was pivotal: it was discovered soon after that Borel--Moore homology \cite{BorelMoore} can be described as homology with coefficients in the (topological) dualizing complex.
The theory of $\ell$-adic sheaves developed in SGA4 \cite{SGA4} was the first complete incarnation of the six functors
 formalism, and for a long time the only one available in algebraic geometry.
 A key aspect of the six functor formalism that was highlighted in the seminar SGA5 \cite{SGA5} is the \emph{absolute purity} property.
 Stated in \emph{op. cit.} as a conjecture, it was partially solved by Thomason \cite{ThomasonPur}, and completely settled later by Gabber \cite{FujiwaraPurity,Gabber}.

More recently, Morel and Voevodsky introduced motivic homotopy theory \cite{MV,Voe1}.
As in algebraic topology, the stable motivic homotopy category classifies
 cohomology theories which satisfy homotopy invariance with respect to
 the algebraic affine line $\AA^1$.
The stable motivic homotopy category also satisfies
 the six functors formalism (see \cite{Ayo1}).
Moreover, it satisfies a suitable universal property \cite{Robalo} and contains the classical theories of algebraic geometry,
 such as Betti cohomology, \'etale $\ell$-adic cohomology, algebraic de Rham cohomology
 in characteristic $0$, and rigid cohomology in positive characteristic. It also incorporates newer theories such as motivic cohomology,
 algebraic K-theory and algebraic cobordism.
 The latter-mentioned theories share the common property of being
 \emph{oriented}, like their respective topological analogues, singular cohomology,
 complex K-theory, and complex cobordism.
 However, a salient feature of the motivic homotopy category is that it also contains theories which are \emph{not} oriented,
 such as Chow--Witt groups \cite{BarMor, Fasel1, Fasel2}
 and Milnor--Witt motivic cohomology \cite{DF2}, Balmer's higher Witt groups \cite{Balmer},
 hermitian K-theory (also called higher Grothendieck groups, \cite{Hornb, PaninWalter}), certain variants of algebraic cobordism \cite{PaninWalter2}, and the stable cohomotopy groups, represented by the motivic sphere spectrum.

The formalism of six operations gives rise to a great deal of structure at the level of cohomology and Borel--Moore homology groups.
Parts of this structure were axiomatized by Bloch and Ogus \cite{BO}, via their notion of \emph{Poincar\'e duality theory}, and later through the \emph{bivariant theories} of Fulton and MacPherson \cite{FMP}.
The key element of these axiomatizations was the notion of the \emph{fundamental class}, which was used to express duality isomorphisms.

\vspace{-0.2cm}
\subsection*{Main problematic.}
 Our goal in this paper is to incorporate Fulton
 and MacPherson's ideas into stable motivic homotopy theory, thereby obtaining a universal
 bivariant theory. In order to treat oriented and non-oriented spectra
 in a single theory, we have to replace Tate twists, as used for example in the Bloch--Ogus axiomatic, by ``Thom twists'',
 i.e., twists with respect to vector bundles (or more generally, with respect to virtual
 vector bundles). Let us explain the justification for this idea.

Our first inspiration is Morel and Voevodsky's homotopy purity theorem, which asserts that, for smooth closed pairs
 $(X,Z)$, the homotopy type of $X$
 with support in $Z$ is isomorphic to $\Th_Z(N_ZX)$, 
 the Thom space of the normal bundle of $Z$ in $X$.
 Here the homotopy type of $\Th_Z(N_ZX)$ should be understood as the homotopy type of $Z$
 twisted by the vector bundle $N_ZX$.
Another motivation is Morel's work \cite{MorelLNM} on computations of homotopy groups, in which a crucial role is played by the construction of good transfer maps
 for finite field extensions in the unstable homotopy category. In this work, twists are usually
 avoided but at the cost of choosing orientations. Similar constructions enter into play in
 Voevodsky's theory of framed correspondences \cite{VoevodskyFramed,GP,EHKSY}, where the ``framing'' provides a chosen 
 trivialization of the normal bundle. Finally,  Calm\`es and Fasel have introduced the notion
 of MW-correspondences, based on Chow--Witt theory, where transfers do appear
 with a twist.
 These examples show the utmost importance of having good transfer or Gysin morphisms in $\AA^1$-homotopy theory.
 The last indication which points out to our central construction
 is the extension obtained in \cite{Jin}
 of the finer operations of Fulton's intersection theory, such as refined Gysin morphisms,
 in the motivic Borel-Moore theory (see in particular \cite[Def.~3.1]{Jin}).
 The translation becomes possible once one recognizes Borel-Moore motivic homology as a particular
 instance of bivariant theory.

%As the reader can guess, 
Another fundamental model for our theory is that of Chow--Witt groups.
 In this theory, the necessity of considering twists appears most notably when Gysin morphisms are
 at stake (see \cite{Fasel1, Fasel2}).
 Much of the interest in these groups comes from the fact that they are natural receptacles for \emph{Euler classes} of vector bundles.
The Euler class provides an obstruction for a vector bundle to split off a trivial summand of rank one (see \cite{BargeMorel,MorelLNM,FaselSrinivas}).
In this paper we also develop a general theory of Euler classes in $\mathbb{A}^1$-homotopy theory, and show that they satisfy the expected obstruction-theoretic property (Corollary~\ref{cor:Euler obstruction}). Our motivation to introduce these Euler classes is to formulate excess intersection formulas in our bivariant theories, see the following theorem.

%: the Euler class of the excess bundle measures the failure of the base change formula for a non-transverse intersection.
%This formula specializes as usual to self-intersection and blow-up formulas.
%In Borel--Moore motivic homology, for instance, the Euler class coincides with the top Chern class, and we recover the usual excess intersection formula.
%In the example of Milnor--Witt bivariant theory, this formula is new, even for just the Chow--Witt groups.

\subsection*{Main construction.}
The Thom space functor $\Th_X$,
 associating to a vector bundle $E$ over a scheme $X$ its (stable) Thom space $\Th_X(E) \in \SH(X)$, canonically extends to the Picard groupoid of virtual vector bundles over $X$ (see \cite[4.1]{Riou}).
Given any motivic spectrum $\E \in \SH(S)$, we can define the \emph{(twisted) bivariant theory} with coefficients in $\E$, graded by integers $n \in \ZZ$ and virtual vector bundles $v$ on $X$, as the following group:
\begin{align*}
\Ebiv_{n}(X/S,v) := \Hom_{\SH(S)}(\Th_X(v)[n], p^!(\E)),
\end{align*}
for any morphism $p : X \to S$ that is separated of finite type. For $\E=\un_S$, we simply write $H_n(X/S,v)$ and call this \emph{bivariant $\AA^1$-theory}. %Note that the above construction
 %is functorial in $\E$. 
The construction is functorial in $\E$ so that, given any ring spectrum $\E$ over $S$ with
 unit $\nu:\un_S \rightarrow \E$, we get a canonical map:
\begin{align*}
H_*(X/S,v) \rightarrow \Ebiv_*(X/S,v),
\end{align*}
expressing the universal role of bivariant $\AA^1$-theory.\footnote{We call this map
 the $\AA^1$-regulator map, by extension of Beilinson's terminology; see Definition \ref{df:regulator}.}
These bivariant theory groups satisfy a rich functoriality, including covariance for proper maps and contravariance for \'etale
 maps; see Paragraph~\ref{num:bivar_ppties} for details.
 There is also a composition product, which takes the form
\begin{align*}
\Ebiv_*(Y/X,w) \otimes \Ebiv_*(X/S,v) \rightarrow \Ebiv_*(Y/S,w+q^*v), (y,x) \mapsto y.x,
\end{align*}
for schemes $Y/X/S$, and virtual vector bundles $v/X$, $w/Y$ (see again Paragraph~\ref{num:bivar_ppties}). That being
 given, here is the central construction of this paper.
\begin{thmi}[Theorem \ref{thm:fdl_qplci} and Proposition \ref{prop:lci_ex_inters}]
For any \emph{smoothable lci}\footnote{A morphism of schemes is smoothable if it admits a (global) factorization 
 into a closed immersion followed by a smooth morphism; such a morphism is lci (a local complete intersection) if the closed immersion is a regular immersion; see our conventions.}
 morphism $f : X \to Y$, there exists a canonical class $\fdl_f$,
 called the \emph{fundamental class} of $f$:
\begin{align*}
    \fdl_f \in H_0(X/Y,\vb{L_f}),
\end{align*}
where $\vb{L_f}$ is the virtual tangent bundle of $f$
 (equivalently, the virtual bundle associated with the cotangent complex of $f$,
 which is perfect under our assumption).

These classes satisfy the following properties:
\begin{thmlist}
\item \emph{Associativity}. Consider morphisms  $f : X \to Y$ and $g : Y \to Z$
 such that $f$, $g$, and $g \circ f$ are smoothable and lci. Then one has:
\begin{align*}
\eta_g.\eta_f \simeq \eta_{g \circ f}
\end{align*}
in $H_0(X/Z, \vb{L_f}+f^*\vb{L_g}) \simeq H_0(X/Z,\vb{L_{g\circ f}})$.
\item \emph{Excess intersection}. Consider a cartesian square
\begin{equation}
    \begin{tikzcd}[matrix scale=0.7]
      Y \ar{d}{v}\ar{r}{g}\ar[phantom]{rd}{\scriptstyle\Delta}
        & T\ar{d}{u}
      \\
      X \ar[swap]{r}{f}
        & S
    \end{tikzcd}
\end{equation}
such that $f$ and $g$ are smoothable lci. Consider a factorization
 $X \xrightarrow i P \xrightarrow p S$ of $f$ such that $p$ is smooth and $i$ is
 a closed immersion. Let $k$ be the pullback of $i$ along $u$ and denote by $\xi$
 the quotient bundle of the natural monomorphism of vector bundles:
 $N_k \rightarrow u^{-1}N_i$. Then there exists an Euler class
 $e(\xi) \in  H_0(Y/Y,\vb \xi)$ such that the following formula holds:
\begin{align*}
\Delta^*(\fdl_f) \simeq e(\xi).\fdl_g
\end{align*}
in $H_0(Y/T,v^*\vb{L_f}) \simeq H_0(Y/T, \vb{L_g}-\vb{\xi})$.

In particular, if the square $\Delta$ is tor-independent, we get
 $\Delta^*(\fdl_f) \simeq \fdl_g$.
\end{thmlist}
\end{thmi}
This construction is universal in the stable motivic homotopy category. Indeed, given a motivic ring spectrum $\E$
 the regulator map (Definition~\ref{df:regulator}) gives rise to fundamental classes with coefficients in $\E$.
 This in turn yields Gysin homomorphisms
\begin{align*}
    f^! : \E_n(Y/S, e) \to \E_n(X/S, f^*(e)+\vb{L_f}), \quad y \mapsto \fdl_f.y,
\end{align*}
 for a smoothable lci morphism $f : X \to Y$,
and the associativity (resp. excess intersection) property above corresponds
 to the compatibility with composition (resp. excess intersection formula) satisfied
 by these Gysin morphisms, as in Chow theory. Note that for oriented theories like Borel--Moore motivic homology, the orientation provides a Thom isomorphism that replaces the virtual twist $\vb{L_f}$ with a shift by the relative virtual dimension $\chi(L_f)$, so that the Gysin homomorphism takes a more familiar shape.\footnote{Similar simplifications of twists occur for the so-called \emph{Sp-oriented theories} of Panin and Walter. We leave the general formulation for future works. The reader can also consult Example \ref{ex:MW cohomology}.}

From a categorical point of view, the universality of our construction is best stated in
 the language of the six operations: the fundamental class $\fdl_f$ as above corresponds
 to a natural transformation of functors
\begin{align*}
    \pur_f : f^*(-) \otimes \Th_X(L_f) \to f^!
\end{align*}
that we call the \emph{purity transformation} associated to $f$
 (see Section~\ref{sec:applications/purity}). Then by adjunction, we get 
 \emph{trace} and \emph{cotrace} maps, extending the classical construction of SGA4:
\begin{align*}
    \tr_f&:f_!\Sigma^{L_f}f^* \rightarrow \Id \\
    \cotr_f&:\Id \rightarrow f_*\Sigma^{-L_f}f^!.
\end{align*}
These natural transformations can be considered as a natural extension,
 and in fact an important part, of the six functors formalism.\footnotemark
 \footnotetext{In fact, we show that our construction allows one
 to define these transformations in a greater generality, say for arbitrary
 motivic $(\infty,1)$-category of coefficients in the sense of \cite[Chap.~2, Def.~3.5.2]{KhanThesis}.
 See Paragraph \ref{num:motivic_categories}.} The Gysin morphisms in bivariant theory above are immediate consequences of these transformations, when
 applied to the ring spectrum $\E$. Moreover, one gets Gysin morphisms (wrong-way variance)
 for the traditional package of co/homological theories, with and without proper support,
 associated with a spectrum
 (even without a ring structure).
 We refer the interested reader to Paragraph \ref{num:Gysin four} for details,
 and to Section \ref{sec:applications/examples} for a list of concrete examples.

When $f$ is \emph{smooth}, the map $\pur_f$ considered above is invertible and 
 coincides with the classical purity isomorphism of the six functors,
 as defined by Ayoub (see Paragraph~\ref{num:smooth purity}).
In general, the purity transformation measures the failure of a given motivic spectrum $\E \in \SH(Y)$ to satisfy the \emph{purity} property with respect to $f$ (see Definition~\ref{df:f-pure}); when $f$ is a closed immersion between regular schemes, this property corresponds to the notion of \emph{absolute purity} (see Definition~\ref{df:absolute_purity}), axiomatizing the original
 conjecture of Grothendieck.  Such an axiomatization is not new (see \cite[1.3]{Deg12},
 \cite[A.2]{CD4}). However, the formulation we obtain here is more flexible and has a number of advantages (see Example \ref{rem:absolute_purity}).
 The absolute purity is essential for arithmetic applications, and has already been obtained
 in several contexts (rational motives, \'etale motives, $\KGL$-modules). We believe that new examples
 will be obtained in the future (and have work in progress in that direction).
 Finally, our construction has been applied in \cite[Appendix~A]{FranklandSpitzweck} to prove a new absolute purity result for motivic cohomology.

\subsection*{Further applications.}
We finally give several applications of our formalism. The first one is the possibility
 of extending Fulton's theory of \emph{refined Gysin morphism} to an arbitrary ring spectrum
 (Definition \ref{df:ref_fdl}). These new refined Gysin morphisms are used to
 define \emph{specialization maps} in any representable theory, on the model of Fulton's definition
 of specialization for the Chow group. In fact, our specialization maps can be lifted to
 natural transformations of functors (see Paragraph \ref{num:specialization 2} for details).
 Most interestingly, the theory can be applied to Chow--Witt groups and give specializations
 of quadratic cycles (see Example \ref{ex:specialization 1 HZtilde}).

Note that the idea of refining classical formulas to the quadratic setting has been explored recently by many authors \cite{Fasel1,FaselExcess,HoyoisLefschetz,KassWickelgren,LevineEnumerative}.
In this direction another application of the theory we develop is a motivic refinement of the classical Gauss-Bonnet formula \cite[VII~4.9]{SGA5}.
Given a smooth proper $S$-scheme $X$, the \emph{categorical Euler characteristic} $\chi^{cat}(X/S)$ is the endomorphism of the motivic sphere spectrum $\un_S$ given by the trace of the identity map of $\Sigma^\infty_+(X)\in\SH(S)$.
%there is a so-called  %, which lives in the bivariant theory represented by the .
A simple application of our excess intersection formula then computes this invariant as the degree of the Euler class of the tangent bundle $T_{X/S}$ (see Theorem~\ref{thm:Gauss-Bonnet}).
This result is a generalization of a theorem of Levine \cite[Thm. 1]{LevineEnumerative}, which applies when $S$ is the spectrum of a field and $p : X \to S$ is smooth projective. It also recovers the $\mathrm{SL}$-oriented variant in \cite[Theorem~1.5]{LevineRaksit}, used in \emph{op.~cit.} to prove a certain explicit formula for the quadratic Euler characteristic conjectured by Serre.

\subsection*{Related work and further developments.}
Bivariant theories represented by \emph{oriented} motivic spectra were studied in detail in \cite{Deg16}. In that setting, the fundamental class of a regular closed immersion
 is given by a construction of Navarro \cite{Navarro}, which itself is based on a construction of Gabber in the setting of \'etale cohomology \cite[Exp.~XVI]{Gabber}.
 A similar construction to our fundamental class for closed immersions,
 in the context of equivariant stable $\AA^1$-homotopy,
 has been developed independently in recent work of Levine \cite{LevineVirtual}.

An immediate but important consequence of our work here is that the cohomology theory represented by any motivic spectrum $\E$ admits a canonical structure of \emph{framed transfers}; that is, it extends to a presheaf on the category of framed correspondences (see \cite{VoevodskyFramed,EHKSY}).
It is proven in \cite{EHKSY} that this structure can be used to recognize infinite $\PP^1$-loop spaces, in the same way that $\mathcal{E}_\infty$-structures can be used to recognize infinite loop spaces in topology (see also \cite{GP}).
In \cite{EHKSY2,DruzhininKolderup}, framed transfers are applied to construct categories of finite $\E$-correspondences, for any motivic spectrum $\E$, together with canonical functors from the category of framed correspondences.
As explained in \emph{op.~cit.}, such functors play an important role in the yoga of motivic categories.
Another application of the existence of framed transfers is a topological invariance statement for the motivic homotopy category, up to inverting the exponential characteristic of the base field \cite{ElmantoKhanPerfection}.

An application of our theory of Euler classes can be found in \cite{Jin18}, where it is used to give a characterization of the characteristic class of a motive.

Finally, our constructions can be extended to the setting of \emph{quasi-smooth} morphisms in derived algebraic geometry.
This yields a formalism of motivic virtual fundamental classes, see \cite{KhanVirtual}.

\subsection*{Contents.}

In Section~\ref{sec:bivariant} we construct the bivariant theory and cohomology theory associated to a motivic ring spectrum, and study their basic properties.
Following Fulton--MacPherson \cite{FMP}, we also introduce the abstract notion of \emph{orientations} of morphisms in this setting; fundamental classes will be examples of orientations.
We then show how any choice of orientation gives rise to a purity transformation.

The heart of the paper is Section~\ref{sec:construction}, where we construct fundamental classes and verify their basic properties.
In the case where $f$ is smooth, the fundamental class comes from the purity theorem (see Definition~\ref{df:fdl_smooth}).
For the case of a regular closed immersion we use the technique of deformation to the normal cone.
We then explain how to glue these to obtain a fundamental class for any quasi-projective lci morphism.
Throughout this section, we restrict our attention to the bivariant theory represented by the motivic sphere spectrum.

Finally in Section~\ref{sec:applications}, we return to the setting of the bivariant theory represented by any motivic ring spectrum.
We show how the fundamental class gives rise to Gysin homomorphisms.
We prove the excess intersection formula in this setting (Proposition~\ref{prop:excess with coefficients}).
We also discuss the purity transformation, the absolute purity property, and duality isomorphisms (identifying bivariant groups with certain cohomology groups).
We then import some further constructions from Fulton's intersection theory, including refined Gysin maps and specialization maps.
Finally, we conclude with a proof of the motivic Gauss-Bonnet formula mentioned above.

\subsection*{Conventions.}

The following conventions are in place throughout the paper:
\begin{enumerate}
  \item All schemes in this paper are assumed to be quasi-compact and quasi-separated.
  \item The term \emph{s-morphism} is an abbreviation for ``separated morphism of finite type''.
  Similarly, an \emph{s-scheme over $S$} is an $S$-scheme whose structural morphism is an s-morphism.
  \item We write $\AA^1$ for the affine line over $\spec{\ZZ}$ and $\GG$ for the complement of the origin.
  For a scheme $X$, we write $\AA^1 X$ and $\GG X$ for $\AA^1\times{X}$ and $\GG\times{X}$, respectively.
  \item We follow \cite[Exps. VII--VIII]{SGA6} for our conventions on regular closed immersions and lci morphisms.
  Recall that if $X$ and $Z$ are regular schemes, or are both smooth over some base $S$, then any closed immersion $Z \to X$ is regular.
  Given a regular closed immersion $i:Z \to X$, we write $N_i$, $N_ZX$, or occasionally $N(X,Z)$, for its normal bundle.
  Recall that a morphism of schemes $X \to S$ is \emph{lci} (= a local complete intersection) if it admits, Zariski-locally on the source, a factorization $X \xrightarrow{i} Y \xrightarrow{p} S$, where $p$ is smooth and $i$ is a regular closed immersion.
  If it is also \emph{smoothable}, then it admits such a factorization globally.
  This is for example the case if $f$ is quasi-projective (in the sense that it factors through an immersion into some projective space $\PP^n_S$).
  \item
  A cartesian square of schemes
   \begin{equation}
     \begin{tikzcd}[matrix scale=0.7]
       X' \ar{r}\ar{d}
         & Y' \ar{d}{p}
       \\
       X \ar{r}{f}
         & Y
     \end{tikzcd}
   \end{equation}  
  % $$
  % \xymatrix@=10pt{
  % X'\ar^-{}[r]\ar_-{}[d] & Y'\ar^-{p}[d] \\
  % X\ar^-{f}[r] & Y
  % }
  % $$
  is \emph{tor-independent} if the groups $\mathrm{Tor}^{\cO_Y}_i(\cO_X, \cO_{Y'})$ vanish for $i>0$.
  In this case we also say that $p$ is \emph{transverse} to $f$.
  Recall that if $p$ or $f$ is flat, then this condition is automatic.

  \item
  We will make use of the language of stable \inftyCats \cite{HA}.
  Given a stable \inftyCat $\C$, we write $\Maps_\C(X, Y)$ for the mapping spectrum of any two objects $X$ and $Y$.
  We write $\Hom_\C(X,Y)$ or simply $[X, Y]$ for the abelian group of connected components $\pi_0 \Maps_\C(X,Y)$.

  \item
  Given a topological $S^1$-spectrum $E$, we write $x\in E$ to mean that $x$ is a point in the infinite loop space $\Omega^\infty(X)$ (or an object in the corresponding \inftyGrpd).
\end{enumerate}

\subsection*{Acknowledgments.}
We would like to thank Denis-Charles Cisinski for many helpful discussions.
In particular, we learned the interpretation of the Gysin map in terms of the six operations from him.
We also thank Martin Frankland for comments on a previous revision.
F. Jin learned from Marc Levine his motivic Gauss-Bonnet formula during the lectures at the the 2017 Trimester Program \emph{$K$-theory and related fields} at the Hausdorff Research Institute for Mathematics in Bonn.
The third author would like to thank Elden Elmanto, Marc Hoyois, Vova Sosnilo and Maria Yakerson for their collaboration on the papers \cite{EHKSY,EHKSY2}, which served as a strong motivation in writing this paper.
We would also like to thank the referees for their careful reading and attentive comments which helped clarify some missing points.
Finally, we thank the FRIAS for its hospitality while part of this work was completed.

F. D\'eglise received support from the French ``Investissements d'Avenir'' program, project ISITE-BFC (contract ANR-lS-IDEX-OOOB) and from a Marie Sk\l{}odowska-Curie FCFP Senior Fellowship hosted by the FRIAS. F. Jin is partially supported by the DFG Priority Programme SPP 1786 Project ``Motivic filtrations over Dedekind domains''. 

\changelocaltocdepth{2}

%!TEX root = ../virtual.tex

\section{Bivariant theories and cohomology theories}
\label{sec:bivariant}

%!TEX root = ../../virtual.tex

\subsection{The six operations}
\label{ssec:six}

Given any (quasi-compact quasi-separated) scheme $S$, we write $\SH(S)$ for the stable \inftyCat of motivic spectra, as in \cite[Appendix~C]{HoyoisLefschetz} or \cite{KhanThesis}.
When $S$ is noetherian and of finite dimension, then the homotopy category of $\SH(S)$ is equivalent, as a triangulated category, to the stable $\AA^1$-homotopy category originally constructed by Voevodsky \cite{Voe1}.
As $S$ varies, these categories are equipped with the formalism of Grothendieck's six operations \cite{Ayo1,CD3}.
In this subsection we briefly recall this formalism, and its $\infty$-categorical refinement as constructed in \cite[Chap.~2]{KhanThesis} (see also \cite{RobaloThesis} for another approach).

\begin{num}
First, the stable presentable \inftyCat $\SH(S)$ is symmetric monoidal, and we denote the monoidal product and monoidal unit by $\otimes$ and $\un_S$, respectively.
It also admits internal hom objects $\uHom(\E, \F) \in \SH(S)$ for all $\E, \F \in \SH(S)$.
For any morphism of schemes $f : T \to S$, we have a pair of adjoint functors
  \begin{equation*}
    f^* : \SH(S) \to \SH(T),
    \qquad
    f_* : \SH(T) \to \SH(S),
  \end{equation*}
called the functors of \emph{inverse} and \emph{direct image} along $f$, respectively.
If $f$ is an \emph{s-morphism}\footnote{Using Zariski descent, the operations $(f_!,f^!)$ can be extended to the case where $f$ is locally of finite type; assuming this extension, the reader can globally redefine the term ``s-morphism'' as ``locally of finite type morphism''.}, i.e., a separated morphism of finite type, then there is another pair of adjoint functors
  \begin{equation*}
    f_! : \SH(T) \to \SH(S),
    \qquad
    f^! : \SH(S) \to \SH(T),
  \end{equation*}
called the functors of \emph{exceptional} direct and inverse image along $f$, respectively.
Each of these operations is $2$-functorial.
\end{num}

\begin{num}\label{num:six operation compatibilities}
The six operations $(\otimes, \uHom, f^*, f_*, f_!, f^!)$ satisfy a variety of compatibilities.
These include:
\begin{enumerate}
  \item\label{item:six operation compatibilities/f^* monoidal}
  For every morphism $f$, the functor $f^*$ is symmetric monoidal.

  \item\label{item:six operation compatibilities/proper}
  There is a natural transformation $f_! \to f_*$ which is invertible when $f$ is proper.

  \item\label{item:six operation compatibilities/open immersion}
  There is an invertible natural transformation $f^* \to f^!$ when $f$ is an open immersion.

  \item\label{item:six operation compatibilities/projection formula}
  The operation of exceptional direct image $f_!$ satisfies a projection formula against inverse image.
  That is, there is a canonical isomorphism
    \begin{equation*}
      \E \otimes f_!(\F) \to f_!(f^*(\E) \otimes \F)
    \end{equation*}
  for any s-morphism $f : T \to S$ and any $\E \in \SH(S)$, $\F \in \SH(T)$.

  \item\label{item:six operation compatibilities/base change}
  The operation $f_!$ satisfies base change against inverse images $g^*$, and similarly $f^!$ satisfies base change against direct images $g_*$.
  That is, for any cartesian square
    \begin{equation}\label{eq:cart}
      \begin{tikzcd}[matrix scale=0.7]
        T' \ar{r}{g}\ar[swap]{d}{q}
          & S' \ar{d}{p}
        \\
        T \ar[swap]{r}{f}
          & S
      \end{tikzcd}
    \end{equation}
  where $f$ and $g$ are s-morphisms, there are canonical isomorphisms
    \begin{equation*}
      p^*f_! \to g_!q^*,
      \qquad
      q_*g^! \to f^!p_*.
    \end{equation*}
\end{enumerate}
All the above data are subject to a homotopy coherent system of compatibilities (see \cite[Chap.~2, Sect.~5]{KhanThesis}).
\end{num}

\begin{num}\label{num:six/A^1-invariance}
The $\AA^1$-homotopy invariance property of $\SH$ is encoded in terms of the six operations as follows.
For a scheme $S$ and any vector bundle $\pi : E \to S$, the functor $\pi^*:\SH(S) \to \SH(E)$ is \emph{fully faithful}.
In particular, the unit $\Id \to \pi_*\pi^*$ is invertible.
\end{num}

\begin{num}\label{num:Thom suspension}
Given a locally free sheaf $\cE$ of finite rank over $S$, let $E = \specx{S}{\Sym(\cE)}$ denote the associated vector bundle\footnote{Throughout the text we will generally not distinguish between a locally free sheaf $\cE$ and its associated vector bundle $E$.  Thus for example we will also write $\Sigma^{E}$ instead of $\Sigma^{\cE}$, or similarly $\vb{E} \in \K(S)$ instead of $\vb{\cE} \in \K(S)$ (see Paragraph~\ref{num:Thom_virtual} below).  It should always be clear from the context what is intended.}.
There is an auto-equivalence
  \begin{equation*}
    \Sigma^\cE : \SH(S) \to \SH(S),
  \end{equation*}
called the \emph{$\cE$-suspension functor}, with inverse denoted $\Sigma^{-\cE}$.
This functor is compatible with the monoidal product $\otimes$ via a projection formula that provides canonical isomorphisms $\Sigma^\cE(\E) \simeq \E \otimes \Sigma^\cE(\un_S)$ for any $\E \in \SH(S)$.
It is also compatible with the other operations in the sense that we have canonical isomorphisms
  \begin{equation}\label{eq:Thom suspension compatibilities}
    f^*\Sigma^{\cE} \simeq \Sigma^{f^*(\cE)} f^*,
    \quad
    f_*\Sigma^{f^*(\cE)} \simeq \Sigma^{\cE} f_*,
    \quad
    f_!\Sigma^{f^*(\cE)} \simeq \Sigma^{\cE} f_!,
    \quad
    f^!\Sigma^{\cE} \simeq \Sigma^{f^*(\cE)} f^!.
  \end{equation}

The motivic spectrum $\Sigma^\cE(\un_S)\in \SH(S)$ is (the suspension spectrum of) the \emph{Thom space} of $\cE$, and is denoted $\Th_S(\cE)$.
Its $\otimes$-inverse $\Sigma^{-\cE}(\un_S)$ is denoted $\Th_S(-\cE)$.

Given a motivic spectrum $\E \in \SH(S)$, we denote by $\E(n) \in \SH(S)$ the motivic spectrum $\Sigma^{\cO^n_S}(\E)[-2n]$ for each $n \ge 0$.
The assignment $\E \mapsto \E(n)$ defines then another auto-equivalence of $\SH(S)$, inverse to $\E \mapsto \E(-n) = \Sigma^{-\cO^n_S}(\E)[2n]$.
\end{num}

\begin{num}\label{num:Thom_virtual}
Let $\Vect(S)$ denote the groupoid of locally free sheaves on $S$ of finite rank, and $\Pic(\SH(S))$ the \inftyGrpd of $\otimes$-invertible objects in $\SH(S)$.
The assignment $\cE \mapsto \Th_S(\cE)$ determines a map of presheaves of \inftyGrpds
  \begin{equation*}
    \Th : \Vect \to \Pic(\SH).
  \end{equation*}
Moreover, if $\K$ denotes the presheaf sending $S$ to its Thomason--Trobaugh K-theory space $\K(S)$, this extends to a map of $\cE_\infty$-groups
  \begin{equation*}
    \Th : \K \to \Pic(\SH),
  \end{equation*}
see \cite[Subsect.~16.2]{BachmannHoyois}.
In particular, any perfect complex $\cE$ on $S$ defines a K-theory class\footnotemark~$\vb{\cE} \in \K(S)$ and thus an auto-equivalence $\Sigma^{\cE} : \SH(S) \to \SH(S)$ and a Thom space $\Th_S(\cE)\in\SH(S)$.
The formulas \eqref{eq:Thom suspension compatibilities} also extend.
Moreover, any exact triangle $\cE' \to \cE \to \cE''$ of perfect complexes induces canonically a path $\vb{\cE} \simeq \vb{\cE'} + \vb{\cE''}$ in the space $\K(S)$, hence also identifications
  \begin{equation}\label{eq:Thom exact sequence}
    \Sigma^{\cE} \simeq \Sigma^{\cE'}\Sigma^{\cE''}
      \simeq \Sigma^{\cE''}\Sigma^{\cE'}
  \end{equation}
and an isomorphism $\Th_S(\cE) \simeq \Th_S(\cE') \otimes \Th_S(\cE'')$ in $\SH(S)$. 
\footnotetext{An abuse of language we will commit often is to say ``K-theory class'' when it would be more precise to say ``point of the K-theory space''.}
\end{num}

\begin{num}\label{num:localization}
If $i : Z \to S$ is a closed immersion, then the direct image functor $i_* : \SH(Z) \to \SH(S)$ is fully faithful.
Moreover, if the complementary open immersion $j : U \to S$ is quasi-compact, then by the Morel--Voevodsky localization theorem there is an exact triangle
  \begin{equation*}
    i_*i^! \to \Id \to j_*j^*.
  \end{equation*}
\end{num}

\begin{num}\label{num:smooth purity}
If $f$ is a smooth s-morphism, then by Ayoub's purity theorem, there is a canonical isomorphism of functors
  \begin{equation*}
    \pur_f : \Sigma^{T_f} f^* \to f^!,
  \end{equation*}
where $T_{f}$ is the relative tangent bundle.
It follows in particular that $f^*$ admits a left adjoint
  \begin{equation*}
    f_\sharp = f_! \Sigma^{T_{f}}
  \end{equation*}
which satisfies base change and projection formulas against inverse images $g^*$.
If $f$ is \'etale then we get an isomorphism $\pur_f : f^* \simeq f^!$, generalizing Paragraph~\ref{num:six operation compatibilities}\ref{item:six operation compatibilities/open immersion}.
\end{num}

\begin{num}\label{num:relative purity}
Let $f : X \to Y$ be a closed immersion of s-schemes over $S$.
Suppose that $X$ and $Y$ are \emph{smooth} over $S$, with structural morphisms $p : X \to S$ and $q : Y \to S$.
Then by the relative purity theorem of Morel--Voevodsky, there exist isomorphisms of functors
  \begin{equation*}\label{eq:relative purity}
    q_\sharp f_* \simeq p_\sharp \Sigma^{N_XY},
    \qquad
    \Sigma^{-N_XY} p^* \simeq f^! q^*,
  \end{equation*}
where $N_XY$ denotes the normal bundle.
\end{num}

Many further compatibilities can be derived from the ones already listed.
A few that will be especially useful in this paper are as follows:

\begin{num}\label{num:Ex^*!}
Given a cartesian square as in \eqref{eq:cart} where $f$ and $g$ are s-morphisms, the base change formula (Paragraph~\ref{num:six operation compatibilities}\ref{item:six operation compatibilities/base change}) induces a natural transformation
  \begin{equation}\label{eq:Ex^*!}
    Ex^{*!} : q^*f^! \to g^!p^*.
  \end{equation}
It follows from the purity theorem (Paragraph~\ref{num:smooth purity}) that if $f$ or $p$ is smooth, then $Ex^{*!}$ is invertible.
For example if $i : Z \to S$ is a closed immersion, then the cartesian square
  \begin{equation*}
    \begin{tikzcd}[matrix scale=0.7]
      Z \ar[equals]{r}\ar[equals]{d}
        & Z \ar{d}{i}
      \\
      Z \ar{r}{i}
        & S
    \end{tikzcd}
  \end{equation*}
gives rise to a canonical natural transformation
  \begin{equation}\label{eq:i^! -> i^*}
    i^! \to i^*.
  \end{equation}
\end{num}

\begin{num}\label{num:Ex^!*_otimes}
For any s-morphism $f : X \to S$ and any pair of motivic spectra $\E,\F\in\SH(S)$, there is a canonical morphism
  \begin{equation*}
    Ex^{!*}_\otimes : f^!(\E) \otimes f^*(\F) \to f^!(\E \otimes \F)
  \end{equation*}
induced by adjunction from the projection formula (Paragraph~\ref{num:six operation compatibilities}\ref{item:six operation compatibilities/projection formula}).
% For a fixed $\E$ (resp. a fixed $\F$), these define a natural transformation.
If $\F$ is \emph{$\otimes$-invertible}, then $Ex^{!*}_\otimes$ is invertible.
\end{num}
%!TEX root = ../../virtual.tex

\subsection{Bivariant theories}
\label{sec:bivariant/bivariant}

In this subsection we construct the \emph{bivariant theory} represented by a motivic spectrum, and state its main properties.
In fact, bivariant theory is only one of the ``four theories'' associated to a motivic spectrum (cf. \cite[Chap.~4, Sect.~9]{FSV}).
For sake of completeness we define them all now:

\begin{df}\label{df:bivariant}
Let $S$ be a scheme and $\E \in \SH(S)$ a motivic spectrum.

\begin{remlist}
\item\emph{Bivariant theory.}
For any s-morphism $p : X \to S$ and any K-theory class $v \in \K(X)$, we define the \emph{$v$-twisted bivariant spectrum} of $X$ over $S$ as the mapping spectrum
  \begin{align*}
    % \Ebiv(X/S, v) = \Maps_{\SH(S)}(p_!(\Th_X(v)),\E).
    \Ebiv(X/S, v) &= \Maps_{\SH(S)}(\un_S, p_*(p^!(\E) \otimes \Th_X(-v)))\\
      &\simeq \Maps_{\SH(S)}(p_!(\Th_X(v)),\E).
  \end{align*}
We also write
  \begin{equation*}
    \Ebiv_n(X/S, v) = \pi_{n} \Ebiv(X/S,v) %= \left[p_!(\Th_X(v))[n],\E\right]
      = \left[\un_S[n], p_*(p^!(\E) \otimes \Th_X(-v))\right]
  \end{equation*}
for each integer $n\in\ZZ$.

\item\emph{Cohomology theory.}
For any morphism $p : X \to S$ and any $v \in \K(X)$, we define the \emph{$v$-twisted cohomology spectrum} of $X$ over $S$ as the mapping spectrum
  \begin{align*}
    \E(X, v)
      &= \Maps_{\SH(S)}(\un_S, p_*(p^*(\E) \otimes \Th_X(v)))\\
      &\simeq \Maps_{\SH(X)}(\un_X, p^*\E \otimes \Th_X(v)).
  \end{align*}
We also write
  \begin{equation*}
  \E^n(X, v) = \pi_{-n} \E(X,v) = \left[\un_S, p_*(p^*(\E) \otimes \Th_X(v))[n]\right].
  \end{equation*}
for each integer $n\in\ZZ$.

\item\emph{Bivariant theory with proper support (or homology).}
For any s-morphism $p : X \to S$ and any K-theory class $v \in \K(X)$, we define the spectrum of \emph{$v$-twisted bivariant theory with proper support} of $X$ over $S$ as the mapping spectrum
  \begin{equation*}
    \E^c(X/S,v)=\Maps_{\SH(S)}\big(\un_S,f_!(f^!(\E_S) \otimes \Th_X(-v))\big)
  \end{equation*}

\item\emph{Cohomology with proper support.}
For any s-morphism $p : X \to S$ and any K-theory class $v \in \K(X)$, we define the spectrum of \emph{$v$-twisted cohomology with proper support} of $X$ over $S$ as the mapping spectrum
  \begin{equation*}
    \E_c(X/S,v)=\Maps_{\SH(S)}\big(\un_S,f_!(f^*(\E) \otimes \Th_X(v))\big).
  \end{equation*}
\end{remlist}
\end{df}

\begin{rem}\label{rem:cohomology in terms of bivariant theory}
Note that we have canonical identifications
  \begin{equation*}
    \E(X, v) \simeq \Ebiv(X/X, -v)
  \end{equation*}
for any s-scheme $X$ over $S$ and $v \in \K(X)$.
\end{rem}

\begin{rem}
The bivariant groups $\Ebiv_\ast(X/S, \ast)$ were previously called \emph{Borel--Moore homology} groups in \cite{Deg16}. This terminology is justified when $S$ is the spectrum of
 a field, and coincides with that of \cite[Chap.~4, Sect.~9]{FSV}. However, in the case where $S$
 is an arbitrary scheme, and especially singular, the homology groups $\Ebiv_\ast(X/S, \ast)$
 are no longer given by the cohomology with coefficients in a dualizing object, which is a
 characteristic property of the original theory of Borel and Moore. For that reason we find the
 terminology ``bivariant'' more suitable.
\end{rem}

\begin{rem}
Given a morphism $f : T \to S$, we can consider the inverse image $\E_T = f^*(\E) \in \SH(T)$ and the associated bivariant theory $\Ebiv_T(-/T, \ast)$ over $T$.
When there is no risk of confusion, we will usually abuse notation by writing $\Ebiv(-/T, \ast) = \Ebiv_T(-/T, \ast)$.
\end{rem}

\begin{rem}
Note that any isomorphism $v\simeq w$ in $\K(X)$ induces an isomorphism of bivariant spectra $\Ebiv(X/S, v) \simeq \Ebiv(X/S, w)$ and of cohomology spectra $\E(X, v) \simeq \E(X, w)$.
More precisely, the assignments $v \mapsto \Ebiv(X/S, v)$ and $v \mapsto \E(X, v)$ are functors on $\K(X)$ (viewed as an \inftyGrpd).
\end{rem}

\begin{notat}
In the notation $\E(X/S, v)$ and $\E(X, v)$, we will sometimes implicitly regard $v$ as a class in $\K(X)$ even if it is actually defined over some deeper base.
For example we will write $\E(X/S, v) = \E(X/S, f^*(v))$ where $f : X \to S$ and $v \in \K(S)$.
\end{notat}

\begin{num}\label{num:bivar_ppties}
The bivariant theory represented by a motivic spectrum $\E \in \SH(S)$ satisfies the following axioms, which are $\K$-graded and spectrum-level refinements of the axioms of Fulton and MacPherson \cite{FMP}:

\begin{enumerate}
\item\label{item:bivar_ppties/base change}\emph{Base change.}
For any cartesian square
  \begin{equation*}
    \begin{tikzcd}[matrix scale=0.7]
      X_T \ar{r}{g}\ar[swap]{d}{q}\ar[phantom]{rd}{\scriptstyle\Delta}
        & X \ar{d}{p}
      \\
      T \ar[swap]{r}{f}
        & S
    \end{tikzcd}
  \end{equation*}
% $$
% \xymatrix@=10pt{
% Y\ar^g[r]\ar_q[d]\ar@{}|\Delta[rd] & X\ar^p[d] \\
% T\ar_f[r] & S,
% }
% $$
there is a canonical base change map
$$
\Delta^*:\Ebiv(X/S,v) \rightarrow \Ebiv(X_T/T,g^*v).$$
This is induced by the natural transformation
    \begin{equation*}
      p_*\Sigma^{-v}p^!
        \xrightarrow{\mathrm{unit}} p_*\Sigma^{-v}p^!f_* f^*
        \simeq p_*\Sigma^{-v}g_* q^! f^*
        \simeq p_*g_*\Sigma^{-g^*v} q^! f^*
        \simeq f_*q_*\Sigma^{-g^*v} q^! f^*,
    \end{equation*}
where we have used the base change formula (Paragraph~\ref{num:six operation compatibilities}\ref{item:six operation compatibilities/base change}) and the formula \eqref{eq:Thom suspension compatibilities}.

\item\label{item:bivar_ppties/proper direct image}\emph{Proper covariance.}
For any proper morphism $f:X \rightarrow Y$ of s-schemes over $S$, there is a direct image map
  $$f_*:\Ebiv(X/S,f^*v) \rightarrow \Ebiv(Y/S,v).$$
This covariance is induced by the unit map $f_!f^! \rightarrow \Id$ and the identification $f_! \simeq f_*$ as $f$ is proper (Paragraph~\ref{num:six operation compatibilities}\ref{item:six operation compatibilities/proper}).

\item\label{item:bivar_ppties/etale inverse image}\emph{\'Etale contravariance.}
For any \'etale s-morphism $f:X \rightarrow Y$ of s-schemes over $S$, there is an inverse image map
  $$f^!:\Ebiv(Y/S,v) \rightarrow \Ebiv(X/S,f^*v).$$
This contravariance is induced by the purity isomorphism $\pur_f : f^! \simeq f^*$ (Paragraph~\ref{num:smooth purity}).

\item\label{item:bivar_ppties/product}\emph{Product.}
If $\E$ is equipped with a multiplication map $\mu_\E : \E \otimes \E \to \E$, then for s-morphisms $p : X \to S$ and $q : Y \to X$, and any K-theory classes $v \in \K(X)$ and $w \in \K(Y)$, there is a map
  $$
  \Ebiv(Y/X,w) \otimes \Ebiv(X/S,v) \rightarrow \Ebiv(Y/S,w+q^*v).
  $$
Given maps
 $y:\Th_Y(w)[m] \rightarrow q^!\E_X$ and
 $x:\Th_X(v)[n] \rightarrow p^!\E_S$, the product $y.x$ is defined as follows:
\begin{equation*}
  \begin{split}
  \Th_Y&(w+q^*v)[m+n]
    \xrightarrow{y \otimes \Id} q^!\E_X \otimes \Th_Y(q^*v)[n]
    \xrightarrow{Ex^{!*}_\otimes} q^!(\E_X \otimes \Th_Y(v)[n])
  \\
    &\xrightarrow{q^!(\Id\otimes x)} q^!(\E_X \otimes p^!\E_S)
    \xrightarrow{Ex^{!*}_\otimes} q^!p^!(\E_S \otimes \E_S)
    \xrightarrow{\mu_\E} q^!p^!(\E_S) = (pq)^!(\E_S).
  \end{split}
\end{equation*}
\end{enumerate}

These structures satisfy the usual properties stated by Fulton and
 MacPherson (functoriality,
 base change formula both with respect to base change and \'etale contravariance,
 compatibility with pullbacks and projection formulas; see \cite[1.2.8]{Deg16} for the precise formulation).
\end{num}

\begin{rem}
One of the main objectives of this paper concerns the extension of contravariant functoriality from \emph{\'etale} morphisms to \emph{smoothable lci} morphisms.
This will be achieved in Theorem~\ref{thm:Gysin_coef_basics}.
\end{rem}

\begin{rem}
Note that a particular case of the product of Paragraph~\ref{num:bivar_ppties}\ref{item:bivar_ppties/product} is the cap-product:
$$
\cap : \E(X,v) \otimes \E(X/S,w) \rightarrow \E(X/S,w-v).
$$
\end{rem}

The localization theorem (Paragraph~\ref{num:localization}) gives the following direct corollary:

\begin{prop}\label{prop:localization}
Let $i:Z \rightarrow X$ be a closed immersion of s-schemes over $S$,
 with quasi-compact complementary open immersion $j:U \rightarrow X$.
Then there exists, for any $e \in \K(X)$, a canonical exact triangle of spectra
 \begin{equation*}
   \E(Z/S, e)
    \xrightarrow{i_*} \E(X/S, e)
    \xrightarrow{j^*} \E(U/S, e)
 \end{equation*}
called the \emph{localization triangle}.
% \begin{equation}
% \label{eq:loc_les}
% \cdots\to
% \E_n(Z/S,e) \xrightarrow{i_*}
%  \E_n(X/S,e) \xrightarrow{j^{*}}
%  \E_n(U/S,e) \xrightarrow{\partial_i}
%  \E_{n-1}(Z/S,e)
%  \to \cdots
% \end{equation}
% for any $e \in \K(X)$,
Moreover, this triangle is natural with respect to the contravariance in $S$,
 the contravariance in $X/S$ for \'etale $S$-morphisms,
 and the covariance in $X/S$ for proper $S$-morphisms (see parts \ref{item:bivar_ppties/base change}--\ref{item:bivar_ppties/etale inverse image} of Paragraph~\ref{num:bivar_ppties}).
\end{prop}

A special case of the naturality in Proposition~\ref{prop:localization} is the following:

\begin{cor}\label{cor:crossing_localization_les}
Suppose given a commutative square
  \begin{equation*}
    \begin{tikzcd}[matrix scale=0.7]
      T \ar{r}{l}\ar[swap]{d}{k}
        & Z' \ar{d}{j}
      \\
      Z \ar[swap]{r}{i}
        & X
    \end{tikzcd}
  \end{equation*}
of closed immersions of s-schemes over $S$.
Assume that the respective complementary open immersions $i'$, $j'$, $k'$, and $l'$ are quasi-compact.
Then the localization triangles (Proposition~\ref{prop:localization}) assemble into a commutative diagram of spectra
  \begin{equation*}
    \begin{tikzcd}[matrix scale=0.7]
      \Ebiv(T/S, e)\ar{r}{l_*}\ar{d}{k_*}
        & \Ebiv(Z'/S, e)\ar{r}{l'^*}\ar{d}{j_*}
        & \Ebiv(Z'-T/S, e)\ar{d}{\tilde{j}_*}
      \\
      \Ebiv(Z/S, e)\ar{r}{i_*}\ar{d}{k'^*}
        & \Ebiv(X/S, e)\ar{r}{i'^*}\ar{d}{j'^*}
        & \Ebiv(X-Z/S, e)\ar{d}{\tilde{j}'^*}
      \\
      \Ebiv(Z-T/S, e)\ar{r}{\tilde{i}_*}
        & \Ebiv(X-Z'/S, e)\ar{r}{\tilde{i}'^*}
        & \Ebiv(X-Z \cup Z'/S, e)
    \end{tikzcd}
  \end{equation*}
for any $e \in \K(X)$.
Here $\tilde j$, $\tilde i$ denote the obvious closed immersions obtained by restriction, and $\tilde{i}'$, $\tilde{j}'$ the complementary open immersions.
\end{cor}

\begin{prop}\label{prop:multiplication_localization_les}
Suppose $\E \in \SH(S)$ is equipped with a multiplication map $\mu_\E : \E \otimes \E \to \E$.
Consider cartesian squares of s-schemes over $S$:
  \begin{equation*}
    \begin{tikzcd}[matrix scale=0.7]
      T \ar{r}{k}\ar{d}\ar[phantom]{rd}{\scriptstyle\Delta_Z}
        & Y \ar{d}
        & V \ar[swap]{l}{k'}\ar{d}\ar[phantom]{ld}{\scriptstyle\Delta_U}
      \\
      Z \ar[swap]{r}{i}
        & X
        & U \ar{l}{i'}
    \end{tikzcd}
  \end{equation*}
such that $i$ and $k$ are closed immersions with quasi-compact complementary open immersions $i'$ and $k'$, respectively.
For any $\pi \in \E(Y/X,e')[r]$ with $e' \in \K(Y)$, $r \in \ZZ$, set:
$$
\pi_Z=\Delta_Z^*(\pi) \in \E(T/Z,e')[r],
\quad
\pi_U=\Delta_U^*(\pi) \in \E(V/U,e')[r].
$$
Then the following diagram of localization triangles is commutative:
  \begin{equation*}
    \begin{tikzcd}[matrix scale=0.7]
      \E(Z/S, e)\ar{r}{i_*}\ar{d}{\gamma_{\pi_Z}}
        & \E(X/S,e)\ar{r}{i'^*}\ar{d}{\gamma_{\pi}}
        & \E(U/S,e)\ar{d}{\gamma_{\pi_U}}
      \\
      \E(T/S, e+e')[r]\ar{r}{k_*}
        & \E(Y/S, e+e')[r]\ar{r}{k'^*}
        & \E(V/S, e+e')[r]
    \end{tikzcd}
  \end{equation*}
% $$
% \xymatrix@R=12pt@C=18pt{
% \E_n(Z/S,e)\ar^-{i_*}[r]\ar_{\gamma_{\pi_Z}}[d] & \E_n(X/S,e)\ar^-{i^{\prime*}}[r]\ar^{\gamma_{\pi}}[d]
%  & \E_n(U/S,e)\ar^-{\partial_{i}}[r]\ar_{\gamma_{\pi_U}}[d] & \E_{n-1}(Z/S,e)\ar^{\gamma_{\pi_Z}}[d] \\
% \E_{n+r}(T/S,e+e')\ar^-{k_*}[r] & \E_{n+r}(Y/S,e+e')\ar^-{k^{\prime*}}[r] & \E_{n+r}(V/S,e+e')\ar^-{\partial_k}[r] & \E_{n+r-1}(T/S,e+e')
% }
% $$
where $\gamma_x$ denotes multiplication by $x\in\{\pi, \pi_Z, \pi_U\}$.
\end{prop}

\begin{proof}
The left-hand square commutes by the projection formula.
The right-hand square commutes since products are compatible with base change.
\end{proof}

\begin{num}\label{num:bivar_product}
Suppose that $\E$ is equipped with a multiplication $\mu_\E : \E \otimes \E \to \E$ as in Paragraph~\ref{num:bivar_ppties}\ref{item:bivar_ppties/product}.
If this multiplication is unital, resp. associative, resp. commutative, then the bivariant theory represented by $\E$ inherits the same property.
This is in particular the case when $\E$ is equipped with an $\cE_\infty$-ring structure.

For example, assume that the multiplication is commutative in the sense that it is further equipped with a commutative diagram in $\SH(S)$:
  \begin{equation*}
    \begin{tikzcd}[matrix scale=0.7]
      \E\otimes\E \ar{r}{\sim}[swap]{\tau}\ar[swap]{d}{\mu_\E}
        & \E\otimes\E \ar{d}{\mu_\E}
      \\
      \E \ar[equals]{r}
        & \E,
    \end{tikzcd}
  \end{equation*}
where $\tau$ is the isomorphism swapping the two factors.
Given s-schemes $p:X\to S$ and $q:Y\to S$ and $x \in \E(X/S,v)[-m]$ and $y \in \E(Y/S,w)[-n]$ (where $m,n\in\ZZ$ and $v\in\K(X)$, $w\in\K(Y)$), consider the cartesian square
  \begin{equation*}
    \begin{tikzcd}[matrix scale=0.7]
      X\times_S Y \ar{r}{\pi_2}\ar[swap]{d}{\pi_1}\ar[phantom]{rd}{\scriptstyle\Delta}
        & Y \ar{d}{q}
      \\
      X \ar[swap]{r}{p}
        & S.
    \end{tikzcd}
  \end{equation*}
Then there is an identification
  $$
  \Delta^*(x).y \simeq (-1)^{m+n}.\big(\Delta^*(y).x\big)
  $$
in $\E(X \times_S Y/S,\pi_1^*v+\pi_2^*w)[-m-n]\simeq \E(X \times_S Y/S,\pi_2^*w+\pi_1^*v)[-m-n]$ using the permutation isomorphism
 $\pi_2^*w+\pi_1^*v \simeq \pi_1^*v+\pi_2^*w$.
\end{num}

\begin{num}
For any $\E \in\SH(S)$ and $v \in \K(S)$, the functor $X \mapsto \E(X/S, v)$ satisfies descent with respect to Nisnevich squares and abstract blow-up squares (hence satisfies cdh descent), on the category of s-schemes over $S$.
\end{num}

%!TEX root = ../../virtual.tex

\subsection{Orientations and systems of fundamental classes}
\label{sec:orientation}

Following Fulton--MacPherson, we now introduce the notion of \emph{orientation} of a morphism $f$.
As we recall in the next subsection, any choice of orientation gives rise to a Gysin map in bivariant theory (Paragraph~\ref{df:orientation gives Gysin}).
The fundamental classes we construct in Section~\ref{sec:construction} will be examples of orientations.

For simplicity, throughout this discussion we will restrict our attention to the bivariant theory represented by the sphere spectrum $\un$:

\begin{notat}\label{notat:bivariant theory for E=S}
We set 
  \begin{equation*}
    H(X/S, v) := \un(X/S, v) = \Maps_{\SH(X)}(\Th_X(v),p^!(\un_S))
  \end{equation*}
for any s-morphism $p : X \to S$ and any $v \in\K(X)$.
We will refer to this simply as \emph{bivariant $\AA^1$-theory}.
Similarly, we set $H(X,v) := \un(X,v)$ and more generally $H_Z(X,v) := \un_Z(X, v) = \un(Z/X, -v)$, where $Z$ is a closed subscheme of $X$.
\end{notat}

\begin{df}\label{df:orientation}
Let $S$ be a scheme and $f:X \rightarrow S$ an s-morphism.
An \emph{orientation} of $f$ is a pair $(\eta_f, e_f)$, where
$$
\eta_f \in H(X/S,e_f)
$$
and $e_f \in \K(X)$.
When there is no risk of confusion, we write simply $\eta_f$ instead of $(\eta_f,e_f)$.
\end{df}

\begin{rem}
The above use of the term ``orientation'' is taken from \cite{FMP}.
We warn the reader however that it is \emph{unrelated} to the notion of ``oriented motivic spectrum'' (see Definition~\ref{df:oriented_spectra}).
\end{rem}

\begin{ex}\label{ex:orientation for smooth}
Let $f:X \rightarrow S$ be a smooth s-morphism with tangent bundle $T_f$.
The purity isomorphism $\pur_f : \Sigma^{T_f} f^* \xrightarrow{\sim} f^!$ (Paragraph~\ref{num:smooth purity}) induces a canonical isomorphism
$$
\fdl_f:\Th_X(T_f) \xrightarrow \sim f^!(\un_S).
$$
This defines a canonical orientation $\fdl_f \in H(X/S,\vb{T_f})$.

It will be useful to have the following description of the purity isomorphism $\pur_f$ (cf. \cite[Def.~2.4.25, Cor.~2.4.37]{CD3}).
We begin by considering the cartesian square $\Delta$:
  \begin{equation*}
    \begin{tikzcd}[matrix scale=0.7]
      X\times_S X \ar{r}{f_1}\ar[swap]{d}{f_2}
        & X \ar{d}{f}
      \\
      X \ar[swap]{r}{f}
        & S
    \end{tikzcd}
  \end{equation*}
% $$
% \xymatrix@=14pt{
% X\times_S X\ar^-{f_1}[r]\ar_-{f_2}[d] & X\ar^f[d] \\
% X\ar^f[r] & S
% }
% $$
Write $\delta:X \rightarrow X\times_S X$ for the diagonal embedding.
Then $\Sigma^{-T_f}\pur_f$ is inverse to the composite:
  \begin{equation}\label{eq:pur_smooth}
    \Sigma^{-T_f} f^!
      \xrightarrow{\sim} \delta^!f_1^*f^!
      \xrightarrow{Ex^{*!}} \delta^!f_2^!f^* = f^*
  \end{equation}
Here the first isomorphism is induced by the relative purity isomorphism $\Sigma^{-T_f} \simeq \delta^!f_1^*$ (Paragraph~\ref{num:relative purity}), modulo the tautological identification between $T_f$ and the normal bundle $N_\delta$.
The exchange transformation $Ex^{*!}$ is invertible because $f$ is smooth (see Paragraph~\ref{num:Ex^*!}).
\end{ex}

\begin{df}\label{df:fdl_smooth}
Let $f : X \to S$ be a smooth s-morphism.
The \emph{fundamental class} of $f$ is the orientation $\eta_f \in H(X/S, \vb{T_f})$ defined in Example~\ref{ex:orientation for smooth}.
\end{df}

\begin{df}\label{df:fdl_classes_abstract}
Let $S$ be a scheme and let $\C$ be a class of morphisms between s-schemes over $S$.
A \emph{system of fundamental classes} for $\C$ consists of the following data:
\begin{remlist} 
  \item \textit{Fundamental classes}.
  For each morphism $f : X \to Y$ in $\C$, there is an orientation $(\eta^\C_f,e_f)$.

  \item\label{item:fdl_classes_abstract/identity}\textit{Normalisation}.
  For $f = \Id_S$ the identity morphism, there is an isomorphism $e_f \simeq 0$ in $\K(S)$, and an isomorphism $\eta^\C_f \simeq 1$ in $H(S/S, e_f) \simeq H(S/S, 0)$.

  \item\label{item:fdl_classes_abstract/functorial}\textit{Associativity formula}.
  Let $f : X \to Y$ and $g : Y \to Z$ be morphisms in $\C$ such that the composite $g\circ f$ is also in $\C$.
  Then there are identifications
    \begin{equation}\label{eq:comp_virt_coef}
      e_{g\circ f} \simeq e_f+f^*(e_g)
    \end{equation}
  in $\K(X)$ and
    $$ \eta^\C_g.\eta^\C_f \simeq \eta^\C_{g \circ f} $$
  in $H(X/Z, e_{g\circ f})$.
\end{remlist}

We say that a system of fundamental classes $(\eta^\C_f)_f$ is \emph{stable under transverse base change} if it is equipped with the following further data:
\begin{remlist}[resume]
\item\label{item:fdl_classes_abstract/base change}\emph{Transverse base change formula.}
For any \emph{tor-independent} cartesian square
  \begin{equation*}
    \begin{tikzcd}[matrix scale=0.7]
      Y\ar{r}{g}\ar[swap]{d}{q}\ar[phantom]{rd}{\scriptstyle\Delta}
        & T \ar{d}{p}
      \\
      X\ar[swap]{r}{f}
        & S
    \end{tikzcd}
  \end{equation*}
% $$
% \xymatrix@=14pt{
% Y\ar^g[r]\ar_q[d]\ar@{}|\Delta[rd]
%  & T \ar^p[d] \\
% X\ar_{f}[r] & S.
% }
% $$
such that $f$ and $g$ are in $\C$, there are identifications $e_g \simeq q^*(e_f)$ in $\K(Y)$ and
  $\Delta^*(\eta^\C_f) \simeq \eta^\C_g$
in $H(X/S, e_g)$.
\end{remlist}
\end{df}

\begin{rem}\label{rem:orientation for general E}
The previous definition admits an obvious extension
 to general bivariant theories (i.e., the contexts of Definition~\ref{df:bivariant} and Paragraph~\ref{num:motivic_categories}),
 and we will freely use this extension.
 Then our definition is both a generalization of \cite[I, 2.6.2]{FMP}
 and of \cite[2.1.9]{Deg16}.
\end{rem}

\begin{rem}\label{rem:compose orientations}
Let $S$ be a scheme and let $\C$ be a class of morphisms between s-schemes over $S$.
Suppose $(\eta_f)_f$ is a system of fundamental classes for $\C$ as in Definition~\ref{df:fdl_classes_abstract}.
For any morphisms $f : X \to Y$ and $g : Y \to Z$ in $\C$ such that the composite $g\circ f$ is also in $\C$, the isomorphisms \eqref{eq:comp_virt_coef} induce canonical isomorphisms of functors $\SH(Z) \to \SH(X)$
  \begin{equation*}
    (g\circ f)^*\Sigma^{e_{gf}} \simeq \Sigma^{e_{f}}f^*\Sigma^{e_{g}}g^*,
  \end{equation*}
and isomorphisms of Thom spaces
  \begin{align*}
    \Th(e_{gf}) &\simeq \Th(f^*e_{g}) \otimes \Th(e_{f})\\
      &\simeq f^*\Th(e_{g}) \otimes \Th(e_{f})
  \end{align*}
in $\SH(X)$.
\end{rem}

\begin{ex}\label{ex:Ayoub&assoc_sm}
It follows from \cite[1.7.3]{Ayo1} that the  
 family of orientations $\eta_f$ for $f$ smooth (Definition~\ref{df:fdl_smooth})
 forms a system of fundamental classes for the class of smooth s-morphisms.
Moreover, this system is stable under (arbitrary) base change: explicit homotopies $\Delta^*(\eta_f) \simeq \eta_g$ as in Definition~\ref{df:fdl_classes_abstract}\ref{item:fdl_classes_abstract/base change} are provided by the deformation to the normal cone space, as in the proof of \cite[Lem.~2.3.13]{Deg16} (where the right-hand square (3) can be ignored).
\end{ex}

\begin{ex}\label{ex:lci basics}
In Section~\ref{sec:construction}, we will extend Example~\ref{ex:Ayoub&assoc_sm} to the class of \emph{smoothable lci} s-morphisms.
Recall from \cite{IllusieCotangentI} that an lci morphism $f:X \rightarrow S$ admits a perfect cotangent complex $L_f = L_{X/S}$, which induces a point $\vb{L_f} \in \K(X)$ (which represents the ``virtual tangent bundle'' of $f$ in $\K_0(X)$ in the sense of \cite[Exp. VIII]{SGA6}).
For example, if $f$ is smooth then $\vb{L_f} = \vb{T_f}$ is the class of the relative tangent bundle; if $f$ is a regular closed immersion then $\vb{L_f} = -\vb{N_f}$, where $\vb{N_f}$ is the class of the normal bundle.
Every \emph{smoothable} lci morphism $f$ factors through a regular closed immersion $i$ followed by a smooth morphism $p$, and such a factorization induces an identification $\vb{L_f} \simeq i^*\vb{T_p} - \vb{N_i}$ in K-theory.
More generally, given lci morphisms $f : X \to Y$ and $g : Y \to Z$, the composite $g\circ f$ is lci with cotangent complex canonically identified in $\K(Y)$ with
  \begin{equation}\label{eq:comp_virt_tg}
    \vb{L_{g\circ f}} \simeq \vb{L_f}+f^{*}\vb{L_g}.
  \end{equation}
The fundamental class of a smoothable lci s-morphism $f : X \to S$ will then be an orientation in $H(X/S, \vb{L_f})$ (see Theorem~\ref{thm:fdl_qplci}).
\end{ex}

We finish with a discussion of \emph{strong} orientations and duality isomorphisms.

\begin{df}\label{df:strong orientation}
Let $f:X \rightarrow S$ be an s-morphism, and $(\eta_f,e_f)$ an orientation of $f$.

\begin{enumerate}
\item
We say that $(\eta_f,e_f)$ is \emph{strong} if for any $v \in \K(X)$, cap-product with $\eta_f$ induces an isomorphism
$$
\gamma_{\eta_f}:H(X,v) \rightarrow H(X/S,e_f-v), \quad x \mapsto x.\eta_f.
$$
In that case, we refer to $\gamma_{\eta_f}$ as the \emph{duality isomorphism} associated with the strong orientation $\eta_f$.

\item
We say that $(\eta_f,e_f)$ is \emph{universally strong} if the morphism $\eta_f: \Th_X(e_f) \rightarrow f^!(\un_S)$ is an isomorphism.
\end{enumerate}
\end{df}

\begin{rem}
It follows immediately from the construction of the cap product that a universally strong orientation is strong.
\end{rem}

\begin{ex}\label{ex:strong orientation of smooth}
If $f$ is \emph{smooth}, then the orientation $\fdl_f$ of Definition~\ref{df:fdl_smooth} is universally strong by the purity theorem (Paragraph~\ref{num:smooth purity}).
\end{ex}

The following lemma explains the terminology ``universally strong''.

\begin{lm}\label{lm:universally strong orientation}
Let $f:X \rightarrow S$ be an s-morphism.
Let $(\eta_f, e_f)$ be a \emph{universally strong} orientation for $f$.
Then for any cartesian square
  \begin{equation*}
    \begin{tikzcd}[matrix scale=0.7]
      X_T \ar{r}{q}\ar[swap]{d}{g}\ar[phantom]{rd}{\scriptstyle\Delta}
        & X \ar{d}{f}
      \\
      T \ar[swap]{r}{p}
        & S
    \end{tikzcd}
  \end{equation*}
% $$
% \xymatrix@=10pt{
% X_T\ar^q[r]\ar_g[d]\ar@{}|\Delta[rd] & X\ar^f[d] \\
% T\ar_p[r] & S
% }
% $$
with $p$ smooth, the orientation $\Delta^*(\eta_f) \in H(X_T/T, q^*(e_f))$ of $g$ is universally strong.
\end{lm}

\begin{proof}
Since $p$ is smooth, the exchange transformation $Ex^{*!} : p^*f^! \rightarrow g^!q^*$ (Paragraph~\ref{num:Ex^*!}) is invertible.
It follows then from the construction of the change of base map $\Delta^*$ that $\Delta^*(\eta_f) : \Th_{X_T}(q^*(e_f)) \to g^!p^*(\un_S)$ is an isomorphism, as claimed.
\end{proof}

%!TEX root = ../../virtual.tex

\subsection{Gysin maps}
\label{sec:gysin}

\begin{df}\label{df:orientation gives Gysin}
Let $f : X \to Y$ be a morphism of s-schemes over $S$.
Then any orientation $\eta_f \in H(X/Y,e_f)$ gives rise to a \emph{Gysin map}:
$$
\eta_f^! : H(Y/S,e) \rightarrow H(X/S,e_f+f^*(e)),
  \quad x \mapsto\fdl_f.x
$$
using the product in bivariant $\AA^1$-theory, for all $e\in \K(Y)$.
 When the orientation $\eta_f$ is clear, we simply put: $f^!=\eta_f^!$.
\end{df}

\begin{prop}\label{prop:abstract Gysin}
Let $S$ be a scheme and let $\C$ be a class of morphisms between s-schemes over $S$.
Suppose $(\eta_f)_f$ is a system of fundamental classes for $\C$ as in Definition~\ref{df:fdl_classes_abstract}.
\begin{thmlist}
\item\label{item:abstract Gysin/functoriality}\emph{Functoriality.}
Let $f$ and $g$ be morphisms in $\C$ such that the composite $g\circ f$ is also in $\C$.
Then for every $e \in \K(Z)$ there is an induced identification
  \begin{equation*}
    (g\circ f)^! \simeq f^! \circ g^!
  \end{equation*}
of Gysin maps $H(Z/S, e) \to H(X/S, e_{g\circ f}+(g\circ f)^*(e))$, modulo the identification of Remark~\ref{rem:compose orientations}.

\item\label{item:abstract Gysin/base change}\emph{Transverse base change.}
Assume that the system $(\eta_f)_f$ is stable under transverse base change.
Suppose given a \emph{tor-independent} cartesian square of s-schemes over $S$ of the form
  \begin{equation*}
    \begin{tikzcd}[matrix scale=0.7]
      X' \ar{r}{g}\ar{d}[swap]{u}
        & Y'\ar{d}{v}
      \\
      X \ar{r}{f}
        & Y,
    \end{tikzcd}
  \end{equation*}  
where $f$ and $g$ are in $\C$, and $u$ and $v$ are proper.
Then for every $e \in \K(Y)$ there is an induced identification
  \begin{equation*}
    f^! \circ v_* \simeq u_* \circ g^!
  \end{equation*}
of maps $H(Y'/S, v^*(e)) \to H(X/S, e_f+f^*(e))$, where we use the identification $e_g+g^*v^*(e) \simeq u^*(e_f + f^*(e))$ in $\K(X')$.
\end{thmlist}
\end{prop}

\begin{ex}
If $f$ is a smooth s-morphism, then the fundamental class $\eta_f$ (Definition~\ref{df:fdl_smooth}) gives rise to canonical Gysin maps
  \begin{equation*}
    f^! : H(Y/S,e) \rightarrow H(X/S,\vb{T_f}+f^*(e)).
  \end{equation*}
This extends the contravariant functoriality from \'etale morphisms to smooth morphisms.
\end{ex}

\begin{lm}\label{lm:htp_inv_biv}
Let $X$ be an s-scheme over $S$ and let $p : E \to X$ be a vector bundle.
Then the tangent bundle $T_p$ is identified with $p^{-1}(E)$ and the Gysin map
  $$p^!:H(X/S,e) \rightarrow H(E/S,p^*\vb{E}+p^*e)$$
is invertible.
\end{lm}

\begin{proof}
In view of the construction of the Gysin map, the claim follows directly from the facts that the morphism $\fdl_p : \Th_X(p^*E) \to p^!(\un_X)$ is invertible, and that the functor $p^* : \SH(X) \to \SH(E)$ is fully faithful (Paragraph~\ref{num:six/A^1-invariance}).
\end{proof}

\begin{df}\label{df:Thom}
In the context of Lemma~\ref{lm:htp_inv_biv}, we define the \emph{Thom isomorphism}\footnote{Not to be confused with the Thom isomorphism of \eqref{eq:Thom_iso}.}
$$
\phi_{E/X} : H(E/S,e) \rightarrow H(X/S,e-\vb{E}),
$$
associated to $E/X$, as the inverse of the Gysin map $p^! : H(X/S,e-\vb{E}) \to H(E/S,e)$.
\end{df}

\begin{rem}
The Thom isomorphism satisfies the properties of compatibility with base change and with direct sums (that is, $\phi_{E \oplus F/X}\simeq \phi_{E \oplus F/F} \circ \phi_{F/X}$ for vector bundles $E$ and $F$ over $X$).
These follow respectively from the compatibility of the Gysin morphisms $p^!$ with base change and with composition.
\end{rem}

We conclude this subsection by recording the naturality
 of the localization sequences (Proposition~\ref{prop:localization})
 with respect to Gysin maps of smooth morphisms.

\begin{prop}\label{prop:local&Gysin_sm}
Consider a commutative diagram of cartesian squares of $S$-schemes
  \begin{equation*}
    \begin{tikzcd}[matrix scale=0.7]
      T \ar{r}{k}\ar[swap]{d}{g}
        & Y \ar{d}{f}
        & V\ar[swap]{l}{l}\ar{d}{h}
      \\
      Z \ar[swap]{r}{i}
        & X
        & U \ar{l}{j}
    \end{tikzcd}
  \end{equation*}
% $$
% \xymatrix@R=10pt@C=20pt{
% T\ar^k[r]\ar_g[d] & Y\ar^f[d] & V\ar^h[d]\ar_l[l] \\
% Z\ar^i[r] & X  & U\ar_j[l]
% }
% $$
such that $f$ is smooth, $i$ is a closed immersion
 and $j$ is the complementary open immersion.
Then for any $e \in \K(X)$, the diagram of spectra
  \begin{equation*}
    \begin{tikzcd}[matrix scale=0.7]
      H(Z/S,e)\ar{r}{i_*}\ar{d}{g^!}
        & H(X/S,e)\ar{r}{j^!}\ar{d}{f^!}
        & H(U/S,e)\ar{d}{h^!}
      \\
      H(T/S,\vb{T_g}+e)\ar{r}{k_*}
        & H(Y/S,\vb{T_f}+e)\ar{r}{l^!}
        & H(V/S,\vb{T_h}+e)
    \end{tikzcd}
  \end{equation*}
% $$
% \xymatrix@=12pt{
% H(Z/S,e_Z)\ar^-{i_*}[r]\ar_{g^!}[d] & H(X/S,e)\ar^-{j^!}[r]\ar_{f^!}[d]
%  & H(U/S,e_U)\ar^-{\partial_{Z/X}}[r]\ar_{h^!}[d] & H(Z/S,e_Z)[-1]\ar^{g^!}[d] \\
% H(T/S,e_T+\vb{T_g})\ar^-{k_*}[r] & H(Y/S,e_Y+\vb{T_f})\ar^-{l^!}[r] & H(V/S,e_V+\vb{T_h})\ar^-{\partial_{T/Y}}[r]
%  & H(T/S,e_T+\vb{T_g})[-1]
% }
% $$
commutes.
\end{prop}
\begin{proof}
The right-hand square commutes by the associativity property of fundamental classes of smooth morphisms (Example~\ref{ex:Ayoub&assoc_sm}).
The left-hand square commutes by the naturality of the relative purity isomorphism of Morel--Voevodsky \cite[Prop. A.4]{HoyoisLefschetz}, in view of the construction of the fundamental classes $\eta_f$ and $\eta_g$ (Example~\ref{ex:orientation for smooth}).
\end{proof}

%!TEX root = ../../virtual.tex

\subsection{Purity transformations}
\label{sec:bivariant/purity}

The notion of orientation seen in the preceding subsection is part of our twisted version
 of the 
 bivariant formalism of Fulton and MacPherson. We state in this subsection a variant,
 or a companion, of this notion in the spirit of Grothendieck's six functors formalism.

\begin{num}\label{num:purity_transfo}
Let us fix an s-morphism $f:X \rightarrow S$, and an orientation $(\eta_f, e_f)$.  According to our definitions,
 the class $\eta_f \in H(X/S, e_f)$ can be seen as a morphism in $\SH(X)$:
$$
\fdl_f : \Th_X(e_f) \to f^!(\un_S).
$$
This gives rise to a natural transformation
  \begin{equation}
    \pur(\eta_f) : \Sigma^{e_f} f^{*} \to f^!
  \end{equation}
associated to the orientation $\eta_f$, defined as the following
 composite:
$$
f^{*}(-) \otimes \Th_X(e_f)
\xrightarrow{\Id\otimes\fdl_f} f^*(-) \otimes f^!(\un_S)
\xrightarrow{Ex^{!*}_\otimes} f^!(- \otimes \un_S) \simeq f^!,
$$
where the exchange transformation $Ex^{!*}_\otimes$ is as in Paragraph~\ref{num:Ex^!*_otimes}.
\end{num}

\begin{rem}\label{rem:purity_transfo}\leavevmode
\begin{remlist}
\item Suppose $f$ is smooth and consider the canonical orientation $\eta_f$ of Definition~\ref{df:fdl_smooth}.
It follows by construction that in this case the associated purity transformation $\pur(\eta_f)$ is nothing else than the purity isomorphism $\pur_f$ \eqref{eq:pur_smooth}.
In particular, $\pur(\eta_f)$ is an isomorphism in this case.

% The remark below only makes sense for the specific orientations we will define later...
% \item For non-smooth morphisms $f$, the purity transformation of the above definition
%  is in general not an isomorphism. A very interesting question is to know when
%  it becomes an isomorphism when applied to a specific object
%  (see Definitions \ref{df:f-pure} and 
%  \ref{df:absolute_purity}).

\item Note that the datum of an orientation $(\eta_f,e_f)$ and that of the associated purity transformation $\pur(\eta_f)$ are essentially interchangeable.
Indeed we recover $\eta_f$ by evaluating $\pur(\eta_f)$ at the unit object $\un_S$.
\end{remlist}
\end{rem}

\begin{num}\label{num:trace_cotrace}
Consider the notation of the previous definition. Then one associates to $\pur(\eta_f)$,
 using the adjunction properties, two natural transformations:
\begin{align*}
\tr_f&:f_!\Sigma^{e_f}f^* \rightarrow \Id \\
\cotr_f&:\Id \rightarrow f_*\Sigma^{-e_f}f^!.
\end{align*}
The first (resp. second) natural transformation 
 will be called the \emph{trace map} (resp. \emph{co-trace map}) associated with the orientation
 $\eta_f$, following the classical usage in the literature.
 These two maps are functorial incarnations of the Gysin map defined earlier
 (Paragraph~\ref{df:orientation gives Gysin}), as we will see later (see Paragraph~\ref{num:Gysin four}).
\end{num}

The notion of a system of fundamental classes (Definition \ref{df:fdl_classes_abstract})
 was introduced to reflect the functoriality of Gysin morphisms.
 For completeness, we now formulate the analogous functoriality property
 for the associated purity transformations.

\begin{prop}\label{prop:2-functor}
Let $S$ be a scheme and let $\C$ be a class of morphisms between s-schemes over $S$.
Suppose $(\eta_f)_f$ is a system of fundamental classes for $\C$ as in Definition~\ref{df:fdl_classes_abstract}.
\begin{thmlist}
\item\label{item:2-functor/functoriality}\emph{Functoriality.}
Let $f$ and $g$ be morphisms in $\C$ such that the composite $g\circ f$ is also in $\C$.
Then there is a commutative square
  \begin{equation*}
    \begin{tikzcd}[matrix scale=0.7, column sep=60]
      \Sigma^{e_{g\circ f}} (g\circ f)^* \ar{r}{\pur(\eta_{g\circ f})}\ar[equals]{d}
        & (g\circ f)^! \ar[equals]{d}
      \\
      \Sigma^{e_{f}}f^*\Sigma^{e_{g}}g^* \ar{r}{\pur(\eta_f) \ast \pur(\eta_g)}
        & f^!g^!
    \end{tikzcd}
  \end{equation*}
where the left-hand vertical isomorphism is from Remark~\ref{rem:compose orientations}, and the lower horizontal arrow is the horizontal composition of the $2$-morphisms $\pur(\eta_f)$ and $\pur(\eta_g)$.

\item\label{item:2-functor/base change}\emph{Transverse base change.}
Assume that the system $(\eta_f)_f$ is stable under transverse base change.
Suppose given a \emph{tor-independent} cartesian square of the form
  \begin{equation*}
    \begin{tikzcd}[matrix scale=0.7]
      X' \ar{r}{g}\ar{d}{u}
        & Y'\ar{d}{v}
      \\
      X \ar{r}{f}
        & Y,
    \end{tikzcd}
  \end{equation*}  
where $f$ and $g$ are in $\C$.
Then there are commutative squares of natural transformations:
  \begin{equation*}
    \begin{tikzcd}[matrix xscale=0.7]
      u^*\Sigma^{e_f}f^* \ar[equals]{r}\ar{d}{u^*\ast\pur_f}\ar[phantom]{rrd}{\scriptstyle\mathrm{(a)}}
        & \Sigma^{e_g}u^*f^* \ar[equals]{r}
        & \Sigma^{e_g}g^*v^* \ar{d}{\pur_g\ast v^*}
      \\
      u^*f^! \ar[swap]{rr}{Ex^{*!}}
        &
        & g^!v^*
    \end{tikzcd}
    \quad
    \begin{tikzcd}[matrix xscale=0.7]
      \Sigma^{e_f}f^*v_* \ar{r}{Ex^*_*}\ar{d}{\pur_f\ast v_*}\ar[phantom]{rrd}{\scriptstyle\mathrm{(b)}}
        & \Sigma^{e_f}u_*g^* \ar[equals]{r}
        & u_*\Sigma^{e_g}g^* \ar{d}{u_*\ast\pur_g}
      \\
      f^!v_* \ar{rr}{\sim}
        &
        & u_*g^!
    \end{tikzcd}
  \end{equation*}
If $u$ and $v$ are s-morphisms, then there are also commutative squares
  \begin{equation*}
    \begin{tikzcd}[matrix xscale=0.7]
      u^!\Sigma^{e_f}f^* \ar[equals]{r}\ar{d}{u^!\ast\pur_f}\ar[phantom]{rrd}{\scriptstyle\mathrm{(c)}}
        & \Sigma^{e_g}u^!f^* \ar[equals]{r}
        & \Sigma^{e_g}g^*v^! \ar{d}{\pur_g\ast v^!}
      \\
      u^!f^! \ar[equals]{rr}
        &
        & g^!v^!
    \end{tikzcd}
    \quad
    \begin{tikzcd}[matrix xscale=0.7]
      u_!\Sigma^{e_g}g^* \ar[equals]{r}\ar{d}{u_!\ast\pur_g}\ar[phantom]{rrd}{\scriptstyle\mathrm{(d)}}
        & \Sigma^{e_f}u_!g^* \ar{r}{\sim}
        & \Sigma^{e_f}f^*v_! \ar{d}{\pur_f\ast v_!}
      \\
      u_!g^! \ar[swap]{rr}{Ex^!_!}
        &
        & f^!v_!
    \end{tikzcd}
  \end{equation*}
\end{thmlist}
\end{prop}

\begin{proof}
Claim~\ref{item:2-functor/functoriality} follows from axioms \ref{item:fdl_classes_abstract/identity} and \ref{item:fdl_classes_abstract/functorial} of Definition~\ref{df:fdl_classes_abstract}.

In claim~\ref{item:2-functor/base change}, commutativity of (a) follows directly from axiom~\ref{item:fdl_classes_abstract/base change} of Definition~\ref{df:fdl_classes_abstract}.
The square (b) can be derived from (a) by applying $u_*$ on the left and $v_*$ on the right, and using the naturality of the unit and counit of the adjunctions $(u^*, u_*)$ and $(v^*, v_*)$, respectively.
Similarly, commutativity of square (c) will follow similarly from (d).
For (d), we may unravel the definition of the purity transformation (Paragraph~\ref{num:purity_transfo}) to write the square as follows:
\begin{equation*}
  \begin{tikzcd}[matrix xscale=0.5]
    u_!\Sigma^{e_g}g^* \ar[equals]{r}\ar[equals]{d}
      & \Sigma^{e_f}u_!g^* \ar{r}{\sim}
      & \Sigma^{e_f}f^*v_! \ar[equals]{d}
    \\
    u_!(g^*(-) \otimes \Sigma^{e_g}g^*(\un)) \ar{d}{u_!(g^*(-)\otimes\pur_g(\un))}
      & 
      & f^*v_!(-) \otimes \Sigma^{e_f}f^*(\un) \ar{d}{f^*v_!(-)\otimes\pur_f(\un)}
    \\
    u_!(g^*(-) \otimes g^!(\un)) \ar{d}{Ex^{*!}_\otimes}
      & 
      & f^*v_!(-) \otimes f^!(\un) \ar{d}{Ex^{*!}_\otimes}
    \\
    u_!g^! \ar[swap]{rr}{Ex^!_!}
      &
      & f^!v_!
  \end{tikzcd}
\end{equation*}
Observing that $\pur_g(\un) = \pur_g(v^*(\un))$, we can use square (a) to decompose the left-hand middle arrow into a natural transformation induced by $\pur_f(\un)$ and an exchange transformation.
The commutativity of the resulting square is then a formal exercise.
\end{proof}

\begin{rem}
Suppose that the class $\C$ contains all identity morphisms and is closed under composition, and let $\base^{\C}$ denote the subcategory of the category of schemes $\base$ whose morphisms all belong to $\C$.
At the level of homotopy categories, Proposition~\ref{prop:2-functor}\ref{item:2-functor/functoriality} implies that the assignment $f \mapsto \pur(\eta_f)$ defines a natural transformation of contravariant pseudofunctors
  \begin{equation*}
    \pur : \Ho(\SH)^{e*} \to \Ho(\SH)^!
  \end{equation*}
on the category $\base^{\C}$, where the notation is as follows:
  \begin{itemize}
    \item
    $\Tri$ denotes the $(2,1)$-category of large triangulated categories, triangulated functors, and invertible triangulated natural transformations.
    \item $\Ho(\SH)^!$ is the pseudofunctor $(\base^\C)^\op \to \Tri$, given by the assignments
$$
S \mapsto \Ho(S\mathscr H(S)), f \mapsto f^!.
$$
    \item
    Similarly $\Ho(\SH)^{e*}$ is the pseudofunctor $(\base^\C)^\op \to \Tri$ given by the assignments 
$$
S \mapsto \Ho(\SH(S)), f \mapsto \Sigma^{e_f} f^*.
$$
  \end{itemize}
The expected enhancement to a natural transformation at the level of \inftyCats requires further work that we do not undertake in this paper.
\end{rem}

We can also reformulate the transverse base change property (Proposition~\ref{prop:2-functor}\ref{item:2-functor/base change}) in terms of the (co)trace maps.

\begin{cor}\label{cor:base change in terms of traces}
Under the assumptions of Proposition~\ref{prop:2-functor}\ref{item:2-functor/base change}, the following diagrams commute:
  \begin{align*}
    &\hspace{.75em}
    \begin{tikzcd}[ampersand replacement=\&]
      f_!\Sigma^{e_f}f^*v_! \ar[swap]{d}{\tr_f \ast v_!}\ar{r}{\sim}
        \& f_!\Sigma^{e_f}u_!g^* \ar[equals]{r}
        \& f_!u_!\Sigma^{e_g}g^* \ar[equals]{r}
        \& v_!g_!\Sigma^{e_g}g^* \ar{d}{v_!\ast\tr_g}
      \\
      v_! \ar[equals]{rrr}
        \& \& \& v_!
    \end{tikzcd}
    \\
    &
    \begin{tikzcd}[ampersand replacement=\&, matrix xscale=0.723]
      v^! \ar[equals]{rrr}\ar[swap]{d}{\cotr_g\ast v^!}
        \& \& \& v^! \ar{d}{v^!\ast\cotr_f}
      \\
      g_*\Sigma^{-e_g}g^!v^! \ar[equals]{r}
        \& g_*\Sigma^{-e_g}u^!f^! \ar[equals]{r}
        \& g_*u^!\Sigma^{-e_f}f^! \ar{r}{\sim}
        \& v^!f_*\Sigma^{-e_f}f^!
    \end{tikzcd}
  \end{align*}
\end{cor}

%!TEX root = ../virtual.tex

\section{Construction of fundamental classes}
\label{sec:construction}

%!TEX root = ../../virtual.tex

\subsection{Euler classes}
\label{sec:Euler}

Before proceeding to our construction of fundamental classes, we begin with a discussion of Euler classes in the setting of bivariant theories.
Our basic definition is very simple and can be formulated unstably.

\begin{num}
Let $X$ be a scheme and $E$ be a vector bundle over $X$.
 Recall that the Thom space $\Th_X(E) \in \SH(X)$ is in fact the suspension spectrum of a pointed motivic space in $\Hot(X)$.
 Moreover, the latter can be modelled by the pointed Nisnevich sheaf of sets
$$
\Th_X(E):=\mathrm{coKer}(E^\times \rightarrow E),
$$
where $E^\times$ is complement of the zero section.

Note that Thom spaces are functorial with respect to monomorphisms of vector bundles.
That is, if $\nu:F \rightarrow E$ is a monomorphism of vector bundles over $X$, one gets a canonical morphism of pointed sheaves:
$$
\nu_*:\Th_X(F) \rightarrow \Th_X(E).
$$
\end{num}

\begin{df}\label{df:unstable_euler}
Let $E$ be a vector bundle over a scheme $X$,
 and $s$ be its zero section. We can regard $s$ as a monomorphism
 of vector bundles $s:X \rightarrow E$.
 We define the \emph{Euler class} $e(E)$ of $E/X$ as the induced map in $\Hot(X)$:
$$
s_*:X_+=\Th_X(X) \rightarrow \Th_X(E).
$$
\end{df}

\begin{rem}\label{rem:unstable_euler}
We will often view the Euler class as a class
  \begin{equation*}
    e(E) \in H(X, \vb{E}) \simeq H(X/X, -\vb{E}),
  \end{equation*}
via the canonical map
  \begin{equation*}
    \Maps_{\Hot(X)}(X_+, \Th_X(E))
      \to \Maps_{\SH(X)}(\un_X, \Th_X(E)).
  \end{equation*}
It can then also be realized as a class in the (twisted) cohomology spectrum of any motivic ring spectrum $\E$ (Definition~\ref{df:Euler_class_coef}).
When $\E$ is oriented, our Euler class coincides with the top Chern class (see Paragraph~\ref{num:oriented_spectra}).
When $\E$ is the Milnor--Witt spectrum (Example~\ref{ex:MW cohomology}) it recovers the classical Euler class in the Chow--Witt group.
\end{rem}

It is easy to see that Euler classes commute with base change:
\begin{lm}
\label{lm:Euler&bc}
For any morphism $f:Y\to X$ and any vector bundle $E$ over $X$, the following diagram is commutative:
  \begin{equation*}
    \begin{tikzcd}[matrix scale=0.7, column sep=60]
      Y_+ \ar{r}{e(f^{-1}(E))}\ar[equals]{d}
        & \Th_Y(f^{-1}(E)) \ar[equals]{d}
      \\
      f^*(X_+) \ar{r}{f^*(e(E))}
        & f^*(\Th_X(E))
    \end{tikzcd}
  \end{equation*}
% $$
% \xymatrix@C=40pt@R=16pt{
% Y_+\ar^-{e(f^{-1}E)}[r]\ar@{=}[d]
%  & \Th(f^{-1}E) \ar^{\wr}[d] \\
% f^*X_+\ar^{f^*e(E)}[r] & f^*\Th(E).
% }
% $$
\end{lm}
\proof
This follows from the fact that the functor $f^*$ commutes with cokernels, and the base change of the zero section of $E$ is the zero section of $f^{-1}(E)$.
\endproof

\begin{num}
 By construction, the Thom space fits into a cofiber sequence in $\Hot(X)$:
$$
(E^\times)_+ \to E_+ \to \Th_X(E).
$$
By $\AA^1$-homotopy invariance, the projection map
 $E \rightarrow X$ induces an isomorphism in $\Hot(X)$,
 whose inverse is induced by the zero section $s:X \rightarrow E$.
 It follows from our construction that the following diagram commutes:
  \begin{equation*}
    \begin{tikzcd}[matrix scale=0.7, column sep=50]
      (E^\times)_+\ar{r}\ar[equals]{d}
        & E_+\ar{r}
        & \Th_X(E)\ar[equals]{d}
      \\
      (E^\times)_+\ar{r}
        & X_+\ar{r}{e(E)}\ar{u}{\wr}[swap]{s_*}
        & \Th_X(E).
    \end{tikzcd}
  \end{equation*}
\end{num}
\begin{df}
For any vector bundle $E$ over $X$, the \emph{Euler cofiber sequence} is the cofiber sequence
$$
(E^\times)_+ \to X_+ \xrightarrow{e(E)} \Th_X(E)
$$
in $\Hot(X)$.
\end{df}

The Euler cofiber sequence immediately yields the following characteristic property of Euler classes:
\begin{prop}\label{prop:euler_zero}
Let $E$ be a vector bundle over $X$.
Any nowhere vanishing section $s$ of $E\to X$ induces a null-homotopy of the morphism $X_+ \to\Th_X(E)$ corresponding to the Euler class $e(E)$.
In particular, $s$ induces an identification $e(E) \simeq 0$ in $H(X, \vb{E})$.
\end{prop}
\begin{proof}
Such a section $s$ induces a section of the morphism $(E^\times)_+ \to X_+$ in the Euler cofiber sequence.
\end{proof}

\begin{cor}\label{cor:Euler obstruction}
Let $E$ be a vector bundle over $X$.
If $E$ contains the trivial line bundle $\AA^1_X$ as a direct summand, then there is an identification $e(E) \simeq 0$ in $H(X, \vb{E})$.
\end{cor}

\begin{num}\label{num:unstable_Thom}
Suppose we have an exact sequence of vector bundles
  \begin{equation}\label{eq:unstable_Thom}
    0 \rightarrow E' \xrightarrow{\nu} E \rightarrow E'' \rightarrow 0
  \end{equation}
over a scheme $X$.
Then the isomorphism $\Th_X(E) \simeq \Th_X(E') \otimes \Th_X(E'')$ in $\SH(X)$ (Paragraph~\ref{num:Thom_virtual}) in fact exists already unstably:
  \begin{equation*}
    \Th_X(E) \simeq \Th_X(E') \wedge \Th_X(E'')
  \end{equation*}
in $\Hot(X)$.
This is clear when the sequence \eqref{eq:unstable_Thom} is split, and one reduces to this case by pulling back to the $\uHom_X(E'', E')$-torsor of splittings of the sequence (which is fully faithful by $\AA^1$-invariance).
\end{num}

\begin{lm}\label{lm:Euler&mono}
Given an exact sequence of vector bundles as in \eqref{eq:unstable_Thom}, the following diagram is commutative:
  \begin{equation*}
    \begin{tikzcd}[matrix scale=0.7, column sep=60]
      \Th_X(E') \ar{r}{\Id\wedge e(E'')}\ar[equals]{d}
        & \Th_X(E') \wedge \Th_X(E'') \ar[equals]{d}
      \\
      \Th_X(E') \ar{r}{\nu_*}
        & \Th_X(E)
    \end{tikzcd}
  \end{equation*}
\end{lm}

\begin{proof}
Argue as above to reduce to the case where \eqref{eq:unstable_Thom} is split.
\end{proof}

The additivity property of Euler classes is then a direct corollary:
\begin{prop}
Given an exact sequence of vector bundles as in \eqref{eq:unstable_Thom}, the  following diagram is commutative:
  \begin{equation*}
    \begin{tikzcd}[matrix scale=0.7, column sep=80]
      X_+ \ar{r}{e(E')\wedge e(E'')}\ar[equals]{d}
        & \Th_X(E') \wedge \Th_X(E'') \ar[equals]{d}
      \\
      X_+ \ar{r}{e(E)}
        & \Th_X(E)
    \end{tikzcd}
  \end{equation*}
\end{prop}
\begin{proof}
Compose the diagram of Lemma~\ref{lm:Euler&mono} with the map $e(E'):X_+ \rightarrow \Th_X(E')$ (on the left).
\end{proof}

%!TEX root = ../../virtual.tex

\subsection{Fundamental classes: regular closed immersions}
\label{sec:closed}

In this section we construct the fundamental class of a regular closed  immersion and demonstrate its expected properties.
Before proceeding, we make a brief digression to consider a certain preliminary construction.

\begin{num}\label{num:G_m splitting}
Let $X$ be a scheme and consider the diagram
 \begin{equation*}
   \begin{tikzcd}[matrix scale=0.7]
     \GG X \ar{rr}{j}\ar[swap]{rd}{\pi}
       & & \AA^1 X \ar{ld}{\bar\pi}
     \\
       & X &
   \end{tikzcd}
 \end{equation*}  
For any $e \in \K(X)$, we have a commutative diagram
  \begin{equation*}
    \begin{tikzcd}[matrix scale=0.7]
      H(X,1-e) \ar[equals]{r}\ar[swap]{d}{\pi^*}
        & H(X/X, e-1) \ar{r}{\bar{\pi}^!}\ar[swap]{d}{\pi^!}
        & H(\AA^1 X/X, e) \ar[swap]{d}{j^!}
      \\
      H(\GG X, 1-e) \ar{r}{\gamma_{\eta_\pi}}
        & H(\GG X/X, e) \ar[equals]{r}
        & H(\GG X/X, e)
    \end{tikzcd}
  \end{equation*}
using the identifications $\vb{T_\pi} \simeq 1$ in $\K(\GG X)$, $\vb{T_{\bar\pi}} \simeq 1$ in $\K(\AA^1 X)$.
(Here, as usual, $1 \in \K(X)$ denotes the unit, i.e., the class $\vb{\AA^1X}$, for any scheme $X$.)
The right-hand square consists of Gysin maps and commutes by Example~\ref{ex:Ayoub&assoc_sm}; moreover, the morphism $\bar{\pi}^!$ is invertible (Lemma~\ref{lm:htp_inv_biv}).
In the left-hand square, the morphism $\gamma_{\eta_\pi}$ is the duality isomorphism (Definition~\ref{df:strong orientation}) associated to the strong orientation (Example~\ref{ex:strong orientation of smooth}) of $\pi$, and the square evidently commutes by construction of the morphisms involved.
Since the left vertical arrow $\pi^*$ admits a retraction $s_1^*$, given by the inverse image by the unit section $s_1 : X \to \GG X$ in cohomology, we also get a canonical retraction $\nu_t$ of the %upper 
right vertical arrow $j^!$.
\end{num}

\begin{num}\label{num:fdl_unit_section}
Consider now the localization triangle
  \begin{equation*}
    H(\AA^1X/X, e)[-1]
      \xrightarrow{j^!} H(\GG X/X, e)[-1]
      \xrightarrow{\partial_{s_0}} H(X/X, e)
  \end{equation*}
associated with the zero section $s_0:X \rightarrow \AA^1_X$.
By Paragraph~\ref{num:G_m splitting} it is canonically split, and we get an $H(X/X, e)$-linear section
  $$ \gamma_t:H(X/X,e) \rightarrow H(\GG X/X,e)[-1]$$
of $\partial_{s_0}$.
By linearity, this map is determined uniquely by a morphism
  $$
  \{t\}:\un_{\GG X} \rightarrow \pi^!(\un_X)[-1]
  $$
in $\SH(\GG X)$; that is, $\gamma_t$ is multiplication by $\{t\} \in H(\GG X/X, 0)[-1]$.
By construction, $\{t\}$ is stable under arbitrary base changes.
If $X$ is an s-scheme over some base $S$, we will abuse notation and write $\gamma_t$ also for the map
$$
\gamma_t:H(X/S,e) \rightarrow H(\GG X/S,e)[-1],
$$
given again by the assignment $x \mapsto \{t\}.x$.
\end{num}

We now proceed to the construction of the fundamental class.

\begin{num}\label{num:deformation}
Let $X$ be an $S$-scheme and $i:Z \rightarrow X$ a regular closed immersion.
We write $D_ZX$ or $D(X,Z)$ for the (affine) deformation space $B_{Z \times 0}(X\times \AA^1)-B_{Z\times 0}(X\times 0)$, as defined by Verdier (denoted $M(Z/X)$ in \cite[\textsection 2]{VerdierRR}); here $B_ZX$ denotes the blow-up of $X$ along $Z$.
This fits into a diagram of tor-independent cartesian squares
\begin{equation*}
  \begin{tikzcd}
    N_ZX \ar{r}{k}\ar{d}
      & D_ZX \ar{d}{r}
      & \GG X \ar[swap]{l}{h}\ar[equals]{d}
    \\
    X \ar{r}{\{0\}}
      & \AA^1 X
      & \GG X \ar{l}
  \end{tikzcd}
\end{equation*}
where $k$ and $h$ are complementary closed/open immersions, and the left-hand arrow is the composite of the projection $p : N_ZX \to Z$ and $i : Z \to X$.
Consider the associated localization triangle (Proposition~\ref{prop:localization}):
  \begin{equation*}
    H(N_ZX/S,e) \xrightarrow{k_*}
    H(D_ZX/S,e) \xrightarrow{h^!}
    H(\GG X/S,e) \xrightarrow{\partial_{N_ZX/D_ZX}}
    H(N_ZX/S,e)[1]
  \end{equation*}
for any $e\in\K(X)$.
\end{num}

\begin{df}\label{df:specialization}
With notation as above, the \emph{specialization to the normal cone} map associated to $i$ is the composite
$$
\sigma_{Z/X}: H(X/S,e) \xrightarrow{\gamma_t} H(\GG X/S,e)[-1]
  \xrightarrow{\partial_{N_ZX/D_ZX}} H(N_ZX/S,e),
$$
where $\gamma_t$ is the map constructed in Paragraph~\ref{num:fdl_unit_section}.
\end{df}

\begin{df}\label{df:fdl_closed}
The \emph{fundamental class} $\fdl_i \in H(Z/X,-\vb{N_ZX})$ associated to the regular closed immersion $i$ is the image of $1 \in H(X/X, 0)$ by the composite
$$
H(X/X,0) \xrightarrow{\sigma_{Z/X}} H(N_ZX/X,0)
  \xrightarrow{\phi_{N_ZX/Z}} H(Z/X,-\vb{N_ZX}),
$$
where $\phi_{N_ZX/Z}$ is the Thom isomorphism of $p : N_ZX\to Z$ (Definition~\ref{df:Thom}).
In other words, $\fdl_i = \phi_{N_ZX/Z}\big(\sigma_{Z/X}(1)\big)$.
\end{df}

\begin{rem}\leavevmode
\begin{remlist}
\item
By definition of the Thom isomorphism, the fundamental class $\fdl_i$ is determined uniquely by the property that there is a canonical identification
  \begin{equation*}
    p^!(\eta_i) \simeq \sigma_{Z/X}(1)
  \end{equation*}
in $H(N_ZX/X,0)$, where $p:N_ZX \rightarrow Z$ is the projection of the normal bundle
 and we have used the notation of Lemma \ref{lm:htp_inv_biv}.

\item
%Our definition of the specialization map is formally very close to the corresponding map in Rost's theory of cycle modules, denoted by $J(X,Z)$ in \cite[\textsection 11]{Ros}.
If one applies the $\mathbb{A}^1$-regulator map (see Definition~\ref{df:regulator}) at the motivic ring spectrum defined by a cycle module, the specialization map is sent to the corresponding map on cycle modules defined by Rost, denoted by $J(X,Z)$ in \cite[\textsection 11]{Ros}. Our construction in the general case therefore uses a very similar idea.

\item
For each regular closed immersion $i$, the fundamental class defines an orientation $(\fdl_i, -\vb{N_i})$ in the sense of Definition~\ref{df:orientation}.
It can be viewed equivalently as a morphism
$$
\fdl_i:\Th_Z(-N_i) \rightarrow i^!(\un_X)
$$
in $\SH(Z)$.

\item
Let us assume that $X$ is an s-scheme over a base $S$.
Then the orientation $(\fdl_i, -\vb{N_i})$ gives rise to a Gysin map (Definition~\ref{df:orientation gives Gysin}):
  $$
  i^!:H(X/S,e) \rightarrow H(Z/S,-\vb{N_i}+e),
  \quad x \mapsto \fdl_i.x.
  $$
It follows from the definitions that this map can also be described as the composite:
$$
H(X/S,e) \xrightarrow{\sigma_{Z/X}} H(N_ZX/S,e)
 \xrightarrow{(p^!)^{-1}} H(Z/S,-\vb{N_ZX}+e),
$$
therefore comparing our construction with that of Verdier \cite{VerdierRR}.
Note also that $\sigma_{Z/X} \simeq p^!i^!$ so that the Gysin map and the specialization map
 uniquely determine each other.
\item One can describe the map $\fdl_i$ more concretely as follows. Let us recall
 the deformation diagram:
  \begin{equation*}
    \begin{tikzcd}[matrix scale=0.7]
      \GG X\ar[hookrightarrow]{r}{h}\ar[swap]{rd}{\pi}
        & D_Z X\ar{d}{r}\ar[hookleftarrow]{r}{k}
        & N_Z X\ar{d}{p}
      \\
        & X \ar[hookleftarrow]{r}{i}
        & Z
    \end{tikzcd}
  \end{equation*}
% $$
% \xymatrix@C=30pt@R=10pt{
% \GG X\ar_\pi[rd]\ar@{^(->}^h[r] & D_ZX\ar^r[d] & N_ZX\ar@{^(->}_k[l]\ar^p[d] \\
% & X & Z\ar^i[l]
% }
% $$
First the map $\{t\}:\un_{\GG X}[1] \rightarrow \pi^!(\un_X)$ corresponds by adjunction and
 after one desuspension to a map:
$$
\sigma_t:\un_{D_Z X} \rightarrow h_*\pi^!(\un_X[-1])=h_*h^*r^!(\un_X[-1]).
$$
Then one gets the following composite map:
$$
\un_{D_Z X} \xrightarrow{\sigma_t} h_*h^*r^!(\un_X[-1])
 \xrightarrow{\text{boundary}} k_!k^!r^!(\un_X)=k_!p^!i^!(\un_X)
 \simeq k_*p ^*\big(\Th_Z(N_ZX) \otimes i^!(\un_X)\big)
$$
where the last isomorphism uses the purity isomorphism $\pur_p$.
Using the identification $\un_D \simeq r^*(\un_X)$ and the adjunction $(r^*,r_*)$
 we deduce a map:
$$
\un_Z \rightarrow p_*p ^*\big(\Th_Z(N_ZX) \otimes i^!(\un_X)\big)
 \simeq \Th_Z(N_ZX) \otimes i^!(\un_X)
$$
where the last isomorphism follows from the $\AA^1$-homotopy invariance 
 of $\SH$ (Paragraph~\ref{num:six/A^1-invariance}).
 The latter composite is nothing else than the morphism $\Th_Z(N_ZX) \otimes \fdl_i$.
\end{remlist}
\end{rem}

\begin{num} \label{not_excess}
Consider a cartesian square where $i$ and $k$ are regular closed immersions
  \begin{equation*}
    \begin{tikzcd}[matrix scale=0.7]
      T \ar{r}{k}\ar[swap]{d}{q}\ar[phantom]{rd}{\scriptstyle\Delta}
        & Y \ar{d}{p}
      \\
      Z \ar[swap]{r}{i}
        & X
    \end{tikzcd}
  \end{equation*}
% $$
% \xymatrix@=10pt{
% T\ar^k[r]\ar_q[d]\ar@{}|\Delta[rd] & Y\ar^p[d] \\
% Z\ar_i[r] & X.
% }
% $$
Then we get a morphism of deformation spaces $D_TY \rightarrow D_ZX$
 and similarly a morphism of vector bundles:
$$
N_TY \xrightarrow{\nu} q^{-1}N_ZX \rightarrow N_ZX
$$
where $\nu$ is in general a monomorphism of vector bundles
 (\emph{i.e.} the codimension of $T$ in $Y$ can be strictly
 smaller than that of $Z$ in $X$: there is excess of intersection).
 We put $\xi=q^{-1}N_ZX/N_TY$, the excess intersection bundle.
\end{num}
\begin{prop}[Excess intersection formula]\label{prop:excess}
With notation as above, we have a canonical isomorphism
$$
\Delta^*(\fdl_i) \simeq e(\xi).\fdl_k
$$
in $H(T/Y,-\vb{q^{-1}N_ZX})$, modulo the identification $\vb{\xi}+\vb{N_TY}\simeq\vb{q^{-1}N_ZX}$,
 where $e(\xi) \in H(T/T,\vb \xi)$ is the Euler class of $\xi$
 (Definition \ref{df:unstable_euler} and Remark \ref{rem:unstable_euler}).
\end{prop}
\begin{proof}
Let us put $D'_TY=D_ZX \times_X Y$ and $N'_TY=q^{-1}N_ZX$.
 Then we get the following commutative diagram of schemes, in which each square is cartesian:
  \begin{equation*}
    \begin{tikzcd}[matrix scale=0.7]
      N_TY \ar[hookrightarrow]{r}\ar{d}
        & D_TY \ar[hookleftarrow]{r}\ar{d}
        & \GG Y \ar[equals]{d}
      \\
      N'_TY \ar[hookrightarrow]{r}\ar{d}
        & D'_T Y \ar[hookleftarrow]{r}\ar{d}
        & \GG Y \ar{d}
      \\
      N_ZX \ar[hookrightarrow]{r}
        & D_ZX \ar[hookleftarrow]{r}
        & \GG X
    \end{tikzcd}
  \end{equation*}
% $$
% \xymatrix@R=12pt@C=24pt{
% N_TY\ar@{^(->}[r]\ar[d] & D_TY\ar[d] & \GG Y\ar@{^(->}[l]\ar@{=}[d] \\
% N'_TY\ar@{^(->}[r]\ar[d] & D'_TY\ar[d] & \GG Y\ar@{^(->}[l]\ar[d] \\
% N_ZX\ar@{^(->}[r] & D_ZX & \GG X\ar@{^(->}[l]
% }
% $$
Therefore, one gets the following commutative diagram:
$$
\xymatrix@R=16pt@C=60pt{
H(\GG Y/Y,0)[-1]\ar^-{\partial_{T/Y}}[r]\ar@{=}[d]\ar@{}|{(1)}[rd]
 & H(N_TY/Y,0)\ar^{\phi_{N_TY/T}}[r]\ar[d]\ar@{}|{(3)}[rd] & H(T/Y,-\vb{N_TY})\ar^{\nu_*}[d]  \\
H(\GG Y/Y,0)[-1]\ar[r]\ar@{}|{(2)}[rd] & H(N'_TY/Y,0)\ar|{\phi_{N'_TY/T}}[r]\ar@{}|{(4)}[rd]
 & H(T/Y,-\vb{N'_TY}) \\
H(\GG X/X,0)[-1]\ar_-{\partial_{Z/X}}[r]\ar^{\Delta^*}[u] & H(N_ZX/Z,0)\ar_{\phi_{N_TY/T}}[r]\ar_{\Delta^*}[u]
 & H(Z/X,-\vb{N_ZX}).\ar_{\Delta^*}[u]
}
$$
Here the right-hand arrow labelled $\Delta^*$ is the change of base map (Paragraph~\ref{num:bivar_ppties}\ref{item:bivar_ppties/base change}); we have abused notation by also writing $\Delta^*$ for the two analogous maps on the left and middle (induced by the obvious cartesian squares).
Square (1) (resp. (2)) is commutative because of the naturality of localization triangles with respect
 to the proper covariance (resp. base change).
 Square (3) is commutative by definition of $\nu_*$,
 and square (4) by compatibility of Thom isomorphisms with respect to base change.

Now observe that the image of $\{t\} \in H(\GG X/X,0)[-1]$ by the counter-clockwise composite in the above diagram is nothing else than the class $\Delta^*(\fdl_i) \in H(T/Y, -\vb{N'_TY})$.
Similarly the image by the clockwise composite is the class $\nu_*(\fdl_k) \in H(T/Y, -\vb{N'_TY})$.
We conclude using Lemma~\ref{lm:Euler&mono}.
\end{proof}

\begin{ex}\label{ex:examples of excess}\leavevmode
 We get the following usual applications of the preceding formula.
\begin{remlist}
\item\label{item:examples of excess/transverse}\emph{Transverse base change formula.}
If we assume that $p$ is transverse to $i$, then $\nu$ is an isomorphism and the excess bundle vanishes.
Thus we get a canonical identification $\Delta^*(\fdl_i) \simeq \fdl_k$.

\item\label{ex:self_inters}\emph{Self-intersection formula.}
If we apply the formula to the self-intersection square
  \begin{equation*}
    \begin{tikzcd}[matrix scale=0.7]
      Z \ar[equals]{r}\ar[equals]{d}\ar[phantom]{rd}{\scriptstyle\Delta}
        & Z \ar{d}{i}
      \\
      Z \ar[swap]{r}{i}
        & X
    \end{tikzcd}
  \end{equation*}
% $$
% \xymatrix@=10pt{
% Z\ar@{=}[r]\ar@{=}[d]\ar@{}|\Delta[rd] & Z\ar^-{i}[d] \\
% Z\ar_-{i}[r] & X
% }
% $$
where the excess bundle equals $N_ZX$, we get a canonical identification
  \begin{equation}\label{eq:top_chern}
    \Delta^*(\fdl_i)\simeq e(N_ZX)
  \end{equation}
in $H(Z/Z,-\vb{N_ZX})=H(Z,\vb{N_ZX})$.

\item\emph{Blow-up formula.}
In the case where $p:Y\to X$ is the blow-up along $Z$, we obtain a generalization of the ``key formula" for blow-ups in \cite[6.7]{Ful}.
\end{remlist}
\end{ex}

\begin{rem}
If $X$ is a scheme, $E$ is a vector bundle over $X$ and $s_0:X\to E$ is the zero section, the self-intersection formula \eqref{eq:top_chern} applied to $s_0$ says that we can recover the Euler class of $E$ from the fundamental class of $s_0$ by base change along the self-intersection square.
That is,
  \begin{equation*}
    e(E) \simeq \Delta^*(\fdl_{s_0}),
  \end{equation*}
where $\Delta$ denotes the self-intersection square associated to $s_0$.
More generally, if $s:X\to E$ is an arbitrary section, consider the cartesian square
  \begin{equation*}
    \begin{tikzcd}[matrix scale=0.7]
      Z_s \ar{r}\ar[swap]{d}\ar[phantom]{rd}{\scriptstyle\Delta_s}
        & X \ar{d}{s_0}
      \\
      X \ar[swap]{r}{s}
        & E
    \end{tikzcd}
  \end{equation*}
% $$
% \xymatrix@=10pt{
% Z_s\ar[r]\ar[d]\ar@{}|{\Delta_s}[rd] & X\ar^-{s_0}[d] \\
% X\ar_-{s}[r] & E.
% }
% $$
We define the \emph{Euler class with support} as
$$
e(E;s):=\Delta_s^*(\fdl_s)\in H(Z_s/X,-\vb{E}).
$$
This notion corresponds to \cite[Definition 5.1]{LevineEnumerative}. In particular the usual Euler class is the case $s=s_0$.
\end{rem}

We will now state good properties of our constructions of orientations
 for regular closed immersions, culminating in the associativity formula.
 Note that all these formulas will be subsumed once we will get our final
 construction.

\begin{num}
First consider a cartesian square of $S$-schemes:
  \begin{equation*}
    \begin{tikzcd}[matrix scale=0.7]
      T \ar[swap]{d}{g}\ar{r}{k}
        & Y \ar{d}{f}
      \\
      Z \ar[swap]{r}{i}
        & X
    \end{tikzcd}
  \end{equation*}
% $$
% \xymatrix@=10pt{
% T\ar^k[r]\ar_g[d] & Y\ar^f[d] \\
% Z\ar^i[r] & X
% }
% $$
such that $i$ is a regular closed immersion and $f$ is smooth.
The isomorphisms of vector bundles
 $T_g \simeq T_f|_T$ and $N_TY\simeq N_ZX|_T$
induce an identification
\begin{equation}\label{eq:lm:assoc1}
\vb{T_g}-\vb{N_ZX|_T} \simeq \vb{L_{T/X}} \simeq \vb{T_f|_T}-\vb{N_TY}
\end{equation}
in $\K(T)$, where $L_{T/X}$ is the cotangent complex of $T$ over $X$.
\end{num}

\begin{lm}\label{lm:Gysin_closed1_special}
With notation as above, one has the commutative square
  \begin{equation*}
    \begin{tikzcd}[matrix scale=0.7,column sep=70]
      H(X/S, *) \ar{r}{\sigma_{Z/X}}\ar{d}{f^!}
        & H(N_ZX/S, *) \ar{d}{N_g(f)^!}
      \\
      H(Y/S, \vb{T_f}+*) \ar{r}{\sigma_{T/Y}}
        & H(N_TY/S,\vb{T_f}+*)
    \end{tikzcd}
  \end{equation*}
modulo the canonical isomorphism $k^*(T_f)|_{N_TY} \simeq N_{N_TY/N_ZX}$ of vector bundles on $N_TY$.
\end{lm}

\begin{proof}
It suffices to show that both squares in the following diagram commute:
$$
\xymatrix@C=16pt@R=20pt{
H(X/S,*)\ar^-{\gamma_t}[r]\ar_{f^!}[d]\ar@{}|{(1)}[rd]
 & H(\GG X/S,*)[-1]\ar^-{\partial_{N_ZX/D_ZX}}[rr]\ar|{(1 \times f)^!}[d]\ar@{}|{(2)}[rrd]
 && H(N_ZX/S,*)\ar|{N_g(f)^!}[d] \\
H(Y/S,\vb{T_f}+*)\ar_-{\gamma_t}[r]
 & H(\GG Y/S,\vb{T_f}+*)[-1]\ar_-{\partial_{N_TY/D_TY}}[rr]
 && H(N_TY/S,\vb{T_f}+*)
}
$$
In fact,
 the commutativity of (1) (where we have denoted the canonical functions of $D_ZX$ and $D_TY$ by the same letter $t$)
 is obvious, and (2) follows from Proposition~\ref{prop:local&Gysin_sm}.
\end{proof}

\begin{lm}\label{lm:Gysin_closed1}
With notation as above, one has a canonical identification
$$
\fdl_g.\fdl_i \simeq \fdl_k.\fdl_f
$$
in $H(T/X,\vb{L_{T/X}})$.
\end{lm}
\begin{proof}
Consider the following diagram:
  \begin{equation*}
    \begin{tikzcd}
      H(X/X, *) \ar{r}{i^!}\ar{d}{f^!}
        & H(Z/X, -\vb{N_ZX}+e) \ar{d}{g^!}\ar{r}{p_{N_ZX/Z}^!}
        & H(N_ZX/X, *) \ar{d}{N_g(f)^!}
      \\
      H(Y/X, \vb{T_f}+*) \ar{r}{k^!}
        & H(T/X,\vb{L_{T/X}}+*) \ar{r}{p_{N_TY/T}^!}
        & H(N_TY/X, \vb{T_f})
    \end{tikzcd}
  \end{equation*}
The right-hand square commutes by the associativity formula for Gysin morphisms associated with smooth morphisms (Example~\ref{ex:Ayoub&assoc_sm}).
Furthermore, the horizontal arrows $p_{N_ZX/Z}^!$ and $p_{N_TY/T}^!$ are invertible (Lemma~\ref{lm:htp_inv_biv}), so it suffices to show that the composite square commutes.
But the upper and lower composites are the respective specialization maps $\sigma_{Z/X}$ and $\sigma_{T/Y}$, so we conclude by Lemma~\ref{lm:Gysin_closed1_special}.
\end{proof}

\begin{num}\label{num:fdl_relative}
Next we consider a commutative diagram of schemes:
  \begin{equation*}
    \begin{tikzcd}[matrix scale=0.7]
      Z \ar[hookrightarrow]{rr}{i}\ar[swap]{rd}{q}
        & & X \ar{ld}{p}
      \\
        & S &
    \end{tikzcd}
  \end{equation*}
% $$
% \xymatrix@R=8pt@C=16pt{
% Z\ar^i[rr]\ar_/-4pt/q[rd] && X\ar^/-4pt/p[ld] \\ & S &
% }
% $$
such that $i$ is a closed immersion and $p$, $q$ are smooth s-morphisms.
In this situation, the canonical exact sequence of vector bundles over $Z$
$$
0 \rightarrow T_q \rightarrow T_p|_Z \rightarrow N_ZX \rightarrow 0
$$
gives rise to an identification
\begin{equation}\label{eq:lm:Gysin_closed_assoc2}
\vb{T_p|_Z}-\vb{N_ZX} \simeq \vb{T_q}.
\end{equation}
in $\K(Z)$.
\end{num}

\begin{lm}\label{lm:Gysin_closed2}
With notation as above, one has a canonical identification
$$
\fdl_p.\fdl_i \simeq \fdl_q
$$
in $H(Z/S,\vb{T_q})$.
\end{lm}
\begin{proof}
 Consider the cartesian square
  \begin{equation*}
    \begin{tikzcd}[matrix scale=0.7]
      D_Z X\ar[hookleftarrow]{r}{k}\ar{d}\ar[phantom]{rd}{\scriptstyle\Delta}
        & N_Z X\ar{d}{\pi}
      \\
      \AA^1_S\ar[swap,hookleftarrow]{r}{s}
        & S
    \end{tikzcd}
  \end{equation*}
% $$
% \xymatrix@=10pt{
% D_ZX\ar[d] & N_ZX\ar^\pi[d]\ar_k[l] \\
% \AA^1_S & S\ar_-s[l]
% }
% $$
where $s$ is the zero section and $\pi$ is the composite map
 $N_ZX \xrightarrow{p_N} Z \xrightarrow q S$.
 The claim will follow from the commutativity of the diagram
$$
\xymatrix@C=20pt@R=22pt{
H(S/S,0)\ar^-{\gamma_t}[r]\ar_{p^!}[d]\ar@{}|{(1)}[rd]
 & H(\GG S/S,0)[-1]\ar^-{\partial_s}[rr]\ar|{(1 \times p)^!}[d]\ar@{}|{(2)}[rrd]
 && H(S/S,0)\ar|{\pi^!}[d]\ar@{}|{(3)}[rd]
 & H(S/S,0)\ar@{=}[l]\ar^{q^!}[d] \\
H(X/S,T_p)\ar_-{\gamma_t}[r]
 & H(\GG X/S,T_p)[-1]\ar_-{\partial_{N_TY/D_TY}}[rr]
 && H(N_ZX/s,T_p)
 & H(Z/s,e),\ar^-{p_{N}^!}[l]
}
$$
by considering the image of $1 \in H(S/S, 0)$ (recall that we have $\partial_s \circ \gamma_t \simeq 1$ by construction, see Paragraph~\ref{num:fdl_unit_section}).

The commutativity of square (1) is obvious, that of (2) follows from Proposition~\ref{prop:local&Gysin_sm} applied to the cartesian square $\Delta$,
 and that of (3) follows from the associativity of Gysin morphisms associated
 with smooth morphisms (Example~\ref{ex:Ayoub&assoc_sm}).
\end{proof}

Before proceeding, we draw out some corollaries of the previous lemma.

\begin{cor}
Consider the assumptions of Paragraph~\ref{num:fdl_relative}.
Then the orientation $\fdl_i$ is universally strong (Definition~\ref{df:strong orientation}).
In other words, the morphism $\fdl_i:\Th_Z(-N_i) \rightarrow i^!(\un_X)$ is invertible.
\end{cor}

\begin{proof}
This follows from the previous lemma and the fact that the maps $\fdl_p$ and $\fdl_q$ are isomorphisms (Definition~\ref{df:fdl_smooth}).
\end{proof}

\begin{cor} \label{cor:section}
Let $p:X \rightarrow S$ be a smooth s-morphism and $i:S \rightarrow X$ a section of $p$.
Then there is a canonical identification $\fdl_i.\fdl_p \simeq 1$ in $H(S/S, 0)$, modulo the identification $i^{-1}(T_p) \simeq N_i$.
\end{cor}

\begin{ex}
Let $X/S$ be an s-scheme.
 Consider a vector bundle $p:E \rightarrow X$, $s_0:X \rightarrow E$ its zero section
 and $e \in \K(X)$.
 Then the associated Gysin map
$$
s_0^!:H(E/S,e) \rightarrow H(X/S,-\vb{E}+e)
$$
is precisely the Thom isomorphism (Definition~\ref{df:Thom}).
This follows from Corollary~\ref{cor:section}.
In cohomological terms, we also get the Thom isomorphism:
  $$
  H(X,e)
    \xrightarrow{p^*} H(E,e)
    \simeq H(E/E,-e)
    \xrightarrow{s_0^!} H(X/E,-\vb{E}-e)
    = H_X(E,e+\vb{E}).
  $$
\end{ex} 

\begin{num}\label{num:double deformation}
We now proceed towards the formulation of the associativity formula for the fundamental classes of two composable regular closed immersions
$$
Z \xrightarrow k Y \xrightarrow i X.
$$
Recall that there is a short exact sequence
$$
0 \rightarrow N_ZY \rightarrow N_ZX \rightarrow N_YX|_Z \rightarrow 0
$$
of vector bundles over $Z$, whence an identification $\vb{N_ZX} \simeq \vb{N_ZY}+\vb{N_YX|_Z}$ in $\K(Z)$ (Paragraph~\ref{num:Thom_virtual}).
There is also a canonical isomorphism of vector bundles
$$
N(N_ZX,N_ZY) \simeq N(N_YX,N_YX|_Z)
$$
over $Z$; we will abuse notation and write $N$ for both.

We will make use of the double deformation space (cf. \cite[\textsection 10]{Ros})
$$
D=D(D_ZX,D_ZX|_Y).
$$
That is, $D$ is the deformation space associated to $D_ZX|_Y \to D_ZX$, which is a regular closed immersion because the cartesian square
\begin{equation*}
  \begin{tikzcd}[matrix scale=0.7]
    D_ZX|_Y \ar{r}\ar{d}
        & D_ZX\ar{d}
      \\
      Y \ar{r}
        & X
  \end{tikzcd}
\end{equation*}
is tor-independent.
Indeed, this can be checked locally on $X$, so we can assume that $Z \to Y \to X$ is a (transverse) base change of $\{0\} \to \AA^m \to \AA^n$ for some $0\le m\le n$.
Since the deformation space is stable under transverse base change (as is the question of tor-independence), we may reduce to the latter situation, in which case $D_ZX\to X$ is just the projection $\AA^n \times\AA^1 \to \AA^n$, which is transverse to any morphism.

Note that $D$ is a scheme over $X\times\AA^2$; we write $s$ and $t$ for the first and second projections to $\AA^1$, respectively.
Set
$$
D_1=D|_{\{0\} \times \AA^1},\quad
D_2=D|_{\AA^1 \times \{0\}},\quad
D_0=D|_{\{0\} \times \{0\}}.
$$
The fibres of $D$ over various subschemes of $\AA^2$ are summarized in the following table.
\begin{center}%
  {\renewcommand{\arraystretch}{1.3}%
    \vspace{-0.6em}
    \begin{tabular}{lcclc}
      $\{0\} \times \AA^1$ & $D(N_ZX,N_ZY)$ & \hspace{3em}
        & $\GG \times \GG$ & $X \times \GG \times \GG$\\
      $\AA^1 \times \{0\}$ & $D(N_YX,N_YX|_Z)$ &
        & $\GG \times \{0\}$ & $D_2-D_0=\GG \times N_ZX$\\
      $\GG \times \AA^1$ & $D-D_1=\GG \times D_ZX$ &
        & $\{0\} \times \GG$ & $D_1-D_0=N_YX \times \GG$\\
      $\AA^1 \times \GG$ & $D-D_2=D_YX \times \GG$ &
        & $\{0\} \times \{0\}$ & $N$
    \end{tabular}
  }
\end{center}
\end{num}

\begin{lm}\label{lm:Gysin_closed3}
Under the assumptions and notation of Paragraph~\ref{num:double deformation}, the diagram
$$
\xymatrix@R=16pt@C=60pt{
H(X/X,e)\ar^-{\sigma_{Y/X}}[r]\ar_{\sigma_{Z/X}}[d]
 & H(N_YX/X,e)\ar^{\sigma_{N_YX|_Z/N_YX}}[d] \\
H(N_ZX/X,e)\ar^-{\sigma_{N_ZY/N_ZX}}[r] & H(N/X,e)
}
$$
commutes for every $e \in \K(X)$.
\end{lm}
\begin{proof}
By construction of the specialization maps, the square in question factors as in the following diagram:
$$
\xymatrix@=26pt@C=42pt{
H(X/X,e)\ar^-{\gamma_s}[r]\ar_{\gamma_t}[d]\ar@{}|{(1)}[rd]
 & H(X \GG^s/X,e)[-1]\ar^{\partial_{N_YX/D_YX}}[rr]\ar^{\partial_{N_ZX/D_ZX}}[d]\ar@{}|{(2)}[rrd]
 && H(N_YX/X,e)\ar^{\gamma_t}[d] \\
H(X\GG^t/X,e)\ar|-{\gamma_s}[r]\ar_{\partial_{N_XZ/D_ZX}}[d]\ar@{}|{(3)}[rd]
 & H(X \GG^s\GG^t/X,e)[-2]\ar|-{\partial_{N_YX\GG^t/D_YX\GG^t}}[rr]\ar|{\partial_{\GG^s N_ZX/\GG^s D_ZX}}[d]\ar@{}|{(4)}[rrd]
 && H(N_YX\GG^t/X,e)[-1]\ar^{\partial_{N/D_1}}[d] \\
H(N_ZX/X,e)\ar_-{\gamma_s}[r]
 & H(\GG^s N_ZX/X,e)[-1]\ar_-{\partial_{N/D_2}}[rr]
 && H(N/X,e)
}
$$
Some remarks on the notation are in order.
First of all we have omitted the symbol $\times$ in the diagram.
We have also used exponents $s$ and $t$ to indicate that $\GG^s$, resp. $\GG^t$, is viewed as a subset of the $s$-axis, resp. $t$-axis, in $\AA^2$.
Finally, we have written $\gamma_u$ for multiplication
 with the class $\sigma_\pi \in H(\GG^u/\ZZ,0)[-1]$ with $u\in\{s,t\}$.

Now observe that squares (2) and (3) commute by Proposition~\ref{prop:multiplication_localization_les}.
Square (1) is \emph{anti}-commutative by Paragraph~\ref{num:bivar_product}.
Applying Corollary~\ref{cor:crossing_localization_les} to the commutative square
  \begin{equation*}
    \begin{tikzcd}[matrix scale=0.7]
      D_0 \ar{r}\ar{d}
        & D_1 \ar{d}
      \\
      D_2 \ar{r}
        & D
    \end{tikzcd}
  \end{equation*}
% $$
% \xymatrix@=10pt{
% D_0\ar[r]\ar[d] & D_1\ar[d] \\
% D_2\ar[r] & D
% }
% $$
we deduce that square (4) is also anti-commutative, whence the claim.
\end{proof}

\begin{thm}\label{thm:fdl_closed}
Let $S$ be a scheme.
Then there exists a system of fundamental classes $(\fdl_i)_i$ (Definition~\ref{df:fdl_classes_abstract}) on the class of regular closed immersions between s-schemes over $S$, satisfying the following properties:
\begin{thmlist}
\item
For every regular closed immersion $i : Z \to X$, the orientation $\fdl_i \in H(Z/X, -\vb{N_ZX})$ is the fundamental class defined in Definition~\ref{df:fdl_closed}.
\item
The system is stable under transverse base change (Definition~\ref{df:fdl_classes_abstract}\ref{item:fdl_classes_abstract/base change}).
\end{thmlist}
\end{thm}

\begin{proof}
According to Definition~\ref{df:fdl_classes_abstract} we must give the following data:
\begin{remlist}
\item\emph{Fundamental classes.}
For any regular closed immersion $i : Z \to X$, we take the orientation $(\fdl_i, e_i)$ with $e_i = -\vb{N_ZX} \in \K(Z)$ and $\fdl_i \in H(Z/X, -\vb{N_ZX})$ as in Definition~\ref{df:fdl_closed}.

\item\emph{Normalisation.}
If $i=\Id_S$ for a scheme $S$, then $N_SS=S$ and the specialization map and Thom isomorphism are both the identity maps on $H(S/S, 0)$, so we have a canonical identification $\fdl_i \simeq 1$.

\item\emph{Associativity formula.}
Given regular closed immersions $i : Y \to X$ and $k : Z \to Y$, the composite $k \circ i$ is again a regular closed immersion and we have an identification $N_ZX \simeq k^*\vb{N_YX}+\vb{N_ZY}$ in $\K(Z)$.
We obtain a canonical identification $\fdl_k . \fdl_i \simeq \fdl_{i\circ k}$ in $H(Z/X, -\vb{N_ZX})$ from the following commutative diagram
$$
\xymatrix@R=30pt@C=50pt{
H(X/X,0)\ar_-{\sigma_{Y/X}}[r]\ar^{\sigma_{Z/X}}[d]\ar@/^12pt/^{i^!}[rr]\ar@/_32pt/_{(ik)^!}[dd]\ar@{}|{(1)}[rd]
 & H(N_YX/X,0)\ar|{\sigma_{N_YX|_Z/N_YX}}[d]\ar@{}|{(2)}[rd]
 & H(Y/X,-\vb{N_YX})\ar^{\sigma_{Z/Y}}[d]\ar^-{p^!_{N_YX/Y}}[l]\ar@/^50pt/^{k^!}[dd] \\
H(N_ZX/X,0)\ar^-{\sigma_{N_ZY/N_ZX}}[r]\ar@{}|{(3)}[rd]
 & H(N/X,0)\ar@{}|{(4)}[rd]
 & H(N_YX|_Z/X,-\vb{N})\ar_-{p_{N/N_YX|_Z}^!}[l] \\
H(Z/X,-\vb{N_ZX})\ar_{p_{N_ZX/Z}^!}[u]\ar_{p_{N_ZY/Z}^!}[r]
 & H(N_ZY/X,-\vb{N})\ar^{p_{N/N_ZY}^!}[u]
& H(Z/X,-\vb{N_ZY-N_YX}|_Z).\ar^-{p_{N_YX|_Z/Z}^!}[u]\ar^-{p_{N_ZY/Z}^!}[l]
}
$$
by evaluating at $1 \in H(X/X, 0)$, since the maps $p_{N_ZY/Z}^!$ and $p_{N/N_ZY}^!$ are invertible (Lemma~\ref{lm:htp_inv_biv}).
Note that each square is indeed commutative:
\begin{enumerate}
\item Apply Lemma~\ref{lm:Gysin_closed3}.
\item Apply Lemma~\ref{lm:Gysin_closed1_special} to the cartesian square
  \begin{equation*}
    \begin{tikzcd}[matrix scale=0.7]
      N_YX|_Z \ar[swap]{d}{p_{N_YX|_Z/Z}}\ar[hookrightarrow]{r}
        & N_YX \ar{d}{p_{N_YX/Y}}
      \\
      Z \ar[hookrightarrow]{r}{k}
        & Y
    \end{tikzcd}
  \end{equation*}
% $$
% \xymatrix@=10pt{
% N_YX|_Z\ar[r]\ar_{p_{N_YX|_Z/Z}}[d] & N_YX\ar^{p_{N_YX/Y}}[d] \\
% Z\ar^k[r] & Y
% }
% $$
\item This square factors into two triangles:
  \begin{equation*}
    \begin{tikzcd}[matrix scale=0.7, column sep=80]
      H(N_ZX/X, 0)\ar{r}{\sigma_{N_ZY/N_ZX}}\ar{rd}{N_Z(i)^!}
        & H(N/X,0)
      \\
      H(Z/X,-\vb{N_ZX})\ar[swap]{r}{p^!_{N_ZY/Z}}\ar{u}{p^!_{N_ZX/Z}}
        & H(N_ZY/X, -\vb{N})\ar[swap]{u}{p^!_{N/N_ZY}}
    \end{tikzcd}
  \end{equation*}
The upper-right triangle commutes by construction of $N_Z(i)^!$, the Gysin map associated to $N_Z(i) : N_ZY \to N_ZX$.
The lower-left triangle commutes by Lemma~\ref{lm:Gysin_closed2} applied to the commutative diagram:
  \begin{equation*}
    \begin{tikzcd}[matrix scale=0.7]
      N_ZY \ar[hookrightarrow]{rr}\ar[swap]{rd}{p_{N_ZY/Z}}
        & & N_ZX \ar{ld}{p_{N_ZX/Z}}
      \\
        & Z &
    \end{tikzcd}
  \end{equation*}
\item Apply the associativity of Gysin morphisms associated with smooth morphisms (Example~\ref{ex:Ayoub&assoc_sm}).
\end{enumerate}

\item\emph{Transverse base change formula.}
Suppose given a tor-independent cartesian square
  \begin{equation*}
    \begin{tikzcd}[matrix scale=0.7]
      T \ar{r}{k}\ar[swap]{d}{q}\ar[phantom]{rd}{\scriptstyle\Delta}
        & Y \ar{d}{p}
      \\
      Z \ar[swap]{r}{i}
        & X
    \end{tikzcd}
  \end{equation*}
where $i$ is a regular closed immersion.
Then $k$ is also a regular closed immersion and there is a canonical identification of vector bundles $N_TY \simeq q^{-1}(N_ZX)$.
The canonical identification $\Delta^*(\fdl_i) \simeq \fdl_k$ comes then from Example~\ref{ex:examples of excess}\ref{item:examples of excess/transverse}.
\end{remlist}
\end{proof}

%!TEX root = ../../virtual.tex

\subsection{Fundamental classes: general case}
\label{sec:lci}

In this subsection we conclude our main construction by gluing the system of fundamental classes defined on the class of smooth morphisms (Definition~\ref{df:fdl_smooth}) together with the system defined on the class of regular closed immersions in the previous subsection (Theorem~\ref{thm:fdl_closed}).

\begin{lm}\label{lm:comp_ci_sm}
Suppose given a commutative diagram
  \begin{equation*}
    \begin{tikzcd}[matrix scale=0.7]
      X \ar[hookrightarrow]{r}{i}\ar[hookrightarrow,swap]{rd}{i'}
        & Y \ar{d}{f}
      \\
        & S
    \end{tikzcd}
  \end{equation*}
where $i$ and $i'$ are regular closed immersions and $f$ is a smooth s-morphism.
Then there is a canonical identification
  \begin{equation*}
    \fdl_i.\fdl_f \simeq \fdl_{i'}
  \end{equation*}
in $H(X/S, -\vb{N_XS})$, modulo the identification $-\vb{N_XS} \simeq -\vb{N_XY}+i^*\vb{T_{Y/S}}$ in $\K(X)$.
\end{lm}
\begin{proof}
The diagram in question factors as follows:
  \begin{equation*}
    \begin{tikzcd}[matrix scale=0.7]
      X \ar[hookrightarrow]{r}{\Gamma_i}\ar[equals]{rd}
        & X \times_S Y \ar[hookrightarrow]{r}{p_2}\ar{d}{p_1}
        & Y \ar{d}{f}
      \\
        & X \ar[hookrightarrow,swap]{r}{i'}
        & S,
    \end{tikzcd}
  \end{equation*}
% $$
%   \xymatrix{
%    & X \ar[r]^-{i'}
%    & S \\
%   X \ar@/_10pt/_{i}[rr] \ar@{=}[ru] \ar[r]^-{\ \ \Gamma_i} & X\times_SY \ar[r]^-{p_2} \ar[u]_-{p_1} & Y \ar[u]^-{f}
%   }
% $$
where $\Gamma_i:X\to X\times_SY$ is the graph of $i$, $p_1$ and $p_2$ are the respective projections, and the square is cartesian.
By Lemma~\ref{lm:Gysin_closed1}, Corollary~\ref{cor:section} and Theorem~\ref{thm:fdl_closed} we get canonical identifications
  $$
  \fdl_i.\fdl_f
    \simeq \fdl_{\Gamma_i}.\fdl_{p_2}.\fdl_f
    \simeq \fdl_{\Gamma_i}.\fdl_{p_1}.\fdl_{i'}
    \simeq \fdl_{i'},
  $$
as claimed.
\end{proof}

We are now ready to state the main theorem, defining a system of fundamental classes on the class of smoothable lci morphisms:

\begin{thm} \label{thm:fdl_qplci}
Let $S$ be a scheme.
Then there exists a system of fundamental classes $(\fdl_f)_f$ (Definition~\ref{df:fdl_classes_abstract}) on the class of smoothable lci s-morphisms between s-schemes over $S$, satisfying the following properties:
  \begin{thmlist}
  \item\label{cond_sm}
  The restriction of the system $(\fdl_f)_f$ to the class of smooth s-morphisms coincides with the system of Example~\ref{ex:Ayoub&assoc_sm}.
  \item\label{cond_closed}
  The restriction of the system $(\fdl_f)_f$ to the class of regular closed immersions coincides with the system of Theorem~\ref{thm:fdl_closed}.
  \item\label{cond_base_change}
  The system is stable under transverse base change (Definition~\ref{df:fdl_classes_abstract}\ref{item:fdl_classes_abstract/base change}).
  \end{thmlist}
\end{thm}

\begin{proof}
To define the system $(\fdl_f)_f$, we must give the following data (see Definition~\ref{df:fdl_classes_abstract}):
  \begin{remlist}
  \item\emph{Fundamental classes.}
  Given a smoothable lci s-morphism $f:X\to S$, we may choose a factorization through a regular closed immersion $i: X \to Y$ and a smooth s-morphism $p : Y \to S$.
  We define the fundamental class $\fdl_f = \fdl_i.\fdl_p \in H(X/S, \vb{L_f})$.
  Note that, given another factorization through some $i': X\to Y'$ and $p': Y \to S$, we obtain a canonical identification $\fdl_i.\fdl_p \simeq \fdl_{i'}.\fdl_{p'}$ by applying Lemma~\ref{lm:comp_ci_sm} to the diagram
  $$
    \xymatrix{
     & Y \ar[r]^-{p}
     & S \\
    X \ar@/_10pt/_{i'}[rr] \ar[ru]^-{i} \ar[r]^-{\ \ \ (i,i')} & Y\times_SY' \ar[r]^-{p_2} \ar[u]_-{p_1} & Y' \ar[u]^-{p'}.
    }
  $$

  \item\emph{Normalisation.}
  If $f=\Id_S$ for a scheme $S$, then we choose the trivial factorization $f = \Id_S \circ \Id_S$ and the normalization properties of Example~\ref{ex:Ayoub&assoc_sm} and Theorem~\ref{thm:fdl_closed} give a canonical identification $\fdl_{f} \simeq 1.1 \simeq 1$.

  \item\emph{Associativity formula.}
  If $f:X\to Y$ and $g:Y\to Z$ are two smoothable lci s-morphisms, consider the commutative diagram
    \begin{equation*}
      \begin{tikzcd}[matrix scale=0.7]
        X \ar[hookrightarrow]{r}{i_1}\ar[swap]{rd}{f}
          & P \ar[hookrightarrow]{r}{i_3}\ar{d}{p_1}
          & R \ar{d}{p_3}
        \\
          & Y \ar[hookrightarrow]{r}{i_2}\ar[swap]{rd}{g}
          & Q \ar{d}{p_2}
        \\
          &
          & Z,
      \end{tikzcd}
    \end{equation*}
  % $$
  %   \xymatrix@=10pt{
  %      X \ar[r]^-{i_1} \ar[rd]_-{f} & P \ar[r]^-{i_3} \ar[d]^-{p_1} & R \ar[d]^-{p_3}\\
  %        &  Y  \ar[r]_-{i_2}  \ar[rd]_-{g} & Q \ar[d]^-{p_2} \\
  %         &  & Z,
  %   }
  % $$
  where the $i_k$'s are closed immersions and the $p_k$'s are smooth morphisms, and the square is cartesian.
  By Example~\ref{ex:Ayoub&assoc_sm}, Lemma~\ref{lm:Gysin_closed1} and Theorem~\ref{thm:fdl_closed} we have identifications
  $$
  \fdl_f.\fdl_g
  \simeq
  \fdl_{i_1}.\fdl_{p_1}.\fdl_{i_2}.\fdl_{p_2}
  \simeq
  \fdl_{i_1}.\fdl_{i_3}.\fdl_{p_3}.\fdl_{p_2}
  \simeq
  \fdl_{i_3\circ i_1}.\fdl_{p_2\circ p_3}
  \simeq
  \fdl_{g\circ f}.
  $$
  \item\emph{Transverse base change formula.}
  Suppose given a tor-independent cartesian square
    \begin{equation*}
      \begin{tikzcd}[matrix scale=0.7]
        Y \ar{r}{g}\ar[swap]{d}{v}\ar[phantom]{rd}{\scriptstyle\Delta}
          & T \ar{d}{u}
        \\
        X \ar[swap]{r}{f}
          & S,
      \end{tikzcd}
    \end{equation*}
  where $f$ and $g$ are smoothable lci s-morphisms.
  Choosing a factorization $f = p\circ i$, where $i$ is a regular closed immersion and $p$ is a smooth s-morphism, there is an induced factorization of the square $\Delta$:
    \begin{equation*}
      \begin{tikzcd}[matrix scale=0.7]
        Y \ar{r}{k}\ar[swap]{d}{v}\ar[phantom]{rd}{\scriptstyle\Delta_i}
          & Y' \ar{r}{q}\ar{d}{v'}\ar[phantom]{rd}{\scriptstyle\Delta_p}
          & T \ar{d}{u}
        \\
        X \ar[swap]{r}{i}
          & X' \ar[swap]{r}{p}
          & S,
      \end{tikzcd}
    \end{equation*}
  where $k$ is a regular closed immersion and $q$ is a smooth s-morphism.
  Now by above and by the transverse base change properties of Example~\ref{ex:Ayoub&assoc_sm} and Theorem~\ref{thm:fdl_closed}, we have identifications
    \begin{equation*}
      \Delta^*(\fdl_f)
        \simeq \Delta^*(\fdl_i.\fdl_p)
        \simeq \Delta_i^*(\fdl_i).\Delta_p^*(\fdl_p)
        \simeq \fdl_{k}.\fdl_{q}
        \simeq \fdl_{g}
    \end{equation*}
  as claimed.
  \end{remlist}
\end{proof}

\begin{num} \label{num:lci_ex_inters}
In fact, the transverse base change property of Theorem~\ref{thm:fdl_qplci}\ref{cond_base_change} is a special case of an excess intersection formula generalizing Proposition~\ref{prop:excess}.
Consider a cartesian square
  \begin{equation*}
    \begin{tikzcd}[matrix scale=0.7]
      Y \ar{d}{v}\ar{r}{g}\ar[phantom]{rd}{\scriptstyle\Delta}
        & T\ar{d}{u}
      \\
      X \ar[swap]{r}{f}
        & S
    \end{tikzcd}
  \end{equation*}
% $$
% \xymatrix@=10pt{
% Y\ar^g[r]\ar_q[d]\ar@{}|{\Delta}[rd] & T\ar^p[d] \\
% X\ar_f[r] & S
% }
% $$
where $f$ and $g$ are smootable lci morphisms.
 Factor $f=p\circ i$ as a closed immersion followed by a smooth morphism and consider the diagram of cartesian squares
  \begin{equation*}
    \begin{tikzcd}[matrix scale=0.7]
      Y \ar{r}{k}\ar[swap]{d}{v}
        & Q \ar{r}{q}\ar[swap]{d}{r}
        & T \ar{d}{u}
      \\
      X \ar[swap]{r}{i}
        & P \ar[swap]{r}{p}
        & S
    \end{tikzcd}
  \end{equation*}
% $$
% \xymatrix@=10pt{
% Y\ar^k[r]\ar_q[d] & Q \ar_r[d] \ar^q[r] & T\ar^p[d] \\
% X\ar^i[r] & P \ar^p[r] & S
% }
% $$
where $k$ and $i$ are regular closed immersions and $q$ and $p$ are smooth. By~\ref{not_excess} there is a canonical monomorphism
 of $Y$-vector bundles $N_YQ \xrightarrow{\nu} v^{-1}N_XP$.
 We let $\xi$ be the quotient bundle.
\end{num}
\begin{prop}\label{prop:lci_ex_inters}
With notation and assumptions as above, there is an identification
$$
\Delta^*(\fdl_f) \simeq e(\xi).\fdl_g
$$
in $H(Y/T,v^*\vb{L_f})$, modulo the identification $v^*\vb{L_f}\simeq-\vb{\xi}+\vb{L_g}$ in $\K(Y)$,
 where $e(\xi) \in H(Y/Y,\vb \xi)$ is the Euler class of $\xi$
 (Definition \ref{df:unstable_euler} and Remark \ref{rem:unstable_euler}).
\end{prop}
This follows from Proposition~\ref{prop:excess} and the fact that fundamental classes for smooth morphisms are compatible with any base change (Example~\ref{ex:Ayoub&assoc_sm}).

%!TEX root = ../virtual.tex

\section{Main results and applications}
\label{sec:applications}

%!TEX root = ../../virtual.tex

\subsection{Fundamental classes and Euler classes with coefficients}

\begin{num}\label{num:naturality in E}
Let $S$ be a scheme and $\E \in \SH(S)$ a motivic spectrum.
Observe that the bivariant spectra $\Ebiv(X/S, v)$ of Definition~\ref{df:bivariant} are natural in $\E$.
That is, given any morphism $\varphi : \E \to \F$ in $\SH(S)$, there is an induced map of spectra
  \begin{equation*}
    \varphi_* : \Ebiv(X/S, v) \to \F(X/S, v)
  \end{equation*}
for every s-scheme $X$ over $S$ and every $v \in \K(X)$.
Note that $\varphi_*$ is compatible with the various functorialities of bivariant theory (Paragraph~\ref{num:bivar_ppties}).
Also, if $\E$ and $\F$ are equipped with multiplications which commute with $\varphi$, then the induced map $\varphi_*$ preserves products (as defined in Paragraph~\ref{num:bivar_ppties}\ref{item:bivar_ppties/product}).
\end{num}

\begin{df}\label{df:regulator}
Let $\E\in\SH(S)$ be a motivic spectrum equipped with a unit map $\eta : \un_S \to \E$.
By Paragraph~\ref{num:naturality in E} there is a canonical natural transformation of bivariant theories
  \begin{equation}\label{eq:regulator}
    \rho_{X/S}:H(X/S,v) \rightarrow \Ebiv(X/S,v)
  \end{equation}
that we call the \emph{$\AA^1$-regulator map}.
\end{df}

\begin{df}\label{df:fdl_class_coef}
Let $\E\in\SH(S)$ be a unital motivic spectrum.
Given a smoothable lci s-morphism $f:X \rightarrow S$, let $\fdl_f \in H(X/S,\vb{L_f})$ denote the fundamental class of $f$ as in Theorem~\ref{thm:fdl_qplci}.
We define the \emph{fundamental class of $f$ with coefficients in $\E$}, denoted
  \begin{equation*}
    \fdl_f^\E \in \Ebiv(X/S,\vb{L_f}),
  \end{equation*}
as the image of $\fdl_f$ by the $\AA^1$-regulator map $\rho_{X/S}$ \eqref{eq:regulator}.
\end{df}

If $\E$ is unital, associative and commutative, then because $\rho_{X/S}$ is compatible with products and ``change of base'' maps, the associativity and base change formulas provided by Theorem~\ref{thm:fdl_qplci} are immediately inherited by the fundamental classes of Definition~\ref{df:fdl_class_coef}.
If we extend Definition~\ref{df:fdl_classes_abstract} as indicated in Remark~\ref{rem:orientation for general E}, then we can state more precisely:

\begin{thm} \label{thm:fdl_lci with coefficients}
Let $\E \in \SH(S)$ be a motivic spectrum equipped with a unital associative commutative multiplication.
Then there exists a system of fundamental classes $(\fdl_f^\E)_f$ on the class of smoothable lci s-morphisms between s-schemes over $S$.
This system is stable under transverse base change, and recovers the system of Theorem~\ref{thm:fdl_qplci} in the case $\E = \un_S$.
\end{thm}

\begin{df}\label{df:Euler_class_coef}
Let $\E\in\SH(S)$ be a unital motivic spectrum.
Then for any scheme $X$ and any vector bundle $E/X$,
 one defines the \emph{Euler class of $E/X$ with coefficients in $\E$}, denoted
  \begin{equation*}
    e(E,\E) \in \E(X,\vb E) \simeq \E(X/X, -\vb{E})
  \end{equation*}
 as the image of the class $e(E) \in H(X,\vb E) \simeq H(X/X, -\vb{E})$ of Definition~\ref{df:unstable_euler} by the
 $\AA^1$-regulator map $\rho_{X/S}$ \eqref{eq:regulator}.
\end{df}

\begin{rem}
It is possible to define fundamental classes with coefficients in arbitrary motivic spectra $\E \in \SH(S)$, without using any multiplicative or even unital structure.
Indeed the constructions of Section~\ref{sec:construction} (where we only considered the case $\E = \un_S$ for simplicity) extend immediately to general $\E$ without any difficulty.
Alternatively we can replace the use of the $\AA^1$-regulator map above by using instead the \emph{module} structure, i.e., the canonical action\footnote{See also \cite[1.2.1]{Deg16}.}
$$
H(Y/X,w) \otimes \E(X/S,v) \rightarrow \E(Y/S,w+q^*v)
$$
which is defined by the same formula used to define products in Paragraph~\ref{num:bivar_ppties}, except that the multiplication map $\mu_\E : \E\otimes\E\to\E$ is replaced by the map $\un_S \otimes \E \rightarrow \E$ encoding the structure of $\un_S$-module on $\E$.
As a third approach, Gysin maps with arbitrary coefficients can be obtained from the purity transformation (Paragraph~\ref{num:Gysin four}).
% With either of these approaches, the discussion in the following subsections on Gysin maps and purity transformations also extends immediately to arbitrary $\E\in\SH(S)$.
\end{rem}

%!TEX root = ../../virtual.tex

\subsection{(Refined) Gysin maps with coefficients}
\label{sec:applications/gyscoeff}

Just as in Definition~\ref{df:orientation gives Gysin} and Proposition~\ref{prop:abstract Gysin}, Theorem~\ref{thm:fdl_lci with coefficients} provides Gysin maps with coefficients:

\begin{thm}\label{thm:Gysin_coef_basics}
Let $\E \in \SH(S)$ be a motivic spectrum equipped with a unital associative commutative multiplication.
Then for any smoothable lci s-morphism $f : X \to Y$ of s-schemes over $S$, and any $v \in \K(Y)$, there is a \emph{Gysin map}
  \begin{align}\label{eq:Gysin_E}
  \begin{split}
  f^!:\Ebiv(Y/S,v) &\to\Ebiv(X/S,\vb{L_f}+f^*v)\\
  x &\mapsto \fdl_f^\E.x.
  \end{split}
  \end{align}
These Gysin maps satisfy functoriality and transverse base change formulas that are exactly analogous to those of Proposition~\ref{prop:abstract Gysin}.
% Moreover, they also satisfy the excess intersection formula as in Proposition~\ref{prop:lci_ex_inters}.
\end{thm}

We also have an excess intersection formula with coefficients (generalizing Proposition~\ref{prop:lci_ex_inters}):

\begin{prop}\label{prop:excess with coefficients}
Let $\E \in \SH(S)$ be a motivic spectrum equipped with a unital associative commutative multiplication.
Suppose given a cartesian square of s-schemes over $S$ of the form
  \begin{equation*}
    \begin{tikzcd}[matrix scale=0.7]
      X' \ar{r}{g}\ar[swap]{d}{u}\ar[phantom]{rd}{\scriptstyle\Delta}
        & Y'\ar{d}{v}
      \\
      X \ar[swap]{r}{f}
        & Y,
    \end{tikzcd}
  \end{equation*}
where $f$ and $g$ are smoothable lci s-morphisms.
Let $\xi$ denote the excess bundle as in Paragraph~\ref{num:lci_ex_inters}.
Then we have a canonical identification
  \begin{equation*}
    \Delta^*(\fdl_f) \simeq e(\xi,\E) . \fdl_g
  \end{equation*}
in $\E(X'/Y', v^*\vb{L_f})$.
If $u$ and $v$ are proper, then we also have an identification
  \begin{equation*}
    f^! \circ v_* \simeq u_* \circ \gamma_{e(\xi,\E)} \circ g^!
  \end{equation*}
of maps $\E(Y'/S, v^*(e)) \to \E(X/S, \vb{L_f}+e)$, for any $e\in\K(Y)$, where $\gamma_{e(\xi,\E)}$ denotes multiplication by $e(\xi,\E)$.
\end{prop}

Applying Proposition~\ref{prop:excess with coefficients} to the cartesian square
  \begin{equation*}
    \begin{tikzcd}[matrix scale=0.7]
      Z \ar[equals]{r}\ar[equals]{d}\ar[phantom]{rd}{\scriptstyle\Delta}
        & Z \ar{d}{i}
      \\
      Z \ar[swap]{r}{i}
        & X,
    \end{tikzcd}
  \end{equation*}
we obtain the self-intersection formula with coefficients:

\begin{cor}\label{cor:self-intersection coefficients}
Let $i : Z \to X$ be a regular closed immersion of s-schemes over $S$.
Then we have canonical identifications
  \begin{equation*}
    \Delta^*(\fdl_i) \simeq e(N_ZX,\E)
  \end{equation*}
in $\Ebiv(Z/Z,-\vb{N_ZX})$, and for any $v\in\K(X)$ and identification
  \begin{equation*}
    i^!i_* \simeq \gamma_{e(N_ZX,\E)}
  \end{equation*}
of maps $\E(Z/S, v) \to \Ebiv(Z/S, -\vb{N_ZX}+i^*v)$.
\end{cor}

Applying Corollary~\ref{cor:self-intersection coefficients} to the zero section $s : X \to E$ of a vector bundle, we obtain the following formula to compute Euler classes in $\E^0(X,\vb E)$:

\begin{cor}
Let $X$ be a scheme and $E$ a vector bundle over $X$ with zero section $s : X \to E$.
Then there is a canonical identification
  \begin{equation}\label{eq:euler&Gysin}
    e(E,\E) \simeq s^!s_*(1)
  \end{equation}
in $\Ebiv(X/X, -\vb{E}) \simeq \E(X, \vb{E})$.
\end{cor}

We now introduce the notions of refined fundamental class and refined Gysin maps, following Fulton's treatment in intersection theory (cf. \cite[6.2]{Ful}).

\begin{df}\label{df:ref_fdl}
Suppose given a cartesian square of s-schemes over $S$ of the form
  \begin{equation}\label{eq:ref_fdl}
    \begin{tikzcd}[matrix scale=0.7]
      X' \ar{r}{g}\ar[swap]{d}{u}\ar[phantom]{rd}{\scriptstyle\Delta}
        & Y'\ar{d}{v}
      \\
      X \ar[swap]{r}{f}
        & Y,
    \end{tikzcd}
  \end{equation}
where $f$ is a smoothable lci s-morphism.
\begin{remlist}
\item
The \emph{refined fundamental class} of $f$, with respect to $\Delta$ and with coefficients in $\E$, is the class
  \begin{equation*}
    \fdl_\Delta^\E = \Delta^*(\fdl_f^\E)
  \end{equation*}
in $\Ebiv(X'/Y',u^*\vb{L_f})$.

\item
The \emph{refined Gysin map} of $f$, with respect to $\Delta$ and with coefficients in $\E$, is the Gysin map associated to the orientation $(\fdl_\Delta^\E, u^*\vb{L_f})$ of $g$, in the sense of Definition~\ref{df:orientation gives Gysin}.
That is, it is the induced map of bivariant spectra
  \begin{align*}
    g^!(\fdl_\Delta^\E): \E(Y'/S, e) &\to \E(X'/S, u^*\vb{L_f}+e)\\
    x &\mapsto \Delta^*(\fdl_f^\E) . x
  \end{align*}
for any $e \in \K(Y')$.
We sometimes denote it also by $g^!_\Delta$.
\end{remlist}
\end{df}

In terms of refined fundamental classes, we can reformulate the transverse base change and excess intersection formulas as follows:

\begin{prop}\label{prop:refined}
Suppose given a cartesian square $\Delta$ as in \eqref{eq:ref_fdl}, with $f$ is a smoothable lci s-morphism.
Then we have:
\begin{thmlist}
\item\label{item:refined/tautological}\emph{Tautological base change formula.}
If $u$ and $v$ are proper, then there is a canonical identification
  \begin{equation*}
    f^! \circ v_* \simeq u_* \circ g^!_\Delta
  \end{equation*}
of maps $\E(Y'/S, u^*(e)) \to \E(X/S, \vb{L_f}+e)$, for any $e\in\K(Y)$.

\item
If $g$ is also smoothable lci s-morphism, then we have an identification
  \begin{equation*}
    \fdl_\Delta^\E \simeq e(\xi,\E) . \fdl_g
  \end{equation*}
in $\E(X'/Y', u^*\vb{L_f})$.
In particular the refined Gysin map $g^!_\Delta$ is identified with the composite $\gamma_{e(\xi,\E)} \circ g^!$, where the map $\gamma_{e(\xi,\E)}$ is multiplication by $e(\xi,\E)$.

\item
If the square $\Delta$ is tor-independent, so that in particular $g$ is also smoothable lci s-morphism, then we have an identification
  \begin{equation*}
    \fdl_\Delta^\E \simeq \fdl_g^\E
  \end{equation*}
in $\E(X'/Y', \vb{L_g})$.
In particular the refined Gysin map $g^!_\Delta$ is identified with the Gysin map $g^!$.
\end{thmlist}
\end{prop}

\begin{proof}
The first claim follows from the definitions (Paragraph~\ref{num:bivar_ppties}).
The second is precisely the excess intersection formula (Proposition~\ref{prop:excess with coefficients}).
The third is the stability of fundamental classes under transverse base change (Theorem~\ref{thm:fdl_lci with coefficients}).
\end{proof}

\begin{rem}
Note that the excess intersection formula (Proposition~\ref{prop:excess with coefficients}) and transverse base change formula can be recovered by combining the tautological base change formula of Proposition~\ref{prop:refined}\ref{item:refined/tautological} with parts (ii) and (iii).
\end{rem}

%!TEX root = ../../virtual.tex

\subsection{Purity, traces and duality}
\label{sec:applications/purity}

\begin{num}\label{num:purity_trans}
Let $S$ be a scheme and let $f : X \to Y$ be a smoothable lci s-morphism of s-schemes over $S$.
Then by Theorem~\ref{thm:fdl_lci with coefficients} we obtain, by the construction of Paragraph~\ref{num:purity_transfo}, a \emph{purity transformation}
  \begin{equation*}
    \pur_f: \Sigma^{L_f}f^{*} \to f^!
  \end{equation*}
of functors $\SH(Y) \to \SH(X)$, as well as trace and cotrace maps (Paragraph~\ref{num:trace_cotrace}):
  \begin{align*}
    \tr_f&:f_!\Sigma^{L_f}f^* \rightarrow \Id_{\SH(Y)} \\
    \cotr_f&:\Id_{\SH(Y)} \rightarrow f_*\Sigma^{-L_f}f^!.
  \end{align*}
These natural transformations satisfy $2$-functoriality and transverse base change properties as described in Proposition~\ref{prop:2-functor} and Corollary~\ref{cor:base change in terms of traces}.
\end{num}

\begin{rem}
There are some instances of purity transformations for non-smooth morphisms in the literature.
One example is the case of Grothendieck duality in the setting of quasi-coherent or ind-coherent sheaves (see e.g. \cite[Theorem~4.3.3]{Conrad}, \cite[Corollary~7.2.4]{GaitsgoryRozenblyumII}), where the transformation is indeed invertible for any Cohen-Macaulay morphism.
A second one, much closer to the motivic context, can be found in \cite[Expos\'e XVIII, (3.2.1.2)]{SGA4} in the derived category of \'etale sheaves. This construction is valid for flat s-morphisms, and only involves Tate twists rather than arbitrary Thom spaces, which reflects the fact that the theory developed in SGA 4 is oriented (cf. Example~\ref{ex:oriented_mot_cat} below).
\end{rem}

\begin{num}\label{num:Gysin four}
The purity transformation induces Gysin maps on each of the four theories defined in Definition~\ref{df:bivariant}.
That is, if $f : X \to Y$ is a smoothable lci s-morphism of s-schemes over $S$, then we get Gysin maps (for every $e \in \K(Y)$):
\begin{remlist}
  \item\emph{Bivariant theory:}
  $f^! : \Ebiv(Y/S, e) \to \Ebiv(X/S, \vb{L_f}+f^*(e))$.

  \item\emph{Cohomology with proper support:}
  $f_! : \Ec(X/S, \vb{L_f}+f^*(e)) \to \Ec(Y/S, e)$.
\end{remlist}
If $f$ is \emph{proper}, then we also get Gysin maps
\begin{remlist}[resume]
  \item\emph{Cohomology:}
  $f_! : \E(X/S, \vb{L_f}+f^*(e)) \to \E(Y/S, e)$.

  \item\emph{Bivariant theory with proper support:}
  $f^! : \Ebivc(Y/S, e) \to \Ebivc(X/S, \vb{L_f}+f^*(e))$.
\end{remlist}
In particular, this is another (obviously equivalent) way to realize the Gysin maps considered in Theorem~\ref{thm:Gysin_coef_basics}.
\end{num}

We now observe that the purity transformation can be defined for any motivic \inftyCat of coefficients:

\begin{num}\label{num:motivic_categories}
Let $\T$ be a motivic \inftyCat of coefficients in the sense of \cite[Chap.~2, Def.~3.5.2]{KhanThesis}, defined on the site $\base$ of (qcqs) schemes.
That is, $\T$ is a presheaf of symmetric monoidal presentable \inftyCats on $\base$ satisfying certain axioms that guarantee (see \cite[Chap.~2, Cor.~4.2.3]{KhanThesis}) that $\T$ admits a full homotopy coherent formalism of six operations.

At this point we note that all the definitions and constructions in Sections~\ref{sec:bivariant} and \ref{sec:construction} make sense in the setting of $\T$ (and not only $\SH$), as they only use the six operations.
In particular:
\begin{enumerate}
  \item 
  One can define the four theories (Definition~\ref{df:bivariant}) in this setting.
  For example, the bivariant theory represented by any $\E \in \T(S)$ is given by:
    \begin{equation*}
      \E(X/S,v,\T) := \Maps_{\T(S)}(\unT_S, p_*(p^!(\E) \otimes \Th_X(-v,\T)))\\
    \end{equation*}
  where $p : X \to S$ is an s-morphism and $v \in \K(X)$.
  Here we have written $\unT_S \in \T(S)$ for the monoidal unit, and $\Th_X(-v,\T)$ for the Thom space\footnote{Thom spaces can be defined in terms of the six operations by relative purity (Paragraph~\ref{num:relative purity}).} internal to $\T$.

  \item
  We have fundamental classes
    \begin{equation*}
      \fdl_f^\T \in \E(X/Y, -\vb{L_f}, \T)
    \end{equation*}
  for any smoothable lci s-morphism $f : X \to Y$ of s-schemes over $S$, with coefficients in any $\E\in\T(S)$ for arbitrary $\T$.
  These again form a system of fundamental classes as in Theorem~\ref{thm:fdl_qplci}, satisfying stability under transverse base change.

  \item
  We have Gysin maps in bivariant theory with coefficients in any $\E\in\T(S)$ (as well as in the other three theories) for arbitrary $\T$.
  These Gysin maps are functorial, satisfy transverse base change and excess intersection formulas.

  \item
  We have natural transformations
    \begin{align*}
      \pur_f^\T &: \Sigma^{L_f}f^{*} \to f^!\\
      \tr_f^\T&:f_!\Sigma^{L_f}f^* \rightarrow \Id \\
      \cotr_f^\T&:\Id \rightarrow f_*\Sigma^{-L_f}f^!
    \end{align*}
  in the setting of any $\T$.
\end{enumerate}
\end{num}

\begin{ex}\label{ex:oriented_mot_cat}
Suppose that $\T$ is \emph{oriented} in the sense that there are Thom isomorphisms
$$
\Th_X(v, \T) \rightarrow \unT_X(r)[2r],
$$
for any $v\in\K(X)$ of virtual rank $r$, which are functorial and respect the $\cE_\infty$-group structure on $\K(X)$ up to a homotopy coherent system of compatibilities.
% together with compatibilities under base change and the map $\K \to \Pic(\T)$. In particular,
Then for any smoothable lci s-morphism $f$ of relative virtual dimension $d$,
 the purity transformation takes the form:
\begin{equation}\label{eq:purity_transfo_oriented}
  \pur_f^\T:f^{*}(-)(d)[2d] \rightarrow f^!
\end{equation}
and similarly for the trace and cotrace maps.
\end{ex}

\begin{ex}
Let $S = \spec{\ZZ[1/\ell]}$ for a prime $\ell$, and let $\Lambda$ be one of $\ZZ/\ell^n\ZZ$, $\ZZ_\ell$, or $\QQ_\ell$.
Then, as $X$ varies over $S$-schemes, the stable \inftyCat of \'etale $\Lambda$-sheaves $D(X_\et,\Lambda)$ defines a motivic \inftyCat of coefficients
 (see \cite{SGA4}, \cite{Ekedahl}, \cite{CD4}, \cite{LiuZheng}).
In particular, we obtain purity transformations
 of the form \eqref{eq:purity_transfo_oriented} generalizing the previously known constructions.
\end{ex}

\begin{df}\label{df:f-pure}
Let $S$ be a scheme and $f : X \to S$ a smoothable lci s-morphism of $S$-schemes.
Let $\T$ be a motivic \inftyCat of coefficients.
We say that $\E \in \T(S)$ is \emph{$f$-pure} if the canonical morphism
  \begin{equation*}
    (\pur_f^\T)_\E : \Sigma^{L_f}f^*(\E) \rightarrow f^!(\E)
  \end{equation*}
is invertible.
\end{df}

\begin{rem}\label{rem:f-pure}\leavevmode
\begin{remlist}
\item If $f$ is \emph{smooth}, then every object $\E \in \T(S)$ is $f$-pure.
This is because the purity theorem (Paragraph~\ref{num:smooth purity}) is valid in any motivic \inftyCat of coefficients $\T$.
\item Variants of Definition~\ref{df:f-pure} have been considered previously (for specific examples of $\T$) by several authors (see e.g. \cite[XVI, 3.1.5]{Gabber}, \cite[4.4.2]{BD1}, \cite[1.7]{PLH}).
\item Given a smoothable lci s-morphism $f : X \to Y$, the full subcategory of $\SH(S)$ spanned by the $f$-pure objects satisfies good formal properties:
 it is stable under direct factors, extensions, and tensor products with strongly dualizable objects.
\end{remlist}
\end{rem}

Note in particular that, for $\T=\SH$, the orientation $\eta_f$ is universally
 strong (Definition~\ref{df:orientation}) if and only if $\un_S$ is $f$-pure.
We have the following variant of Lemma~\ref{lm:universally strong orientation}:

\begin{lm}
Let $\T$ be a motivic \inftyCat of coefficients.
Suppose that $f:X \rightarrow S$ is a smoothable lci s-morphism and that $\E\in\T(S)$ is an $f$-pure object.
Then there are duality isomorphisms
\begin{align*}
\E(X,v)
  & \to \E(X/S,\vb{L_f}-v), \\
\E_c(X/S,v)
  & \to \E^c(X/S,\vb{L_f}-v),
\end{align*}
for every $v\in\K(X)$.
\end{lm}

Recall that an \inftyCat of coefficients $\T$ is called \emph{continuous} if, whenever a scheme $S$ can be written as the limit of a filtered diagram $(S_\alpha)_\alpha$ of schemes with affine dominant transition maps, then the canonical functor
  \begin{equation*}
%    \mathop{\lim_{\longrightarrow}}_{\alpha}%\lim_{\stack{\longrightarrow \\ }} \T(S_\alpha) \to \T(S)
    \colim_\alpha \T(S_\alpha) \to \T(S)
  \end{equation*}
is an equivalence, where the colimit is taken in the \inftyCat of \emph{presentable} \inftyCats (and colimit-preserving functors).

\begin{prop}\label{prop:univ_strong_orientation}
Let $S$ be a scheme, $\T$ a motivic \inftyCat of coefficients, and $\E \in \T(S)$ an object.
Let $f : X \to Y$ be a smoothable s-morphism of $S$-schemes and denote by $p : X \to S$ and $q : Y \to S$ the structure morphisms. Assume that one of the following conditions is satisfied:
\begin{thmlist}
\item
$X$ and $Y$ are smooth over $S$.

\item
$X$ and $Y$ are regular, and $S$ is the spectrum of a field $k$.
The \inftyCat of coefficients $\T$ is continuous.
The object $\E$ is defined\footnote{That is, $\E \in \T(\spec k)$ is isomorphic to $g^*(\E')$, where $g : \spec k \to \spec {k'}$ with $k' \subset k$ a perfect subfield, and $\E' \in \T(\spec{k'})$.} over a perfect subfield of $k$.
\end{thmlist}
Then the morphism $f$ is lci, and $q^*\E$ is $f$-pure.
\end{prop}
\begin{proof}
Since $f$ factors through a closed immersion and a smooth morphism, we may reduce to the case of closed immersions, using the associativity formula and the fact that $\pur_p^\T$ is invertible for $p$ smooth.
Moreover, in both cases $f$ is automatically a \emph{regular} closed immersion.
The second case reduces to the first by using the continuity property of $\T$
 together with Popescu's theorem \cite{Swan}.
For the first case, the morphisms $p : X \to S$ and $q : Y \to S$ are smooth.
By construction we have a commutative diagram
  \begin{equation*}
    \begin{tikzcd}[matrix scale=0.7]
      \E_X \otimes \Th_X(L_f) \ar{r}{\fdl_f^\E}\ar[swap]{d}{\Id\otimes \fdl_f}
        & f^!(\E_Y)
      \\
      \E_X \otimes f^!(\un_Y) \ar[swap]{ru}{Ex^{!*}_{\otimes}}
    \end{tikzcd}
  \end{equation*}
where the left-hand vertical arrow is invertible by Lemma~\ref{lm:Gysin_closed2} and the fact that $\fdl_p$ is an isomorphism for $p$ smooth (Definition~\ref{df:fdl_smooth}).
Therefore it suffices to note that the morphism induced by the exchange transformation $Ex^{!*}_{\otimes}$ (Paragraph~\ref{num:Ex^!*_otimes}) is invertible.
After writing $\E_X = p^*(\E)$ and $\E_Y = q^*(\E)$, and using the purity isomorphisms $\pur_p : p^* \simeq \Sigma^{-T_p}p^!$ and $\pur_q : q^* \simeq \Sigma^{-T_q}q^!$ (Paragraph~\ref{num:smooth purity}), this follows from the the $\otimes$-invertibility of Thom spaces.
\end{proof}

The following definition first appears (as a conjecture)
 in the context of \'etale cohomology in \cite[I, 3.1.4]{SGA5}.
In our setting it was already introduced in \cite{Deg12,Deg16}.
The following could be regarded as a more precise formulation, though in fact it is not difficult to see that both definitions are equivalent (cf. \cite[Prop.~4.2.2]{Deg16}).

\begin{df}\label{df:absolute_purity}
Let $S$ be a scheme, $\T$ be a motivic \inftyCat of coefficients, and $\E \in \T(S)$ an object.
We say that $\E$ satisfies \emph{absolute purity} if the following condition holds: given any commutative triangle
  \begin{equation*}
    \begin{tikzcd}
      X \ar{rr}{f}\ar[swap]{rd}{p}
        & & Y \ar{ld}{q}
      \\
        & S &
    \end{tikzcd}
  \end{equation*}
where $f$, $p$ and $q$ are s-morphisms, $f$ is smoothable lci, and $X$ and $Y$ are regular, then the inverse image $q^*(\E) \in \T(Y)$ is $f$-pure.
% Let $f:X \rightarrow Y$ be a smoothable lci s-morphism between \emph{regular} s-schemes $X$ and $Y$ over $S$.
% Denote by $q : Y \to S$ the structural morphism.
% Then the inverse image 
\end{df}

\begin{rem}\label{rem:absolute_purity}\leavevmode
\begin{remlist}
\item In view of the functoriality property of the purity transformation
 (Proposition~\ref{prop:2-functor}), and by part~(i) of Remark~\ref{rem:f-pure}, it suffices to check the absolute purity property for diagrams as above where $f$ is a \emph{closed immersion}.
\item If $S$ is the spectrum of a field and $\T$ is continuous, then it follows from Proposition~\ref{prop:univ_strong_orientation} that every $\E\in\T(S)$ satisfies absolute purity (see \cite[Appendix C]{DFKJ} for more details). 
\item Given the previous definition, the absolute purity property is stable
 under direct factors,
 extensions, tensor product with strongly dualizable objects
 (as in Remark \ref{rem:f-pure}). One also deduces from the projection formula
 that absolutely pure objects are stable under direct image $p_*$ for $p$ smooth
 and proper.
\end{remlist}
\end{rem}

\begin{ex}
It is known that the motivic spectra $\KGL$, $\HM\QQ$, and $\MGL\otimes \QQ$ satisfy absolute purity over $\spec\ZZ$ (see \cite[Rem.~1.3.5]{Deg12}). This is also the case for the hermitian $K$-theory spectrum $\BO$ in Example~\ref{ex:hermitian K-theory}, see \cite[Appendix D]{DFKJ}; this implies that the rational motivic sphere spectrum $\un\otimes\QQ$ over $\spec{\ZZ[1/2]}$-schemes also satisfies absolute purity, see \cite[Theorem 5.5.1]{Borelcharacter}.
It was conjectured in \cite[Conjectures~B and C]{Deg12} that $\MGL$ and $\un$ 
also satisfy absolute purity.
%In a work in preparation, we will show that 

\end{ex}

\begin{rem}
Note that our definition of absolute purity is formally advantageous.
 For example, recall it was deduced in \cite{CD3} that $\KGL$ satisfies absolute purity.
 This implies that $\KGL_\QQ$ satisfies absolute purity. We know that
 $\HB$ is a direct factor of $\KGL_\QQ$. According to our definition,
 this implies that $\HB$ is absolutely pure. This argument allows to bypass the
 proof of \cite[Th. 14.4.1]{CD3}.
\end{rem}
%!TEX root = ../../virtual.tex

\subsection{Examples}
\label{sec:applications/examples}

We first consider the oriented case, which was already worked out in \cite{Deg12,Deg16,Navarro}.
We will show that we do in fact recover the constructions of \emph{op.cit.} in this case.

\begin{df}\label{df:oriented_spectra}
Let $S$ be a scheme and $\E\in\SH(S)$ a motivic spectrum.
An \emph{orientation} of $\E$ consists of the data of \emph{Thom isomorphisms}
  \begin{equation}\label{eq:Thom_iso}
    \tau^c_v: \Ebiv(X/S,v) \simeq \Ebiv(X/S, r),
  \end{equation}
for every s-scheme $X$ over $S$ and every class $v \in \K(X)$ of virtual rank $r$, which are functorial and respect the $\cE_\infty$-group structure on $\K(X)$ up to a homotopy coherent system of compatibilities.
Here we write simply $r \in \K(X)$ for the class of the trivial bundle of rank $r$.
\end{df}

\begin{ex}
An $\MGL$-module structure on $\E$ (where $\MGL$ is viewed as an $\cE_\infty$-ring spectrum) gives rise to an orientation of $\E$.
\end{ex}

\begin{num}\label{num:oriented_spectra}
Let $\E \in \SH(S)$ be an oriented motivic spectrum.
Then for any smoothable lci s-morphism $f : X \to Y$ of s-schemes over $S$, denote by $d_f$ the virtual dimension of $f$ (i.e., $d_f$ is the rank of $\vb{L_f}\in\K(X)$).
Under the Thom isomorphism \eqref{eq:Thom_iso}, the fundamental class $\fdl_f \in \Ebiv(X/Y,\vb{L_f})$ corresponds to a class
$$
  \fdl^c_f= \tau^c_{L_f}(\fdl_f) \in \Ebiv(X/Y, d_f).
$$
This latter class can be viewed as another orientation $(\fdl^c_f,d_f)$, and coincides with the fundamental class defined in \cite[2.5.3]{Deg16}.\footnote{One reduces
 to the case of regular closed immersion and smooth morphisms. The
 case of smooth morphisms is obvious (which reduces to the six
 functors formalism). For the case of regular closed immersions,
 using the deformation to the normal cone and the compatibility of the two
 fundamental classes with transverse base change, we reduce to the case
 of the zero section of a vector bundle. This case follows because both
 fundamental classes gives the (refined) Thom class of the vector bundle,
 by design.}
These orientations also form a system of fundamental classes as in Definition~\ref{df:fdl_classes_abstract}, and they define Gysin maps
$$
f^!_c = f^!(\fdl^c_f) : \Ebiv(X/S,r) \rightarrow \Ebiv(Y/S,d_f+r),
 \quad
 x \mapsto \fdl^c_f.x,
$$
which are related to the Gysin maps of Theorem~\ref{thm:Gysin_coef_basics} via a commutative diagram:
  $$
  \xymatrix@R=14pt@C=30pt{
  \E(X/S,v)\ar^-{f^!}[r]\ar_{\tau^c_{v}}^\wr[d]
   & \E(Y/S,\vb{L_f}+f^*v)\ar^{\tau^c_{L_f+f^*v}}_\wr[d] \\
  \E(X/S,r)\ar^-{f^!_c}[r] & \E(Y/S,d_f+r).
  }
  $$
Therefore, Theorem~\ref{thm:fdl_lci with coefficients} gives in particular a homotopy coherent refinement of the construction of \cite{Deg16}.
Note, by the way, that the diagram above gives a simple proof of the Grothendieck--Riemann--Roch formula
 (cf. \cite[3.2.6 and 3.3.10]{Deg16} for the formulation in $\AA^1$-homotopy).
 Indeed, it boils down to the definition of the Todd class
 (cf. \cite[3.2.4 and 3.3.5]{Deg16}).
 Compared to \emph{op. cit.}, this proof does not require choosing a factorization of $f$.\footnote{For the record, Grothendieck mentioned that such a direct proof of his formula, without going through a factorisation and the use of a blow-up, should exist.}
\end{num}

\begin{ex}\label{ex:HZ}\emph{Higher Chow groups}.
Suppose that $S$ is the spectrum of a field $k$ of characteristic exponent $p$.
Taking $\E\in\SH(S)$ to be the (oriented) motivic cohomology spectrum $\HM\ZZ$, we may identify the resulting bivariant theory with Bloch's higher Chow groups, up to inverting $p$ \cite[Example~1.2.10(1)]{Deg16}.
The construction of Paragraph~\ref{num:oriented_spectra} then gives Gysin maps in higher Chow groups (with $p$ inverted) for arbitrary smoothable lci s-morphisms.
By the results of Sect.~\ref{sec:applications/gyscoeff}, these Gysin maps are functorial and satisfy transverse base change and excess intersection formulas.
\end{ex}

We now proceed to consider some new examples.

\begin{ex}\label{ex:hermitian K-theory}\emph{Hermitian K-theory.}
Let $S = \spec{\ZZ[1/2]}$.
 According to \cite{PaninWalter}\footnote{%The paper of Panin and
% Walter is not yet published. 
In the case where $S$ is the spectrum of a field of characteristic different from $2$, then one can also take the ring spectrum constructed in \cite{Hornb}.}, for any regular
 $S$-scheme $X$, there exists a motivic ring spectrum $\BO_X\in\SH(X)$
 that represents hermitian K-theory of smooth $X$-schemes.
 In view of its geometric model (denoted by $\BO^{geom}$
 in \emph{op.cit.}), $\BO$ is defined over $S$ (in the sense that there are canonical isomorphisms $\BO_X \simeq f^*(\BO_S)$ for every $f : X \to S$).
 Note that for non-regular schemes, $\BO$-cohomology
 is a \emph{homotopy invariant} version of hermitian K-theory (on the model of \cite{Cis}), though this notion has not yet been introduced and worked out as far as we know.

The twisted bivariant theory associated with $\BO$
 as above is new. The Gysin morphisms that one gets on
 $\BO$-cohomology are also new, at least in the generality
 of arbitrary proper smoothable lci s-morphisms, between arbitrary schemes
 (possibly singular and not defined over a base field).
 In the case of regular schemes, our construction
 for some part of hermitian K-theory (namely, that which compares
 to Balmer's higher Witt groups) should be compared to that
 of \cite{HornbCalmes}.
 This would require a similar discussion to that of Paragraph~\ref{num:oriented_spectra} as, according to Panin and Walter,
 hermitian K-theory has a special kind of orientation
 which allows to consider only twists by line bundles
 (see also the next example).
 We intend to come back to these questions in a future work.
Note that the Gysin morphisms for Balmer-Witt groups agree with the construction in \cite{NenashevWitt} in the case of a closed immersion between smooth quasi-projective schemes over a field of characteristic different from $2$.
\end{ex}

\begin{ex}\label{ex:MW cohomology}\emph{Higher Chow--Witt groups}.
Let $k$ be a perfect field.
Introduced by Barge and Morel,
 the theory of Chow--Witt groups was fully developed by Fasel
 \cite{Fasel1,Fasel2}. More recently, the theory was extended
 to ``higher Chow--Witt groups'' in a series of works
 \cite{CalmesFasel, DF1, DF2}. In particular, given any
 coefficient ring $R$, there exists
 a motivic ring spectrum $\HMW R$ in $\SH(k)$
 called the \emph{Milnor--Witt spectrum} (cf. \cite[3.1.2]{DF2}).
 We denote by $H^{MW}(X/k,v,R)$ (resp. $H_{MW}(X,v,R)$)
 its associated bivariant theory (resp. cohomology).

 For any smooth s-scheme $X$ over $k$ and any $v\in\K(X)$ of virtual rank $r$, one has a canonical isomorphism
  $$
  H_{MW}^0(X,v,R) \simeq \wCH^r(X,\det(v)) \otimes R,
  $$
 which is contravariantly functorial in $X$ and covariantly functorial in
 $v$ \cite[4.2.6, 4.2.7]{DF1}.
 In particular, the ring spectrum $\HMW$ is symplectically oriented
 in the sense of Panin and Walter \cite{PaninWalter}.
 When $X$ is possibly non-smooth,
 the bivariant theory $H^{MW}_0(X/k,v)$ can be computed by a 
 Gersten complex with coefficients in the Milnor--Witt ring of 
 the residue fields, so we can put:
$$
H^{MW}_0(X/k,v,R)=\wCH_r(X,\det(v))\otimes R
$$
and view this as the Chow--Witt group of the scheme
 $X$. Similarly, the groups $H^{MW}_i(X/k,v)$ for $i \geq 0$
 can be viewed as the \emph{higher Chow--Witt groups}. In fact, we
 have canonical maps
$$
\varphi_X:H^{MW}_i(X/k,v,R) \rightarrow \CH_n(X,i) \otimes R,
$$
where $n$ is the rank of the virtual bundle $v$,
 which are functorial in $X$ with respect to proper pushforward
 (resp. pullback along open immersions).

The construction of Theorem~\ref{thm:Gysin_coef_basics} gives Gysin maps on these higher Chow--Witt groups, for any smoothable lci s-morphisms.
These Gysin maps are functorial and satisfy transverse base change and excess intersection formulas.
Furthermore, the maps $\varphi_X$ are compatible with Gysin morphisms by construction.
All in all, we get a robust bivariant theory of higher Chow--Witt groups.
\end{ex}

\begin{ex}\label{ex:A1_homology}\emph{$\AA^1$-homology.}
Let $S$ be a scheme.
Recall that for any commutative ring $R$, there is a motivic ring spectrum $\mathsf NR_S \in \SH(S)$ representing $\AA^1$-homology with coefficients.
This is nothing else than the $R$-linearization $\un_S \otimes R$ of the motivic sphere spectrum (see \cite[5.3.35]{CD3} for another description).
It is clear that $\mathsf{N}R$ is stable under base change in the sense that there are tautological isomorphisms $f^*(\mathsf{N}R_S) \simeq \mathsf{N}R_T$ for every morphism $f : T \to S$.

By Theorem~\ref{thm:Gysin_coef_basics} we obtain Gysin morphisms for the associated bivariant
 theories and cohomologies.
Note in particular that this gives a very general notion
 of transfer maps in cohomology, along arbitary finite lci morphisms, extending the definitions of Morel 
 in \cite{MorelLNM}.\footnote{Morel defines transfer maps
 only for finite field extensions, but he works unstably.}
\end{ex}

%!TEX root = ../../virtual.tex

\subsection{Application: specializations}

In this subsection we investigate two of the many applications of the theory of refined Gysin maps (Definition~\ref{df:ref_fdl}).
Throughout this subsection, we fix a motivic \inftyCat of coefficients $\T$, a scheme $S$, and an object $E\in\T(S)$.

% \begin{num}\label{num:special}
% Let $\T$ be a motivic \inftyCat of coefficients, $S$ a scheme, and $E\in\T(S)$ an object.
% Recall that for any cartesian square of s-schemes over $S$
%   \begin{equation*}
%     \begin{tikzcd}[matrix scale=0.7]
%       X' \ar{r}{g}\ar[swap]{d}{q}\ar[phantom]{rd}{\scriptstyle\Delta}
%         & Y' \ar{d}{p}
%       \\
%       X \ar[swap]{r}{f}
%         & Y
%     \end{tikzcd}
%   \end{equation*}
% with $f$ a smoothable lci s-morphism, we have the \emph{refined Gysin map}
%   \begin{equation*}
%     g^!_\Delta : \E(Y'/S, e) \to \E(X'/S, e+q^*\vb{L_f}).
%   \end{equation*}
% \end{num}

\begin{num}\label{num:specialization 1}
Let $S$ be a scheme.
For an s-scheme $X$ over $S$, any section $s : S \to X$ which is a \emph{regular} closed immersion, and any s-morphism $p : Y \to X$, consider the cartesian square
  \begin{equation*}
    \begin{tikzcd}[matrix scale=0.7]
      Y_s \ar{r}{t}\ar[swap]{d}{p_s}\ar[phantom]{rd}{\scriptstyle\Delta}
        & Y \ar{d}{p}
      \\
      S \ar[swap]{r}{s}
        & X.
    \end{tikzcd}
  \end{equation*}
The associated refined Gysin map takes the form
  \begin{equation*}
    t^!_\Delta: \E(Y/S, e) \to \E(Y_s/S, -p_s^*\vb{N_SX}+e)
  \end{equation*}
for any $e\in\K(Y)$.
Thus any class $\alpha\in\Ebiv(Y/S,e)$ determines a family of \emph{specializations}
  \begin{equation*}
    \alpha_s = t^!_\Delta(\alpha) \in \E(Y_s/S, -p_s^*\vb{N_SX}+e)
  \end{equation*}
for every $s$.
\end{num}

\begin{ex}
In the case where $S=\spec{k}$ is the spectrum of a field, we can take $s : S \to X$ to be any regular $k$-rational point.
\end{ex}

\begin{rem}
If $\E$ is oriented, then we can identify $\E(Y_s/S, e-p_s^*\vb{N_SX}) \simeq \E(Y_s/S, e-d)$, where $d$ is the codimension of the point $s$, and thus view $\alpha_s$ as a class in $\E(Y_s/S, e-d)$.
The same can be accomplished in general up to some choice of a trivialization of the normal bundle $N_SX$.
\end{rem}

\begin{ex}\label{ex:specialization 1 HZ}
In the case where $S=\spec{k}$ is the spectrum of a field and $\E = \HM\ZZ$ (Example~\ref{ex:HZ}), the construction of Paragraph~\ref{num:specialization 1} generalizes Fulton's construction in \cite[Sect.~10.1]{Ful} to higher Chow groups.
\end{ex}

\begin{ex}\label{ex:specialization 1 HZtilde}
Suppose that $S=\spec{k}$ is the spectrum of a perfect field and take now $\E$ to be the Milnor--Witt spectrum (Example~\ref{ex:MW cohomology}).
In this case the construction of Paragraph~\ref{num:specialization 1} gives a refinement of Example~\ref{ex:specialization 1 HZ} with ``coefficients in quadratic forms''.
In particular, we can specialize Chow--Witt cycles: let $X$ be a smooth and connected scheme over $k$ of dimension $d$ and let $p:Y\to X$ be an s-morphism. Then $Y$ can be considered as a family of $k$-schemes parametrized by $X$, and given any Chow--Witt cycle $\alpha \in \wCH_d(Y)$ we get specializations $\alpha_s \in \wCH_0\big(Y_s, p_s^*\det(-N_SX)\big)$.

Beware however that the theory is more complicated than the case of usual Chow groups.
For instance, assume that the morphism $p : Y \to X$ is proper.
Then each $p_s : Y_s \to S$ is proper and we may consider the \emph{degree} of $\alpha_s$ in the (twisted) Grothendieck--Witt group, i.e.:
  \begin{equation*}
    \deg(\alpha_s) = (p_s)_*(\alpha_s) \in \wCH_0\big(S, \det(-N_SX)\big) \simeq \GW(k, \det(-N_SX)) \simeq \GW(k).
  \end{equation*}
Now unlike in the Chow groups, these degrees depend in general on the point $s$.
In fact, the Chow--Witt cycle $p_*(\alpha) \in \wCH_n(X)$ corresponds to the class of an unramified quadratic form $\varphi$ in $\GW(\kappa(X))$, and the class $\deg(\alpha_s)$ is the specialization of $\varphi$ at $s$.
\end{ex}

\begin{num}\label{num:specialization 2}
We now construct an analogue of Fulton's specialization map in \cite[Sect.~20.3]{Ful}. Suppose given cartesian squares
  \begin{equation*}
    \begin{tikzcd}[matrix scale=0.7, column sep=50]
      X_Z \ar[hookrightarrow]{r}{i_X}\ar[swap]{d}{f_Z}\ar[phantom]{rd}{\scriptstyle\Delta}
        & X \ar[hookleftarrow]{r}{j_X}\ar{d}{f}
        & X_U \ar{d}{f_U}
      \\
      Z \ar[hookrightarrow,swap]{r}{i}
        & S \ar[hookleftarrow,swap]{r}{j}
        & U
    \end{tikzcd}
  \end{equation*}
where $i$ is a regular closed immersion of codimension $d$, $j$ is the inclusion of the open complement, and $f$ is an s-morphism. Suppose also given the choice of a null-homotopy $e(N_ZS) \simeq 0$ in $H(Z, \vb{N_ZS})$ (for example when the bundle $N_ZS$ is trivial, any choice of trivialization gives rise to such a null-homotopy by Proposition~\ref{prop:euler_zero}). For any object $A\in\T(X)$, the composition
$$
i_{X*}(i_X^*A\otimes f_Z^{*}Th(-N_ZS))
\xrightarrow{}
i_{X*}i_X^!A
\xrightarrow{}
A
\xrightarrow{}
i_{X*}i_X^*A
$$
where the first map is induced by the refined fundamental class of the square $\Delta$ as in \ref{num:purity_trans} and the other maps are obtained from adjunctions, agrees with the multiplication by the class $f_Z^{*}e(N_ZS)$ by the self-intersection formula (Corollary~\ref{cor:self-intersection coefficients}), which is then null-homotopic by our hypothesis. Therefore using the localization triangle (Proposition~\ref{prop:localization}), we obtain a natural transformation of the form
\begin{equation}
\label{eq:iXtojX}
i_{X*}(i_X^*A\otimes f_Z^{*}Th(-N_ZS))
\to
j_{X!}j_X^!A.
\end{equation}

Now let $\E\in\T(S)$ be an $i$-pure spectrum. Then for any $e \in \K(X)$, the map~\eqref{eq:iXtojX} induces the following \emph{specialization map}:
  \begin{equation}
  \label{eq:specialization_cycles}
    \sigma:\E(X_U/U,e)\simeq\E(X_U/S,e)\to\E(X_Z/S,e-f_Z^{*}\vb{N_ZS})\simeq \E(X_Z/Z,e).
  \end{equation}
%In this situation our goal is to define a \emph{specialization map} of the form
%  \begin{equation}\label{eq:specialization_cycles}
%    \sigma : \E(X_U/U,e) \to \E(X_Z/Z,e).
%  \end{equation}

%To proceed, we will require:
%\begin{enumerate}
%\item\label{item:specialization/pure}
%That the object  (Definition~\ref{df:f-pure}).
%\item\label{item:specialization/euler}
%\end{enumerate}

%We proceed as follows.
%First consider the refined Gysin map
%  \begin{equation*}
%    (i_X)^!_\Delta : \E(X/S,e)\to \E(X_Z/S,e-f_Z^{*}\vb{N_ZS}).
%  \end{equation*}
%From  and assumption~\ref{item:specialization/euler}, we obtain a null-homotopy of the composite
%  \begin{equation*}\label{eq:comp_gysin_top_chern}
%  \E(X_Z/S,e)
%    \xrightarrow{(i_{X})_*}  \E(X/S,e)
%    \xrightarrow{(i_X)^!_\Delta} \E(X_Z/S,e-f_Z^{*}\vb{N_ZS}).
%  \end{equation*}
%, we obtain a factorization of $(i_X)^!_\Delta$ through the spectrum $\simeq\E(X_U/U,e)$.
%We define $\sigma$ as the composite
%  \begin{equation*}
%    \sigma:\E(X_U/U,e)\to\E(X_Z/S,e-f_Z^{*}\vb{N_ZS}) \simeq \E(X_Z/Z,e),
%  \end{equation*}
%where the identification comes from assumption \ref{item:specialization/pure}.
\end{num}

\begin{rem}
In the case where $S=\spec{k}$ is the spectrum of a field and $\E = \HM\ZZ$ (Example~\ref{ex:HZ}), the %construction of Paragraph~\ref{num:specialization 2}
map~\eqref{eq:specialization_cycles} generalizes Fulton's specialization map in \cite[Sect.~20.3]{Ful} to higher Chow groups.
\end{rem}

% \begin{rem}
% Assumption (3) naturally occurs as follows.
%  Let us assume that $S$ is a $k$-scheme for a base field $k$
%  of characteristic exponent $p$.
%  Let $\delta$ be the dimension function given by the Krull dimension.
%  Then we can consider the $\delta$-homotopy t-structure of \cite{BD1}
%  on $\SH[1/p]$. Then according to \cite[3.3.5]{BD1},
%  if $\E$ is non negative with respect to the $\delta$-homotopy t-structure,
%  for any s-scheme $X/S$, one has $\E_p(X/S,e)=0$ as soon as $p<\rk(e)$.
%  In particular we can take $n=\rk(e)$ to get assumption (3).\footnote{This
%  corresponds to the fact specialization map exists on $0$-cycles.}
%  This applies in particular when $\E$ is $\un$, or the ring spectra
%  representing $\AA^1$-homology, motivic cohomology, Milnor--Witt K-theory.

% If $S$ is in addition regular,
%  we get by looking at the relevant degrees of the bivariant
%  theory represented by the Milnor--Witt ring spectrum a well defined
%  specialization map:
% $$
% \wCH_0(X_U,\det(L)) \rightarrow \wCH_0(X_Z,\det(L)),
% $$
% for any line bundle $L$ and any s-scheme $X/S$.
% \end{rem}

\begin{rem}
The construction of Paragraph~\ref{num:specialization 2} is compatible with Ayoub's motivic nearby cycle functor in the following sense.
Let $S$ be the spectrum of a field $k$ of characteristic $0$, $i:S\to\AA^1_S$ the inclusion of the origin, and $j : U=\mathbb{G}_{m,S}\to \AA^1_S$ the complement.
Let $f:X\to\AA^1_S$ be a smooth morphism, $e\in\K(X)$.
Note that any $\E\in\SH(S)$ is $i$-pure by Proposition~\ref{prop:univ_strong_orientation}, and that there is a canonical trivialization of the normal bundle $N_S(\AA^1_S)$.
Therefore the requirements of Paragraph~\ref{num:specialization 2} are satisfied and we obtain a specialization map $\sigma : \E(X_U/U,e) \to \E(X_S/S,e)$ as in~\eqref{eq:specialization_cycles}. Then this map is induced by the canonical natural transformation
  \begin{equation*}
    i_X^*j_{X*}\to \Psi_f
  \end{equation*}
where $\Psi_f:\SH(X_U)\to\SH(X_S)$ is the motivic nearby cycle functor in \cite[3.5.6]{Ayo1}.
\end{rem}
%!TEX root = ../../virtual.tex

\subsection{Application: the motivic Gauss-Bonnet formula}

Let $p:X\to S$ be a smooth proper morphism.
Recall that the spectrum $\Sigma^\infty_+(X) \simeq p_!p^!(\un_S)$ is a strongly dualizable object of $\SH(S)$ \cite[Prop.~2.4.31]{CD3}, so that we may consider the \emph{trace} of its identity endomorphism.
This is an endomorphism $\chi^{cat}(X/S) \in \Maps_{\SH(S)}(\un_S,\un_S)$ that we refer to as the \emph{categorical Euler characteristic}; see \cite[\S~3]{HoyoisLefschetz} for details.
In this subsection, we view $\chi^{cat}(X/S)$ as a class in $H(S/S, 0)$, and compute it as the ``degree of the Euler class of the tangent bundle'':

\begin{thm}\label{thm:Gauss-Bonnet}
Let $p : X \to S$ be a smooth proper morphism.
Then there is an identification $\chi^{cat}(X/S) \simeq p_*(e(T_{X/S}))$ in the group $H(S, 0)$.
\end{thm}

\begin{rem}
Let $S$ be the spectrum of a field of characteristic different from $2$, and let $p : X \to S$ be a smooth projective morphism.
Under these assumptions, a version of Theorem~\ref{thm:Gauss-Bonnet} was proven recently by Levine \cite[Theorem 1]{LevineEnumerative}.
The formulation of \emph{loc. cit.} can be recovered from Theorem~\ref{thm:Gauss-Bonnet} by applying the $\AA^1$-regulator map (Definition~\ref{df:regulator}).
\end{rem}

\begin{num}
In order to prove Theorem~\ref{thm:Gauss-Bonnet}, we begin by giving a useful intermediate description of $\chi^{cat}(X/S)$.
Consider the cartesian square
 \begin{equation*}
   \begin{tikzcd}[matrix scale=0.7]
     X\times_S X \ar{r}{\pi_2}\ar{d}{\pi_1}
       & X \ar{d}{p}
     \\
     X \ar{r}{p}
       & S
   \end{tikzcd}
 \end{equation*}
% $$
% \xymatrix@=14pt{
% X\times_S X\ar^-{\pi_2}[r]\ar_-{\pi_1}[d] & X\ar^p[d] \\
% X\ar^p[r] & S
% }
% $$
and let $\delta : X \to X\times_S X$ denote the diagonal (a regular closed immersion).

\begin{lm}\label{lm:formal description of chi^cat}
The endomorphism $\chi^{cat}(X/S) : \un_S \to \un_S$ is obtained by evaluating the following natural transformation at the monoidal unit $\un_S$:
  \begin{equation*}
    \begin{tikzcd}[matrix scale=0.7]
      \Id \ar[dashed]{rrrr}\ar[swap]{d}{\mathrm{unit}}
        & & & & \Id
      \\
      p_*p^* \ar[equals]{r}
        & p_*\delta^!(\pi_2)^!p^*
        & p_*\delta^!(\pi_1)^*p^! \ar{l}{Ex^{*!}}[swap]{\sim}\ar[swap]{r}{\theta}
        & p_*\delta^*(\pi_1)^*p^! \ar[equals]{r}
        & p_*p^! \ar[swap]{u}{\mathrm{counit}}
    \end{tikzcd}
  \end{equation*}
% $$
% \xymatrix@=14pt{
% \Id\ar@{-->}[rrrr] \ar_-{ad_{(p^*,p_*)}}[d] & & & & \Id \\
% p_*p^* \ar@{=}[r] & p_*\delta^!(\pi_2)^!p^* & p_*\delta^!(\pi_1)^*p^! \ar_-{Ex^{*!}}^-{\sim}[l] \ar^-{\theta}[r] & p_*\delta^*(\pi_1)^*p^! \ar@{=}[r] & p_*p^! \ar_-{ad_{(p_*,p^!)}}[u]
% }
% $$
where $\theta : \delta^! \to \delta^*$ is the exchange transformation $Ex^{*!} : \Id^*\delta^! \to \Id^!\delta^*$.
\end{lm}

\begin{proof}
This follows from the description given in \cite[Prop.~3.6]{HoyoisLefschetz}, in view of the commutativity of the diagram
  \begin{equation*}
    \begin{tikzcd}[matrix scale=0.7]
      \delta^!(\pi_2)^!p^* \ar[equals]{r}
        & p^* \ar[equals]{r}
        & (\pi_1)_*\delta_*\delta^*(\pi_2)^!p^* \ar{r}{\varepsilon}
        & (\pi_1)_*(\pi_2)^!p^*
      \\
      \delta^!(\pi_1)^*p^! \ar{u}{Ex^{*!}}[swap]{\wr}\ar{r}{\theta}
        & p^! \ar[equals]{r}
        & (\pi_1)_*\delta_*\delta^*(\pi_1)^*p^!
        & (\pi_1)_*(\pi_1)^*p^! \ar{u}{Ex^{*!}}[swap]{\wr}\ar{l}{\eta}
    \end{tikzcd}
  \end{equation*}
which the reader will easily verify.
\end{proof}
\end{num}

\begin{proof}[Proof of Theorem~\ref{thm:Gauss-Bonnet}]
By Lemma~\ref{lm:formal description of chi^cat}, it will suffice to show that the following diagram commutes:
  \begin{equation*}
    \begin{tikzcd}[matrix scale=0.7, column sep=60]
      p^* = \delta^!(\pi_2)^!p^* \ar[swap]{rd}{\Sigma^{-T_p}\ast\pur_p}
        & \delta^!(\pi_1)^*p^! \ar[swap]{l}{Ex^{*!}} \ar{r}{\theta}
        & \delta^*(\pi_1)^*p^! = p^!
      \\
        & \Sigma^{-T_p}p^! \ar[swap]{ru}{\mathfrak{e}_p \ast p^!}\ar{u}{\pur_\delta}
        &
    \end{tikzcd}
  \end{equation*}
Here we have written $\mathfrak{e}_p$ for the natural transformation $\Sigma^{-T_p} \to \Id$ induced by the Euler class $e(T_p) : \un_X \to \Th_X(T_p)$.
The commutativity of the left-hand triangle follows by construction of the fundamental class of $p$ (Example~\ref{ex:orientation for smooth}) and Corollary~\ref{cor:section}.
For the right-hand triangle, commutativity follows immediately from the self-intersection formula (Example~\ref{ex:examples of excess}\eqref{ex:self_inters}), which asserts the commutativity of the square
  \begin{equation*}
    \begin{tikzcd}[matrix scale=0.7]
      \Th_X(-T_p) \ar{r}{\fdl_\delta}
        & \delta^!(\pi_1)^*(\un_X) \ar{r}{\theta}
        & \delta^*(\pi_1)^*(\un_X) \ar[equals]{d}
      \\
      \Th_X(-T_p) \ar[equals]{u}\ar{rr}{e(T_p)}
        & & \un_X.
    \end{tikzcd}
  \end{equation*}
\end{proof}

\bibliographystyle{alphamod}
\bibliography{virtual}

\newcommand{\etalchar}[1]{$^{#1}$}
\begin{thebibliography}{EHK{\etalchar{+}}20}

\bibitem[SGA4]{SGA4}
M.~Artin, A.~Grothendieck, and J.-L. Verdier.
\newblock {\em Th\'eorie des topos et cohomologie \'etale des sch\'emas},
  volume 269, 270, 305 of {\em Lecture Notes in Mathematics}.
\newblock Springer-Verlag, 1972--1973.
\newblock S\'eminaire de G\'eom\'etrie Alg\'ebrique du Bois--Marie 1963--64
  (SGA~4).

\bibitem[Ayo07]{Ayo1}
Joseph Ayoub.
\newblock Les six op\'erations de {G}rothendieck et le formalisme des cycles
  \'evanescents dans le monde motivique. {I}.
\newblock {\em Ast\'erisque}, (314):x+466 pp. (2008), 2007.

\bibitem[Bal99]{Balmer}
Paul Balmer.
\newblock Derived {W}itt groups of a scheme.
\newblock {\em J. Pure Appl. Algebra}, 141(2):101--129, 1999.

\bibitem[BD17]{BD1}
M.~Bondarko and F.~D\'eglise.
\newblock Dimensional homotopy t-structures in motivic homotopy theory.
\newblock {\em Adv. Math.}, 311:91--189, 2017.

\bibitem[SGA6]{SGA6}
P.~Berthelot, A.~Grothendieck, and L.~Illusie.
\newblock {\em Th\'eorie des intersections et th\'eor\`eme de
  {R}iemann-{R}och}, volume 225 of {\em Lecture Notes in Mathematics}.
\newblock Springer-Verlag, 1971.
\newblock S\'eminaire de G\'eom\'etrie Alg\'ebrique du Bois--Marie 1966--67
  (SGA~6).

\bibitem[BH17]{BachmannHoyois}
Tom Bachmann and Marc Hoyois.
\newblock Norms in motivic homotopy theory.
\newblock {\em arXiv preprint arXiv:1711.03061}, 2017.

\bibitem[BM60]{BorelMoore}
A.~Borel and J.~C. Moore.
\newblock Homology theory for locally compact spaces.
\newblock {\em Michigan Math. J.}, 7:137--159, 1960.

\bibitem[BM00a]{BarMor}
Jean Barge and Fabien Morel.
\newblock Groupe de {C}how des cycles orient\'{e}s et classe d'{E}uler des
  fibr\'{e}s vectoriels.
\newblock {\em C. R. Acad. Sci. Paris S\'{e}r. I Math.}, 330(4):287--290, 2000.

\bibitem[BM00b]{BargeMorel}
Jean Barge and Fabien Morel.
\newblock Groupe de chow des cycles orient{\'e}s et classe d'euler des
  fibr{\'e}s vectoriels.
\newblock {\em Comptes Rendus de l'Acad{\'e}mie des Sciences-Series
  I-Mathematics}, 330(4):287--290, 2000.

\bibitem[BO74]{BO}
Spencer Bloch and Arthur Ogus.
\newblock Gersten's conjecture and the homology of schemes.
\newblock {\em Ann. Sci. \'{E}cole Norm. Sup. (4)}, 7:181--201 (1975), 1974.

\bibitem[CD16]{CD4}
Denis-Charles Cisinski and Fr{\'e}d{\'e}ric D{\'e}glise.
\newblock \'{E}tale motives.
\newblock {\em Compos. Math.}, 152(3):556--666, 2016.

\bibitem[CD19]{CD3}
Denis-Charles Cisinski and Fr\'{e}d\'{e}ric D\'{e}glise.
\newblock {\em Triangulated categories of mixed motives}.
\newblock Springer Monographs in Mathematics. Springer, Cham, [2019] \copyright
  2019.

\bibitem[CF14]{CalmesFasel}
B.~{Calm{\`e}s} and J.~{Fasel}.
\newblock {Finite Chow-Witt correspondences}.
\newblock arXiv: 1412.2989, December 2014.

\bibitem[CH11]{HornbCalmes}
B.~Calm\`es and J.~Hornbostel.
\newblock Push-forwards for {W}itt groups of schemes.
\newblock {\em Comment. Math. Helv.}, 86(2):437--468, 2011.

\bibitem[Cis13]{Cis}
D.-C. Cisinski.
\newblock Descente par \'eclatements en {$K$}-th\'eorie invariante par
  homotopie.
\newblock {\em Ann. of Math. (2)}, 177(2):425--448, 2013.

\bibitem[Con00]{Conrad}
Brian Conrad.
\newblock {\em Grothendieck duality and base change}, volume 1750 of {\em
  Lecture Notes in Mathematics}.
\newblock Springer-Verlag, Berlin, 2000.

\bibitem[D{\'{e}}g18a]{Deg16}
Fr\'{e}d\'{e}ric D{\'{e}}glise.
\newblock Bivariant theories in motivic stable homotopy.
\newblock {\em Doc. Math.}, 23:997--1076, 2018.

\bibitem[D{\'{e}}g18b]{Deg12}
Fr\'{e}d\'{e}ric D{\'{e}}glise.
\newblock Orientation theory in arithmetic geometry.
\newblock In {\em {$K$}-{T}heory---{P}roceedings of the {I}nternational
  {C}olloquium, {M}umbai, 2016}, pages 239--347. Hindustan Book Agency, New
  Delhi, 2018.

\bibitem[DF17a]{DF1}
F.~{D{\'e}glise} and J.~{Fasel}.
\newblock {MW-motivic complexes}.
\newblock arXiv: 1708.06095, August 2017.

\bibitem[DF17b]{DF2}
F.~{D{\'e}glise} and J.~{Fasel}.
\newblock {The Milnor-Witt motivic ring spectrum and its associated theories}.
\newblock arXiv: 1708.06102, August 2017.

\bibitem[DF19]{Borelcharacter}
F.~D\'eglise and J.~Fasel.
\newblock {The Borel character}.
\newblock arXiv:1903.11679, March 2019.

\bibitem[DFKJ20]{DFKJ}
F.~D\'eglise, J.~Fasel, Adeel~A. Khan, and F.~Jin.
\newblock {On the rational motivic homotopy category}.
\newblock arXiv:2005.10147, May 2020.

\bibitem[DK20]{DruzhininKolderup}
Andrei Druzhinin and H{\aa}kon Kolderup.
\newblock Cohomological correspondence categories.
\newblock {\em Algebr. Geom. Topol.}, 20(3):1487--1541, 2020.

\bibitem[EHK{\etalchar{+}}17]{EHKSY}
Elden Elmanto, Marc Hoyois, Adeel~A. Khan, Vladimir Sosnilo, and Maria
  Yakerson.
\newblock Motivic infinite loop spaces.
\newblock {\em arXiv preprint arXiv:1711.05248}, 2017.

\bibitem[EHK{\etalchar{+}}20]{EHKSY2}
Elden Elmanto, Marc Hoyois, Adeel~A. Khan, Vladimir Sosnilo, and Maria
  Yakerson.
\newblock Framed transfers and motivic fundamental classes.
\newblock {\em Journal of Topology}, 13(2):460--500, 2020.

\bibitem[EK20]{ElmantoKhanPerfection}
Elden Elmanto and Adeel~A. Khan.
\newblock Perfection in motivic homotopy theory.
\newblock {\em Proceedings of the London Mathematical Society}, 120(1):28--38,
  2020.

\bibitem[Eke90]{Ekedahl}
Torsten Ekedahl.
\newblock On the adic formalism.
\newblock In {\em The {G}rothendieck {F}estschrift, {V}ol.\ {II}}, volume~87 of
  {\em Progr. Math.}, pages 197--218. Birkh\"auser Boston, Boston, MA, 1990.

\bibitem[Fas07]{Fasel1}
J.~Fasel.
\newblock The {C}how-{W}itt ring.
\newblock {\em Doc. Math.}, 12:275--312, 2007.

\bibitem[Fas08]{Fasel2}
J.~Fasel.
\newblock Groupes de {C}how-{W}itt.
\newblock {\em M\'em. Soc. Math. Fr. (N.S.)}, (113):viii+197, 2008.

\bibitem[Fas09]{FaselExcess}
J.~Fasel.
\newblock The excess intersection formula for {G}rothendieck--{W}itt groups.
\newblock {\em manuscripta mathematica}, 130(4):411--423, 2009.

\bibitem[FM81]{FMP}
W.~Fulton and R.~MacPherson.
\newblock Categorical framework for the study of singular spaces.
\newblock {\em Mem. Amer. Math. Soc.}, 31(243):vi+165, 1981.

\bibitem[FS09]{FaselSrinivas}
Jean Fasel and Vasudevan Srinivas.
\newblock {C}how-{W}itt groups and {G}rothendieck-{W}itt groups of regular
  schemes.
\newblock {\em Advances in Mathematics}, 221:302--329, 2009.

\bibitem[FS17]{FranklandSpitzweck}
Martin Frankland and Markus Spitzweck.
\newblock Towards the dual motivic {S}teenrod algebra in positive
  characteristic.
\newblock {\em arXiv preprint arXiv:1711.05230}, 2017.

\bibitem[Fuj02]{FujiwaraPurity}
Kazuhiro Fujiwara.
\newblock A proof of the absolute purity conjecture (after {G}abber).
\newblock {\em Algebraic Geometry 2000, Azumino, Hotaka}, pages 153--183, 2002.

\bibitem[Ful98]{Ful}
W.~Fulton.
\newblock {\em Intersection theory}, volume~2 of {\em Ergebnisse der Mathematik
  und ihrer Grenzgebiete. 3. Folge. A Series of Modern Surveys in Mathematics
  [Results in Mathematics and Related Areas. 3rd Series. A Series of Modern
  Surveys in Mathematics]}.
\newblock Springer-Verlag, Berlin, second edition, 1998.

\bibitem[GP18]{GP}
G.~Garkusha and I.~Panin.
\newblock {Framed motives of algebraic varieties (after V. Voevodsky)}.
\newblock arXiv: 1409.4372v4, February 2018.

\bibitem[GR17]{GaitsgoryRozenblyumII}
Dennis Gaitsgory and Nick Rozenblyum.
\newblock {\em A study in derived algebraic geometry. Volume {II}:
  Deformations, {L}ie theory and formal geometry}, volume 221.
\newblock American Mathematical Soc., 2017.

\bibitem[SGA5]{SGA5}
A.~Grothendieck.
\newblock {\em Cohomologie $\ell$-adique et fonctions ${L}$}, volume 589 of
  {\em Lecture Notes in Mathematics}.
\newblock Springer-Verlag, 1977.
\newblock S\'eminaire de G\'eom\'etrie Alg\'ebrique du Bois--Marie 1965--66
  (SGA~5).

\bibitem[Har66]{Hart}
Robin Hartshorne.
\newblock {\em Residues and duality}.
\newblock Lecture notes of a seminar on the work of A. Grothendieck, given at
  Harvard 1963/64. With an appendix by P. Deligne. Lecture Notes in
  Mathematics, No. 20. Springer-Verlag, Berlin-New York, 1966.

\bibitem[Hor05]{Hornb}
J.~Hornbostel.
\newblock {$A^1$}-representability of {H}ermitian {$K$}-theory and {W}itt
  groups.
\newblock {\em Topology}, 44(3):661--687, 2005.

\bibitem[Hoy15]{HoyoisLefschetz}
Marc Hoyois.
\newblock A quadratic refinement of the {G}rothendieck--{L}efschetz--{V}erdier
  trace formula.
\newblock {\em Algebraic \& Geometric Topology}, 14(6):3603--3658, 2015.

\bibitem[Ill06]{IllusieCotangentI}
Luc Illusie.
\newblock {\em Complexe cotangent et d{\'e}formations I}, volume 239.
\newblock Springer, 2006.

\bibitem[ILO14]{Gabber}
L.~Illusie, Y.~Laszlo, and F.~Orgogozo, editors.
\newblock {\em Travaux de {G}abber sur l'uniformisation locale et la
  cohomologie \'etale des sch\'emas quasi-excellents}.
\newblock Soci\'et\'e Math\'ematique de France, Paris, 2014.
\newblock S\'eminaire \`a l'\'Ecole Polytechnique 2006--2008. [Seminar of the
  Polytechnic School 2006--2008], With the collaboration of Fr\'ed\'eric
  D\'eglise, Alban Moreau, Vincent Pilloni, Michel Raynaud, Jo\"el Riou,
  Beno\^\i t Stroh, Michael Temkin and Weizhe Zheng, Ast\'erisque No. 363-364
  (2014) (2014).

\bibitem[{Jin}16]{Jin}
Fangzhou {Jin}.
\newblock {Borel--{M}oore motivic homology and weight structure on mixed
  motives}.
\newblock {\em Math. Z.}, 283(3):1149--1183, 2016.

\bibitem[JY18]{Jin18}
Fangzhou {Jin} and Enlin {Yang}.
\newblock {K\"unneth formulas for motives and additivity of traces}.
\newblock {\em arXiv preprint arXiv:1812.06441}, 2018.

\bibitem[Kha16]{KhanThesis}
Adeel~A. Khan.
\newblock {\em Motivic homotopy theory in derived algebraic geometry}.
\newblock PhD thesis, Universit{\"a}t Duisburg-Essen, 2016.
\newblock Available at \url{https://www.preschema.com/thesis/}.

\bibitem[Kha19]{KhanVirtual}
Adeel~A. Khan.
\newblock Virtual fundamental classes of derived stacks i.
\newblock {\em arXiv preprint arXiv:1909.01332}, 2019.

\bibitem[KW19]{KassWickelgren}
Jesse~Leo Kass and Kirsten Wickelgren.
\newblock The class of {E}isenbud--{K}himshiashvili--{L}evine is the local
  {$\mathbb A^1$}-{B}rouwer degree.
\newblock {\em Duke Math. J.}, 168(3):429--469, 2019.

\bibitem[Lev17a]{LevineVirtual}
Marc Levine.
\newblock The intrinsic stable normal cone.
\newblock {\em arXiv preprint arXiv:1703.03056}, 2017.

\bibitem[Lev17b]{LevineEnumerative}
Marc Levine.
\newblock Toward an enumerative geometry with quadratic forms.
\newblock {\em arXiv preprint arXiv:1703.03049}, 2017.

\bibitem[LR20]{LevineRaksit}
M.~{Levine} and A.~{Raksit}.
\newblock {Motivic {G}auss--{B}onnet formulas}.
\newblock {\em Algebra Number Theory}, 14(7):1801--1851, 2020.

\bibitem[HA]{HA}
Jacob Lurie.
\newblock Higher algebra.
\newblock {\em Preprint, available at
  \url{www.math.harvard.edu/~lurie/papers/HigherAlgebra.pdf}}, 2012.

\bibitem[LZ12]{LiuZheng}
Yifeng Liu and Weizhe Zheng.
\newblock Enhanced six operations and base change theorem for higher {A}rtin
  stacks.
\newblock {\em arXiv preprint arXiv:1211.5948}, 2012.

\bibitem[Mor12]{MorelLNM}
F.~Morel.
\newblock {\em {$\mathbb A^1$}-algebraic topology over a field}, volume 2052 of
  {\em Lecture Notes in Mathematics}.
\newblock Springer, Heidelberg, 2012.

\bibitem[MV99]{MV}
F.~Morel and V.~Voevodsky.
\newblock {${\bf A}^1$}-homotopy theory of schemes.
\newblock {\em Inst. Hautes \'Etudes Sci. Publ. Math.}, (90):45--143 (2001),
  1999.

\bibitem[Nav18]{Navarro}
A.~Navarro.
\newblock Riemann-{R}och for homotopy invariant {$K$}-theory and {G}ysin
  morphisms.
\newblock {\em Adv. Math.}, (328):501–554, 2018.

\bibitem[Nen07]{NenashevWitt}
Alexander Nenashev.
\newblock Gysin maps in balmer-witt theory.
\newblock {\em Journal of Pure and Applied Algebra}, pages 203--221, 2007.

\bibitem[{Pep}15]{PLH}
S.~{Pepin Lehalleur}.
\newblock {Triangulated categories of relative 1-motives}.
\newblock arXiv: 1512.00266, December 2015.

\bibitem[PW18]{PaninWalter2}
I.~{Panin} and C.~{Walter}.
\newblock {On the algebraic cobordism spectra MSL and MSp}.
\newblock arXiv: 1011.0651v2, March 2018.

\bibitem[PW19]{PaninWalter}
I.~{Panin} and C.~{Walter}.
\newblock On the motivic commutative ring spectrum bo.
\newblock {\em St. Petersburg Math. J.}, 30(6):933–972, 2019.

\bibitem[Rio10]{Riou}
J.~Riou.
\newblock Algebraic {$K$}-theory, {${\bf A}^1$}-homotopy and {R}iemann-{R}och
  theorems.
\newblock {\em J. Topol.}, 3(2):229--264, 2010.

\bibitem[Rob14]{RobaloThesis}
Marco Robalo.
\newblock {\em Motivic homotopy theory of non-commutative spaces}.
\newblock PhD thesis, Universit\'e de Montpellier, 2014.
\newblock Available at
  \url{https://webusers.imj-prg.fr/~marco.robalo/these.pdf}.

\bibitem[Rob15]{Robalo}
Marco Robalo.
\newblock {$K$}-theory and the bridge from motives to noncommutative motives.
\newblock {\em Adv. Math.}, 269:399--550, 2015.

\bibitem[Ros96]{Ros}
M.~Rost.
\newblock Chow groups with coefficients.
\newblock {\em Doc. Math. J.}, pages 319--393, 1996.

\bibitem[Swa98]{Swan}
Richard~G. Swan.
\newblock N\'eron-{P}opescu desingularization.
\newblock In {\em Algebra and geometry ({T}aipei, 1995)}, volume~2 of {\em
  Lect. Algebra Geom.}, pages 135--192. Int. Press, Cambridge, MA, 1998.

\bibitem[Tho84]{ThomasonPur}
R.~W. Thomason.
\newblock Absolute cohomological purity.
\newblock {\em Bull. Soc. Math. France}, 112(3):397--406, 1984.

\bibitem[Ver76]{VerdierRR}
Jean-Louis Verdier.
\newblock Expos\'e {IX} : Le th\'eor\`eme de {R}iemann--{R}och pour les
  intersections compl\`etes.
\newblock {\em Ast\'erisque}, 36-37:189--228, 1976.

\bibitem[Voe98]{Voe1}
V.~Voevodsky.
\newblock {$\mathbb A^1$}-homotopy theory.
\newblock In {\em Proceedings of the {I}nternational {C}ongress of
  {M}athematicians, {V}ol. {I} ({B}erlin, 1998)}, number Extra Vol. I, pages
  579--604 (electronic), 1998.

\bibitem[Voe01]{VoevodskyFramed}
Vladimir Voevodsky.
\newblock Notes on framed correspondences.
\newblock unpublished, 2001.

\bibitem[VSF00]{FSV}
V.~Voevodsky, A.~Suslin, and E.~M. Friedlander.
\newblock {\em Cycles, Transfers and Motivic homology theories}, volume 143 of
  {\em Annals of Mathematics Studies}.
\newblock Princeton Univ. Press, 2000.

\end{thebibliography}

\end{document}